\def\bibart#1#2#3#4#5#6#7
{\bibitem{#1}
#2, #4, #5 #6 (#3), #7}

\def\bibcoll#1#2#3#4#5#6#7#8
{\bibitem{#1}
#2, #4, in: #5, #6, #7, #3, pp. #8}

\def\bibbook#1#2#3#4#5#6
{\bibitem{#1} #2, #4, #5, #6, #3}

\def\bibdiss#1#2#3#4#5#6
{\bibitem{#1}
#2 (#3), #4, #5, #6}

\def\mat#1{\mathbb{#1}}

\def\qed{\hfill\ \rule{2mm}{2mm} }

\def\qex{\hfill\ \vbox{\hrule\hbox{\vrule\kern4pt\vbox{\kern4pt{}
\kern4pt}\kern4pt\vrule}\hrule}}

\newcounter{casectr}

\newcounter{reductctr}

\newcounter{claimctr}

\newcounter{classctr}

\documentclass[12pt]{article}

\usepackage{amssymb}

\newtheorem{guess}{Guess}[section]

\newtheorem{define}[guess]{Definition}
\newtheorem{prop}[guess]{Proposition}
\newtheorem{theorem}[guess]{Theorem}

\newtheorem{lem}[guess]{Lemma}

\newtheorem{nota}[guess]{Notation}

\newtheorem{remark}[guess]{Remark}

\newtheorem{oques}[guess]{Open Question}

\newtheorem{conj}[guess]{Conjecture}

\usepackage[pdftitle={Small Width Automorphism Conjecture}, pdfauthor={Bernd Schroeder},
colorlinks=true, pdfstartview=FitV, linkcolor=blue, citecolor=blue, urlcolor=blue]{hyperref}

\usepackage{color}

\begin{document}

\bibliographystyle{plain}

\title{The Automorphism Conjecture for
Ordered Sets of Width $\leq 11$}

\author{
\small Bernd S. W. Schr\"oder \\
\small School of Mathematics and Natural Sciences\\
\small The University of Southern Mississippi\\
\small
118 College Avenue, \#5043\\
\small Hattiesburg, MS 39406\\
}

\date{\small \today}

\maketitle

\begin{abstract}

We
prove the automorphism conjecture
for ordered sets of width less than or equal to $11$.
The proof supports the
meta conjecture that
a large number of
automorphisms
is achievable only as
some type of product of
independent
automorphisms on
highly symmetric subsets.

\end{abstract}


\section{Introduction}

An {\bf ordered set}
consists of an underlying set $P$ equipped with a
reflexive, antisymmetric and transitive relation $\leq $, the order
relation.
An {\bf order-preserving self-map},
or, an {\bf endomorphism},
of an ordered set $P$ is a
self map
$f:P\to P$ such that $p\leq q$ implies
$f(p)\leq f(q)$.
Consistent with standard terminology,
endomorphisms with an inverse that is
an endomorphism, too, are called {\bf automorphisms}.
The set of endomorphisms is denoted ${\rm End} (P)$
and the set of automorphisms is denoted ${\rm Aut} (P)$.
Rival and Rutkowski's
automorphism problem (see \cite{RiRut}, Problem 3)
asks the following.

\begin{oques}

{\bf (Automorphism Problem.)}
Is it true that
$$\lim _{n\to \infty } \max _{|P|=n}
{ |{\rm Aut} (P)|\over |{\rm End} (P)| } =0?$$

\end{oques}

The {\bf Automorphism Conjecture} states that
the
Automorphism Problem
has an affirmative answer.
In light of the facts that, for
almost every ordered set, the identity is the only
automorphism
(see \cite{Proe}, Corollary 2.3a), and that
every ordered set has at least $2^{{2\over3} n} $ endomorphisms
(see \cite{DRSW}, Theorem 1), this conjecture
is quite natural.
Indeed, if, for ordered sets with ``many" automorphisms,
we could show that
there are
``enough" endomorphisms to guarantee the ratio's convergence to zero
(for examples of this technique, see \cite{LRZ,LiuWan}, or
Proposition \ref{manyinmaxlock} here), the conjecture
would be confirmed.
However, the Automorphism Conjecture has been remarkably resilient
against attempts to prove it in general.

Recall that an {\bf antichain} is an ordered set in which no two elements are
comparable and that the {\bf width} $w(P)$
of an ordered set $P$
is the size of the largest antichain contained in $P$.
The Automorphism Conjecture for ordered sets of small width
has recently gathered attention in \cite{BonaMartin}.
It is easy to slightly improve
Theorem 1 in \cite{DRSW}
for ordered sets of bounded width, see Lemma \ref{2tonforbddwidth}.
With such lower bounds available,
it is natural to also consider upper bounds on the number of
automorphisms.
We will see here that the search for
upper bounds on the number of automorphisms
is linked with
numerous insights on
the connection between the combinatorial
structure of an ordered set
and the structure of its automorphism group.

We start our investigation with
ordered sets that have a lot of local symmetry
in Section \ref{maxlocksec}.
Lemma \ref{maxlockeddescr}
essentially shows that, if too much local symmetry
is allowed, then, for any automorphism,
the remainder of the ordered set is locked into
following the automorphism's action on a small subset.
Section \ref{orbclustsec} provides an overall framework
in which this ``transmission of local actions of automorphisms"
can be investigated.
Proposition \ref{getallfromiou}
shows that the actions of automorphisms within
different parts of such a framework,
called interdependent orbit unions,
are independent
of each other.
Section \ref{orbgraphsec} then
provides the framework for the investigation of
``transmission of local actions of automorphisms"
within such interdependent orbit unions.
Theorem \ref{pruneorbit6} is the key to splitting the
automorphism group into two parts, which can then be
analyzed separately.
Bounding the number of automorphisms for the purpose of
tackling ordered sets of small width then
proceeds inductively, using estimates provided in
Section \ref{indprep}. The base step in Section \ref{forbconfsec}
reveals a number of forbidden configurations,
which however, can eventually be defused
with Proposition \ref{manyinmaxlock}.
The induction step in Section \ref{boundnra}
must consider these and other inconvenient
configurations, but it ultimately succeeds in
Theorem \ref{allowedconf}, bounding the number of
automorphisms with a product of factorials of numbers that
are at most one less than the width and whose sum is
at most half of the number of elements minus the height.
(Definition \ref{adequatebound} has the requisite details.)
Section \ref{ACsmwidthproof}
combines all results so far to prove
the Automorphism Conjecture for all ordered sets of width
less than or equal to 11 in
Theorem \ref{ACw<=12}.
Refinements and some further uses of the frameworks developed here
are discussed in
Section \ref{conclusec}.

\section{Max-locked Ordered Sets}
\label{maxlocksec}

An ordered set with ``many" automorphisms must have a high degree of symmetry.
In terms of counting techniques, this means that, even when
the values of an automorphism are known for a ``large" number of points,
the automorphism would still not be uniquely determined by these values.
It is thus natural to try to identify
smaller subsets $S\subset P$
such that
every automorphism
of the ordered set $P$ is uniquely determined
by its values on $S$.

For sets $B,T\subseteq P$, we will write
$B<T$ iff every $b\in B$ is strictly below every $t\in T$.
For singleton sets, we will omit the set braces.
Recall that a nonempty subset $A\subseteq P$ is called
{\bf order-autonomous}
iff, for all $z\in P\setminus A$, we have that
existence of an $a\in A$ with $z<a$ implies $z<A$, and,
existence of an $a\in A$ with $z>a$ implies $z>A$.
An
order-autonomous subset $A\subseteq $ will be called {\bf nontrivial}
iff $|A|\not\in \{ 1,|P|\} $.

Lemma \ref{allbutoneac} below shows the
simple idea that we will iterate throughout this paper:
When there are no nontrivial order-autonomous antichains, then
an automorphism's values
on an antichain are determined by
the
automorphism's
values away from the antichain.
Lemma \ref{allbutoneac} will be generalized in
parts
\ref{pruneorbit7}
and \ref{pruneorbit9}
of Lemma \ref{pruneorbit}.
Recall that
$\uparrow x=\{ p\in P:p\geq x\} $
and
$\downarrow x=\{ p\in P:p\leq x\} $.

\begin{lem}
\label{allbutoneac}

Let $P$ be an ordered set,
let
$A\subseteq P$ be an antichain
that does not contain any nontrivial order-autonomous antichains,
and let
$\Phi , \Psi \in {\rm Aut}
(P)$
be so that both
$\Phi $ and $\Psi $ map $A$ to itself and
$\Phi |_{P\setminus A} =
\Psi |_{P\setminus A} $.
Then $\Phi =\Psi $.

\end{lem}

{\bf Proof.}
Let
$\Phi , \Psi \in {\rm Aut}
(P)$
be so that both
$\Phi $ and $\Psi $ map $A$ to itself and
$\Phi |_{P\setminus A} =
\Psi |_{P\setminus A} $, and let
$x\in A$.
Then
$\uparrow \Phi (x)\setminus \{ \Phi (x)\}
=
\Phi [\uparrow x\setminus \{ x\} ]
=
\Phi |_{P\setminus A} [\uparrow x\setminus \{ x\} ]
=
\Psi |_{P\setminus A} [\uparrow x\setminus \{ x\} ]
=
\uparrow \Psi (x)\setminus \{ \Psi (x)\}
$, and similarly,
$\downarrow \Phi (x)\setminus \{ \Phi (x)\}
=
\downarrow \Psi (x)\setminus \{ \Psi (x)\}
$.
Hence $\{ \Phi (x), \Psi (x)\} \subseteq A$
is an order-autonomous antichain, which, by hypothesis,
implies
$\Phi (x)= \Psi (x)$.
\qed

\vspace{.1in}

Recall that an element $x$ of a finite ordered set $P$ is
said to be {\bf minimal} or of {\bf rank 0}, and we set
${\rm rank} (x):=0$,
iff there is no $z\in P$ such that $z<x$.
Recursively, the element $x$ is said to be of
{\bf rank $k$}, and we set
${\rm rank} (x):=k$,
iff $x$ is minimal in $P\setminus \{ z\in P: {\rm rank} (z)\leq k-1\} $.
It is easy to see that the rank
of a point is preserved by automorphisms.

\begin{define}
\label{inducedonrankdef}

Let $P$ be
an
ordered set.
For every nonnegative integer $j$, we define
$R_j $ to be the set of elements of rank $j$, and we define
$\alpha (R_k ):=|\{ \Phi |_{R_k} :\Phi \in {\rm Aut} (P)\}|$ to be the number of
pairwise distinct restrictions of automorphisms to $R_k $.
Finally, the largest number $h$ such that $R_h \not= \emptyset $
is called the {\bf height} of $P$.

\end{define}

The largest possible value for
$\alpha (R_k )$ in an ordered set of width $w$ is
$w!$.
Because we are interested in bounds for the
number of automorphisms, it is natural to consider
ordered sets for which
$\alpha (R_k )$ is close to $w!$, where ``close" will turn out
to mean $\geq (w-1)!$, except in the case of width $w=4$,
in which it will mean $>8$.
Lemma \ref{maxlockeddescr} will show that, for ordered sets with
large values $\alpha (R_k )$,
every automorphism is determined by its values on $R_k $.
Note that, for consideration of upper bounds, it is
also natural to exclude points that are
common fixed points for all automorphisms.

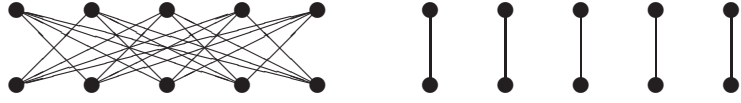
\begin{figure}

\centerline{
\unitlength 1mm 
\linethickness{0.4pt}
\ifx\plotpoint\undefined\newsavebox{\plotpoint}\fi 
\begin{picture}(101,16)(0,0)
\put(45,5){\circle*{2}}
\put(45,15){\circle*{2}}
\put(35,5){\circle*{2}}
\put(35,15){\circle*{2}}
\put(25,5){\circle*{2}}
\put(25,15){\circle*{2}}
\put(70,5){\circle*{2}}
\put(90,5){\circle*{2}}
\put(15,5){\circle*{2}}
\put(70,15){\circle*{2}}
\put(90,15){\circle*{2}}
\put(15,15){\circle*{2}}
\put(80,5){\circle*{2}}
\put(100,5){\circle*{2}}
\put(60,5){\circle*{2}}
\put(5,5){\circle*{2}}
\put(80,15){\circle*{2}}
\put(100,15){\circle*{2}}
\put(60,15){\circle*{2}}
\put(5,15){\circle*{2}}
\put(70,5){\line(0,1){10}}
\put(90,5){\line(0,1){10}}
\put(80,5){\line(0,1){10}}
\put(100,5){\line(0,1){10}}
\put(60,5){\line(0,1){10}}
\put(25,15){\line(-2,-1){20}}
\put(35,15){\line(-2,-1){20}}
\put(45,15){\line(-2,-1){20}}
\put(25,5){\line(-2,1){20}}
\put(35,5){\line(-2,1){20}}
\put(45,5){\line(-2,1){20}}
\put(5,5){\line(1,1){10}}
\put(15,5){\line(1,1){10}}
\put(25,5){\line(1,1){10}}
\put(35,5){\line(1,1){10}}
\put(15,5){\line(-1,1){10}}
\put(25,5){\line(-1,1){10}}
\put(35,5){\line(-1,1){10}}
\put(45,5){\line(-1,1){10}}
\put(45,15){\line(-4,-1){40}}
\put(45,5){\line(-4,1){40}}
\put(5,5){\line(3,1){30}}
\put(15,5){\line(3,1){30}}
\put(5,15){\line(3,-1){30}}
\put(15,15){\line(3,-1){30}}
\end{picture}
}

\caption{The standard example $S_5 $ and the disjoint union $5C_2 $ of
five $2$-chains.
}
\label{st_and_2C}

\end{figure}

\begin{define}

(See Figure \ref{st_and_2C}.)
For $k\geq 3$, we define $S_k $ to be the
standard example of a $k$-dimensional ordered set, that is,
$S_k $ has $k$ minimal elements $\ell _1 ,, \ldots , \ell _k $,
$k$ maximal elements $u_1 , \ldots , u_k $ and
$\ell _i <u_j $ unless $i=j$.
For $k\geq 1$, we define $kC_2 $ to be the disjoint union of
$k$ chains with $2$ elements each.

\end{define}

\begin{lem}
\label{w-1factlem}

Let $P$ be an
ordered set
of width $w\geq 3$
such that no $p\in P$ is a common fixed point for all
automorphisms of $P$,
such that
there is a $j$ such that $|R_j |=w-1$ and
$\alpha (R_j )= (w-1)!$,
and such that $
\{ x\in P: {\rm rank} (x)\geq j+1\} \not= \emptyset $
is not order-autonomous.
Then
$R_j \cup R_{j+1} $ is isomorphic to a set
$S_{w-1}$ or $(w-1)C_2 $,
and consequently
$\alpha (R_{j+1} )= (w-1)!$.

\end{lem}

{\bf Proof.}
First
suppose, for a contradiction,
that there is an $x\in P\setminus R_j $ that is not comparable
to any element of $R_j $.
Then, because the width of $P$ is $w$,
we have that $x$ is the unique element
in
$R_{{\rm rank} (x)} $
that
is not comparable to any element of $R_j $.
Let $\Phi $ be any automorphism of $P$.
Because automorphisms preserve the rank, we have
$\Phi [R_j ]=R_j $.
Now $\Phi (x)$
must be the unique element
in
$R_{{\rm rank} (x)} $
that
is not comparable to any element of $R_j =\Phi [R_j ]$.
Hence
$\Phi (x)=x$, a contradiction to there not being any common fixed points for
all automorphisms.
Hence
every element $x\in P\setminus R_j $ is comparable
to at least one element of $R_j $.

Suppose, for a contradiction, that $R_{j+1} >R_j $.
Let $p\in P$ be an element of rank less than $j$ and let
$z\in \{ x\in P: {\rm rank} (x)\geq j+1\} $.
By the above, there is a $y\in R_j $ such that
$p<y\in R_j <R_{j+1 } $. Moreover, $z$, being of rank at least
$j+1$, is greater than or equal to some element of $R_{j+1} $.
Hence
$p<z$.
However, this means that
$
\{ x\in P: {\rm rank} (x)\geq j+1\} \not= \emptyset $
is
order-autonomous,
a contradiction.
Therefore,
there is an $x\in R_{j+1} $ such that
$|\downarrow x\cap R_j |\in \{ 1, \ldots , w-2\} $.

Let
$x\in R_{j+1} $ be so that
$\ell := |\downarrow x\cap R_j |\in \{ 1, \ldots , w-2\} $.
Let $S$ be any $\ell $-element subset of $R_j $.
Because $|R_j |=w-1$ and
$\alpha (R_j )= (w-1)!$,
there
is
a $\Phi \in {\rm Aut} (P)$ such that $\Phi [\downarrow x\cap R_j ]=S$.
Hence, there is
a $q_S :=\Phi (x)\in R_{j+1} $
with $\downarrow q_S \cap R_j =S$.
Because there are $\pmatrix{ w-1\cr \ell \cr } $
$\ell $-element
subsets of $R_j $
and
$\alpha (R_j )= (w-1)!$, there are
$\pmatrix{ w-1\cr \ell \cr } $ pairwise distinct elements
$q_S \in R_{j+1} $.
Because, for $\ell \not\in \{ 1, w-2\} $, we have
$\pmatrix{ w-1\cr \ell \cr } >w$, we infer that
$\ell \in \{ 1, w-2\} $.
Without loss of generality, we can assume that
$\ell =1$.

Let $R_j =\{ p_1 , \ldots , p_{w-1} \} $.
The set
$R_{j+1} $ contains
$w-1$ pairwise distinct elements $q_1 , \ldots , q_{w-1} $
such that, for $i=1, \ldots , w-1$, we have
$\downarrow q_i \cap R_j = \{ p_i \} $.
In particular, the set
$\{ p_1 , \ldots , p_{w-1} , q_1 , \ldots , q_{w-1} \} $
is a set $(w-1)C_2 $. (Note: $\ell =w-2$ would have produced a set
$S_{w-1} $.)

Finally,
suppose, for a contradiction, that
$R_{j+1} \not= \{ q_1 , \ldots , q_{w-1} \} $.
Then, because the width of $P$ is $w$, there is exactly one
element $z\in
R_{j+1} \setminus \{ q_1 , \ldots , q_{w-1} \} $.
If $|\downarrow z\cap R_j |\not= 1$, then
$z$ would be fixed by all automorphisms of $P$, which was excluded.
Hence
$|\downarrow z\cap R_j |= 1$.
Without loss of generality, we can assume
$\downarrow z\cap R_j = \{ p_1 \} $.
However, then, because
$\alpha (R_{j} )= (w-1)!$,
every $p_i \in R_j $
has at least two upper bounds in
$R_{j+1} $. This is a contradiction, because now every element of
$R_j $ has at least two upper bounds in $R_{j+1} $ and yet, every element of
$R_{j+1} $ has exactly one lower bound in $R_j $.
\qed

\vspace{.1in}

Recall that an ordered set $P$ is called {\bf coconnected}
iff there is no partition into two nonempty subsets
$P=B\cup T$ such that $B<T$.
Also recall that a group $G$ of automorphisms
{\bf acts transitively} on a
set $S$ iff for all $x,y\in S$ there is a $\Phi \in G$ with
$\Phi (x)=y$.

\begin{lem}
\label{maxlockobserv1}

Let $P$ be a
coconnected
ordered set
of width $w\geq 3$
such that no $p\in P$ is a common fixed point for all
automorphisms of $P$, and let
$k$ be such that
$
\alpha (R_k )\geq (w-1)!$.
Then
${\rm Aut} (P)$ acts transitively on $R_k $ and
$|R_k |=w$.

\end{lem}

{\bf Proof.}
First suppose, for a contradiction, that
${\rm Aut} (P)$ does not act transitively
on $R_k $. Then there is a set $A\subset R_k $
such that, for all $\Phi \in {\rm Aut} (P)$,
we have $\Phi [A]=A$ and
$\Phi [R_k \setminus A]=R_k \setminus A$.
Because no element of $P$ is a common fixed point for all
$\Phi \in {\rm Aut} (P)$,
we have that $2\leq |A|\leq w-2$.
However, then
$|{\rm Aut} (P)|\leq |A|!(w-|A|)!< (w-1)!$,
a contradiction.

Now suppose, for a contradiction, that
$|R_k |\not= w$.
Because $\alpha (R_k )\geq (w-1)!$, we infer
$|R_k |= w-1$
and then
$\alpha (R_k )= (w-1)!$.
The same argument as at the start of the proof of Lemma \ref{w-1factlem}
shows that
every $x\in P$ is comparable to an element of $R_k $.

Consider the case that
$\{ x\in P: {\rm rank} (x)\geq k+1\} \not= \emptyset $.
Because $P$ is coconnected
and every $x\in P$ is comparable to an element of $R_k $,
we have that
$\{ x\in P: {\rm rank} (x)\geq k+1\} $
is not order-autonomous.
By Lemma \ref{w-1factlem}, we obtain that
$\alpha (R_{k+1} )= (w-1)!$, and
that
$R_k \cup R_{k+1} $ is isomorphic to a set
$S_{w-1}$ or $(w-1)C_2 $.

Let $h\in {\mat N}$ be the height of $P$.
By repeating the argument above, we establish that,
for every $j\in \{ k, \ldots , h-1\} $,
the set
$R_j \cup R_{j+1} $ is isomorphic to a set
$S_{w-1}$ or $(w-1)C_2 $.
Consequently,
for every $j\in \{ k, \ldots , h\} $,
the set $R_j $ is also the set of
elements $R_{h-j} ^d $ of dual rank $h-j$.

By applying the same argument to the dual of
$P$, we obtain that, for every
$i\in \{ h-k, \ldots , h-1\} $,
the set
$R_i ^d \cup R_{i+1} ^d $ is isomorphic to a set
$S_{w-1}$ or $(w-1)C_2 $.
Consequently,
for every
$i\in \{ h-k, \ldots , h\} $,
the set
$R_i ^d $ is the set $R_{h-i} $.
Therefore,
for every $j\in \{ 0, \ldots , h-1\} $,
the set
$R_j \cup R_{j+1} $ is isomorphic to a set
$S_{w-1}$ or $(w-1)C_2 $.
This, however, implies that the width of $P$ is $w-1<w$, a contradiction.
\qed

\begin{define}
\label{maxlockdef}

Let $P$ be a coconnected ordered set
of width $w$
such that there is no
$p\in P$ that is fixed by all automorphisms of $P$.
$P$
will be called
{\bf max-locked}
iff
one of the following hold.
\begin{enumerate}
\item
$w\geq 3$, $w\not= 4$ and there is
a rank $k$ such that
$\alpha (R_k )\geq (w-1)!$.
\item
$w= 4$ and there is
a rank $k$ such that
$\alpha (R_k )>8$.
\end{enumerate}

\end{define}

Definition \ref{maxlockdef}
looks a little
strange for the case $w=4$,
but an $8$-crown need not satisfy conclusion
\ref{maxlockobserv2} of Lemma \ref{maxlockobserv} below
when it is
topped with a 2-antichain such that each element in the antichain
is above opposite maximal elements of the 8-crown.
See Definition \ref{forbconfdef}
for more on such problematic configurations.
Width $4$ is also the reason why Theorem \ref{allowedconf}
uses
$\max \{ w(U), 5\} $
instead of just $w(U)$.

To prove Lemma \ref{maxlockobserv},
we need a few results on products of factorials.
In this section, Lemma \ref{factprodcombine} is only used to prove
Lemma \ref{noverklem}, but it will be used more extensively in
Section \ref{orbgraphsec}
and beyond, specifically in
Lemmas \ref{boundsforQ} and \ref{addonetomaxlock}.

\begin{lem}
\label{factprodcombine}

Let $(b_1 , \ldots , b_m )$ be a vector
of $m>1$ integers greater than $1$,
sorted in nondecreasing order,
that is not equal to
any of the vectors
$(2,n)$ with $n\geq 2$, $(3,3)$, $(3,4)$, or $(2,2,2)$.
Then
$\prod _{j=1} ^m b_j !
\leq
\left( \sum _{j=1} ^m b_j -m \right) ! $
and the inequality is strict unless
$(b_1 , \ldots , b_m )$
is equal to $(2,2,3)$ or $(3,5)$.

\end{lem}

{\bf Proof.}
We start the proof by considering integers
$b_1 , \ldots , b_m $ that are greater than $1$ and
which need not be sorted in any order.
Note that,
trivially,
because all factors are
less than or equal to
$\max\{ b_j :j=1, \ldots , m\} $
and because there
are $m$ factors equal to $1$ in the product of the factorials,
we have
$\prod _{j=1} ^m b_j !
=
1\cdot \prod _{j=1} ^m \prod _{k=2} ^{b_j } k
<
\left( \sum _{j=1} ^m b_j -(m-1) \right) ! $.

Next, we establish a
simple claim that will be used repeatedly in the
remainder of the proof.

{\em Claim.
Let $k\in \{ 2, \ldots , m\} $ and let
$b_j ' $, $j=1, \ldots , k$ be
integers such that
$\prod _{j=1} ^k b_j '!
\leq
\left( \sum _{j=1} ^k b_j ' -k \right) ! $
and
such that,
for $j=1, \ldots , k$, we have
$2\leq b_j '\leq b_j $.
Then
$\prod _{j=1} ^m b_j !
\leq
\left( \sum _{j=1} ^m b_j  -m \right) ! $.
Moreover,
the second inequality is strict when $k<m$ or when the
first inequality is strict.}

We first
note the following auxiliary
inequality.
$$
\prod _{j=1} ^k b_j !
=
\prod _{j=1} ^k {b_j !
\over
b_j '!}
\prod _{j=1} ^k b_j '!
\leq
\prod _{j=1} ^k {b_j !
\over
b_j '!}
\left( \sum _{j=1} ^k b_j ' -k \right) !
\leq
\left( \sum _{j=1} ^k b_j  -k \right) ! $$

where the last inequality stems from the fact that, for all $j$,
we have that the
smallest factor
$b_{j} '+1$
of ${b_j !
\over
b_j '!}$ is at most
$\left( \sum _{j=1} ^k b_j ' -k \right) +1$.

Now, for $2\leq k<m$,
we can use this inequality to establish the claim, where we
use the trivial observation from the start
for the strict inequality below.
\begin{eqnarray*}
\prod _{j=1} ^m b_j !
& = &
\prod _{j=1} ^{k} b_j !
\prod _{j=k+1} ^m b_j !
\leq
\left( \sum _{j=1} ^k b_j -k \right) !
\prod _{j=k+1} ^{m} b_j !
\\
& < &
\left( \left( \sum _{j=1} ^k b_j -k \right)
+\sum_{j=k+1} ^{m}  b_j -((1+(m-k))-1)\right) !
=
\left( \sum _{j=1} ^m b_j -m \right) !
\end{eqnarray*}

The result is now proved
via the {\em Claim}
by considering the following cases.

{\em Case 1: There are no two
$i\not= k$ such that
$b_i =b_k$.}

In this case,
assume, without loss of generality, that $b_1 <b_2 <\cdots <b_m $.

If $m\geq 3$, we have
$b_1 \geq 2$, $b_2 \geq 3$, $b_3 \geq 4$.
Because $2!\cdot 3!\cdot 4!=288<720=6!=(2+3+4-3)!$,
the result follows from
the {\em Claim}.

This leaves the subcase $m=2$ and $b_1 <b_2 $.
Because the
vectors
$(2,n)$ (with $n\geq 2$)
are excluded, we only need to consider $b_1 \geq 3$.

For $b_1 \geq 4$, we have $b_2 \geq 5$.
Because
$4!\cdot 5! =6!\cdot 4<7!=(4+5-2)!$,
the result follows from
the {\em Claim}.

For $b_1 =3$,
because $(3,4)$ is excluded, we have
$b_2 \geq 5$.
Because $3!\cdot 5! =6!=(3+5-2)!$,
the result follows from
the {\em Claim}.

{\em Case 2: There are
$i\not= k$ such that
$b_i =b_k$.}

Without loss of generality, resort the $b_j $ so that
$i=1$ and $k=2$.

In case
$b_{1} =b_2 \geq 4$,
because $4!\cdot 4!=4!\cdot 4\cdot 6<6!=(4+4-2)!$,
the result follows from
the {\em Claim}.

In case
$b_{1} =b_2 = 3$,
because the vector
$(3,3)$ is excluded from consideration,
we must have $m\geq 3$.
Because
$2!\cdot 3!\cdot 3!=4!\cdot 3<5!=(2+3+3-3)!$,
the result follows from
the {\em Claim}.

In case
$b_{1} =b_2 = 2$,
because the vector
$(2,2)$ is excluded from consideration,
we must have $m\geq 3$.

First consider the case that there is a $b_i \geq 4$.
We can resort the $b_j $ so that $b_{3} \geq 4$.
Because
$2!\cdot 2!\cdot 4!=4!\cdot 4<5!=(2+2+4-3)!$,
the result follows from
the {\em Claim}.

Next consider the case that all $b_i $ are less than $4$
and there is a $b_i = 3$.
Trivially $2!\cdot 2!\cdot 3!=4!=(2+2+3-3)!$,
and
the non-strict inequality follows from
the {\em Claim}.
For the strict inequality,
we
can focus on the case $m\geq 4$ and we can assume
(after possibly resorting) that
$b_{4} =3$.
Because
$2!\cdot 2!\cdot 2!\cdot 3!=4!\cdot 2<5!=(2+2+2+3-4)!$,
the result follows from
the {\em Claim}.

This leaves the case that all the $b_i $ are equal to $2$.
Because $(2,2)$ and $(2,2,2)$ are excluded from
consideration, this means that
$m\geq 4$.
For $m\geq 4$ we have
$(2!)^m < (2m-m)!$, which completes the proof.
\qed

\begin{lem}
\label{noverklem}

Let $n\geq 6$ and let $k$ be a nontrivial divisor of
$n$.
Then $k!\left( \left( {n\over k} \right) ! \right) ^k < (n-1)!$.

\end{lem}

{\bf Proof.}
By Lemma \ref{factprodcombine},
we obtain
\begin{eqnarray*}
k!\left( \left( {n\over k} \right) ! \right) ^k
& < &
\left( k+k {n\over k} -(k+1)\right) !
=
(n-1)!
\end{eqnarray*}
\qed

\begin{lem}
\label{maxlockobserv}

Let $P$ be a max-locked ordered set
of width $w\geq 3$.
In case $w\not= 4$, let
$k$ be such that
$
\alpha (R_k )\geq (w-1)!$.
In case $w= 4$, let
$k$ be such that
$
\alpha (R_k )>8$.
Then the following hold.
\begin{enumerate}
\item
\label{maxlockobserv2}
If $R_k $ contains a non-maximal element,
then no
two elements of $R_k $
have equal sets of upper covers.

\item
\label{maxlockobserv3}
If $R_k $ contains a non-minimal element,
then no
two elements of $R_k $
have equal sets of lower covers.

\end{enumerate}

\end{lem}

{\bf Proof.}
We will only prove
part \ref{maxlockobserv2}.
Part \ref{maxlockobserv3} is proved dually.
To simplify the translation to the analogous proof for
part \ref{maxlockobserv3}, some parts of the proof below
are formulated to dualize easily.

Suppose, for a contradiction, that there are non-maximal
$a,b\in R_k $ that have the same strict upper bounds.
Clearly, $R_k $ can be partitioned into
sets $A_1 , A_2 , \ldots , A_n $
such that, for all
$\ell \in \{ 1, \ldots , n\} $
and for all
$x,y\in A_\ell $, we have
$\uparrow x\setminus \{ x\}  = \uparrow y\setminus \{ y\} $.
Without loss of generality, we can assume $a,b\in A_1 $.
Note that, for every $\Phi \in {\rm Aut} (P)$
and every $\ell \in \{ 1, \ldots , n\} $,
there is a $z\in \{ 1, \ldots , n\} $ such that
$\Phi [A_\ell ]=A_z $.

Because $\alpha (R_k )\geq (w-1)!$ and there is no $x\in P$ that is fixed by all
automorphisms,
by Lemma \ref{maxlockobserv1},
${\rm Aut} (P)$ acts transitively on $R_k $.
Therefore,
for all
$\ell \in \{ 1, \ldots , n\} $, we have
$|A_\ell |=|A_1 |\geq 2$.

If $n$ were equal to $1$, then all elements of $R_k $ would have the same
strict upper bounds.
By Lemma \ref{maxlockobserv1}, we have that $|R_k |=w$.
Because the width of $P$ is $w$, every
element of rank less than $k$ is below an element of rank $k$, and
every element of $\{ x\in P:{\rm rank} (x)>k\} $
is above an element of rank $k$.
Because any two elements of $R_k $ have the same strict upper bounds, we
first infer
$\{ x\in P:{\rm rank} (x)>k\} >R_k $, and then
$\{ x\in P:{\rm rank} (x)>k\} >\{ x\in P:{\rm rank} (x)\leq k\} $,
contradicting that $P$ is coconnected.
Thus $n\geq 2$.
We conclude that $w$ is not prime, and, in particular,
$w\not\in \{ 3,5\} $.

Because
every $\Phi \in {\rm Aut} (P)$
preserves the blocks of the partition into the $A_j $,
we obtain that
$
\alpha (R_k )
\leq
\ell ! \left( \left( {w\over \ell} \right) !\right) ^\ell $.
By Lemma \ref{noverklem}, for $w\geq 6$, this
leads to
$
\alpha (R_k )
\leq
\ell ! \left( \left( {w\over \ell} \right) !\right) ^\ell
<(w-1)!$, which is not possible.
Hence $w=4$.

Because $w=4$, we have
$n=2$ and $|A_1 |=|A_2|=2$.
We conclude
$\alpha (R_k )
\leq 8$,
a contradiction.
\qed

\begin{lem}
\label{maxlock2lev}

Let $P$ be a max-locked ordered set of height $1$ and width $w\geq 3$.
Then
$P$ is isomorphic to $S_w $ or $wC_2 $.

\end{lem}

{\bf Proof.}
By duality, without loss of generality,
for $w\not= 4$, we can assume that
$
\alpha (R_0 )\geq (w-1)!$,
and, for $w= 4$, we can assume that
$
\alpha (R_0 )>8$.
By Lemma \ref{maxlockobserv1},
$|R_0 |=w$ and
${\rm Aut} (P)$ acts transitively on $R_0 $.
By part \ref{maxlockobserv2}
of Lemma \ref{maxlockobserv},
no two minimal elements of $P$
have equal sets of strict upper bounds.

Because no two
elements of $R_0 $
have equal sets of strict upper bounds,
$R_0 $ does not contain any order-autonomous antichains.
By Lemma \ref{allbutoneac}
(with $A=R_0$),
for all $\Phi , \Psi \in {\rm Aut} (P)$, we have
that
$\Phi |_{R_0 } \not= \Psi |_{R_0 } $
implies
$\Phi |_{R_1 } \not= \Psi |_{R_1 } $,
and hence
$
\alpha (R_1 )
\geq
\alpha (R_0 )
\geq (w-1)!
$.
By Lemma \ref{maxlockobserv1},
we obtain that $|R_1 |=w$ and that
${\rm Aut} (P)$ acts transitively on $R_1 $.

Because
$|R_1 |=w$, $P$ is coconnected, and
${\rm Aut} (P)$ acts transitively on
$R_0 $, there is a number $u<w$ such that
every minimal element is below exactly $u$ elements of $R_1 $.
Because
$|R_1 |=w=|R_0 |$ and
${\rm Aut} (P)$ acts transitively on
$R_1 $,
every maximal element is above exactly
$u$ elements of $R_0 $.

Because we can replace $P$ with the ordered set
$P'$ of height $1$ such that, for all elements
$x\in R_0 $ and all elements $z\in R_1
$, we have
$x<z$ in $P'$ iff $x\not\sim z$ in $P$,
we only need to conduct the argument for
$u
\in \left\{ 1, \ldots ,  \left\lfloor {w\over 2} \right\rfloor \right\} $.

Suppose, for a contradiction that $P$ is
a max-locked ordered set of height $1$ and width $w$
that is not isomorphic to either of $S_w $ or $wC_2 $.

In case $u
=1$,
$P$ would be
an ordered set $wC_2 $, contradicting that $P$ is a counterexample.
Thus $u
\geq 2$.

In case $u
=2$,
$P$ would be a pairwise disjoint union of
crowns.
Because $u\leq \left\lfloor {w\over 2} \right\rfloor $, we have
$w\geq 4$.
For $w\geq 5$,
such unions have fewer than $(w-1)!$ automorphisms,
contradicting that
$
\alpha (R_0 )\geq (w-1)!$.
For $w=4$,
because
no two minimal elements have the same
strict upper bounds,
we
obtain that $P$ is a crown with $8$ elements,
contradicting that, for $w=4$, we have
$
\alpha (R_0 )
>8$.
Thus
$u
\in \left\{ 3, \ldots ,  \left\lfloor {w\over 2} \right\rfloor \right\} $.

Let $x$ be a fixed maximal element of $P$ and let
$M:=\downarrow x\setminus \{ x\} \subset R_0 $.
Under an arbitrary automorphism of $P$,
there are
$w$ possible images for $x$ (and hence for $M$), there are
$\leq u!$ ways to map $M$ to its image and
there are $\leq (w-u)!$ ways to
map $R_0 \setminus M$ to its image.
Thus
$(w-1)!\leq
\alpha (R_0 )
\leq
w\cdot u!\cdot (w-u)!$.
For $w\geq 9$ and $3\leq u\leq \left\lfloor {w\over 2} \right\rfloor $,
we have
$w\cdot u!
=
(w\cdot 3!) \cdot 4\cdots u
<(w-1)(w-2)\cdot 4\cdots u
<(w-1)\cdots (w-u+1)!$
and hence
$w\cdot u!\cdot (w-u)!<(w-1)!$, a contradiction.
Thus $w\leq 8$.
Because $3\leq u\leq \left\lfloor {w\over 2} \right\rfloor $,
we obtain $w\in \{ 6,7,8\} $.

Moreover, for
$w=8$ and $u=4$, we have
$8\cdot 4=32<35=7\cdot 5$, so
$w\cdot u!\cdot (w-u)!=
8\cdot 4!\cdot 4!=8\cdot 4\cdot 6\cdot 4!<
7\cdot 5\cdot 6\cdot 4!
=
7!=(w-1)!$, a contradiction.
Thus, in case $w=8$, we have $u\not= 4$,
which means that, in all remaining cases, we have $u=3$.

First consider the case that there are no two distinct elements $x,y\in R_1 $
such that $|\downarrow x ~ \cap \downarrow y|>1$.
Let $b\in R_0 $ and let $t_1 , t_2 , t_3 \in R_1 $ be its upper covers.
Then, for any two distinct $i,j\in \{ 1,2,3\} $, we have
$|\downarrow t_1 ~ \cap \downarrow t_2 ~|=1$.
Hence
$|\downarrow t_1 ~ \cup \downarrow t_2 ~ \cup \downarrow t_3 ~|=7$.
In particular, $w\in \{ 7,8\} $.
Let $\Phi $ be any automorphism of $P$.
There are $w$ possible images $\Phi (b)$, and this image determines the
image of the element (if there is one) of
$R_0 \setminus (\downarrow t_1 ~ \cup \downarrow t_2 ~ \cup \downarrow t_3 )$.
The remaining three doubletons
$(\downarrow t_j ~ \cap R_0 )\setminus \{ b\} $, must be mapped to the
corresponding three doubletons
$(\downarrow \Phi (t_j ) ~ \cap R_0 )\setminus \{ \Phi (b)\} $.
There are $3!\cdot 2!\cdot 2!\cdot 2!$ ways to do this.
In total,
$(w-1)! \leq \alpha(R_0) \leq
w \cdot 3!\cdot 2!\cdot 2!\cdot 2!
=2\cdot w\cdot 4!$.
Because $w\in \{ 7,8\} $ and
$6\cdot 5>2\cdot 7$
and
$7\cdot 6\cdot 5>2\cdot 8$, this is a contradiction.
We conclude that there are
distinct elements $x,y\in R_1 $
such that $|\downarrow x ~ \cap \downarrow y~|>1$,
and hence such that $|\downarrow x ~ \cap \downarrow y~|=2$.

Let $x,y\in R_1 $
such that
$|\downarrow x ~ \cap \downarrow y~|=2$.
Let $M:=\downarrow x\cap R_0$.
As above, under an arbitrary automorphism
$\Phi $
of $P$, there are $w$ possible images for $x$ (and hence for $M$) and there are $\leq 3!$ ways to map $M$ to its image.
Furthermore, $\Phi (y)$ must be above the images of all elements of
$\downarrow  y \cap M
=\{ b_1, b_2\} $ with $b_1\not= b_2$. Because $\Phi (x)>\Phi (b_1) , \Phi (b_2 )$
and $|\uparrow  \Phi (b_1 )\cap \uparrow  \Phi (b_2 )|
=|\uparrow  b_1 \cap \uparrow  b_2 | < 3$, there is
only a single point left for $\Phi (y)$. Hence, $\Phi (y)$ is determined by
$\Phi |_M$,
and {\em a fortiori} so is the image of the
unique point of
$(\downarrow y\setminus \downarrow x)\cap R_0 $.

The remaining
$(w-4)$ points in $R_0$
can be mapped to their images in at most $(w-4)!$ ways.
All together, $(w-1)! \leq \alpha(R_0) \leq w \cdot 3! \cdot (w - 4)!$.
However, $(w-1)(w-2)(w-3)>6w$, which is the final contradiction
that completes this proof.
\qed

\begin{lem}
\label{maxlockeddescr}

Let $P$ be a max-locked ordered set of width $w\geq 3$.
Then, for every $j$ such that $R_{j+1} \not= \emptyset $,
we have that $R_j\cup R_{j+1} $ is isomorphic to
$S_w $ or $wC_2 $.
In particular, $|{\rm Aut} (P)|=w!$

\end{lem}

{\bf Proof.}
Let $k$ be such that
$
\alpha (R_k )\geq (w-1)!$,
and, in case $w=4$, assume that
$\alpha (R_k )>8$.
By Lemma \ref{maxlockobserv1},
${\rm Aut} (P)$ acts transitively on $R_k $ and
$|R_k |=w$.

Consider the case $R_{k+1} \not= \emptyset $.
Because ${\rm Aut} (P)$ acts transitively on $R_k $ and $P$ is coconnected,
no element of $R_k $ is a lower bound of $R_{k+1} $.
Thus $R_k \cup R_{k+1} $ is coconnected, and hence
$R_k \cup R_{k+1} $ is max-locked of height $1$.
By Lemma \ref{maxlock2lev},
$R_k\cup R_{k+1} $ is isomorphic to
$S_w $ or $wC_2 $.
Let $h$ be the height of $P$.
Inductively, for all
$j\in \{ k, \ldots , h-1\} $,
we obtain that $R_j\cup R_{j+1} $ is isomorphic to
$S_w $ or $wC_2 $.

Because $P$ has width $w$, every element of rank
less than $k$ is comparable to an element of
rank $k$.
Hence, for all
$j\in \{ k, \ldots , h\} $, the set
$R_j $ is equal to the set $R_{h-j} ^d $
of elements of dual rank $h-j$.

Applying the same argument to $R_{h-k} ^d $ in the dual ordered set of $P$
proves that, for all
$j\in \{ h-k, \ldots , h-1\} $,
the set $R_j^d \cup R_{j+1} ^d $ is isomorphic to
$S_w $ or $wC_2 $.
Hence, for all
$j\in \{ h-k, \ldots , h\} $, the set
$R_j ^d $ is equal to the set $R_{h-j} $.
We conclude that, for all
$j\in \{ k, \ldots , h-1\} $,
the set $R_j\cup R_{j+1} $ is isomorphic to
$S_w $ or $wC_2 $.

The claim on the number of automorphisms is an
easy induction on the height $h$ of the ordered set $P$ in which
the union of any two consecutive nonempty ranks is isomorphic to
$S_w $ or $wC_2 $.
For the base step, if $h=1$, then
$P$ is isomorphic to
$S_w $ or $wC_2 $ and hence $|{\rm Aut} (P)|
=w!$.
For the induction step, let $h>1$ and assume
the claim holds for all ordered sets as indicated and
of height less than $h$.
Because $P\setminus R_h $ is such that the union of any two
consecutive nonempty ranks is isomorphic to
$S_w $ or $wC_2 $, we have
$|{\rm Aut} (P\setminus R_h )|
=w!$.
Because $R_{h-1} \cup R_h $
is isomorphic to
$S_w $ or $wC_2 $, we can apply Lemma \ref{allbutoneac}
and we obtain that every automorphism of $P\setminus R_h $ has at most
one extension to $P$.
Because every automorphism of
$P\setminus R_h $ has a natural extension to $P$, we obtain
$|{\rm Aut} (P)|
= |{\rm Aut} (P\setminus R_h )|
=w!$.
\qed

\begin{remark}

{\rm
Note that ordered sets such that,
for every $j$ such that $R_{j+1} \not= \emptyset $,
we have that $R_j\cup R_{j+1} $ is isomorphic to
$wC_2 $ need not be unions of pairwise disjoint chains:
There can be $k\geq 0$ and $\ell \geq 2$ such that
$R_k <R_{k+\ell } $.
The ordered set $2C_3 ^* $
from part \ref{forbconffirst1}
of Definition \ref{forbconffirst} is one such example.
}

\end{remark}

\begin{define}

We define coconnected ordered sets of width
$2$ with exactly 2 automorphisms to be {\bf max-locked}, too.

\end{define}

\section{Interdependent Orbit Unions}
\label{orbclustsec}

\def\updownharpoons{\upharpoonleft \! \downharpoonright }

The proof of Lemma \ref{maxlockeddescr}
shows how the action of the automorphism group on
a single set of elements of rank $k$ can, in
natural fashion, ``transmit vertically" though the whole ordered set.
In the situation of Lemma \ref{maxlockeddescr}, the sets $R_k $ happen to be
orbits
(see Definition \ref{grouporb} below) of $P$.
In general, the ``transmission" must
focus on the orbits, not the sets $R_k $, and it
can ``transmit" the action of the automorphism group
on one orbit $O$ to orbits that have no points that are
comparable to any element of $O$.
It should be noted that interdependent orbit unions
(see Definition \ref{interdeporbundef} below)
in graphs have also proven useful
in the set reconstruction of certain graphs, see \cite{SchrSetRec}.
Although the presentation up to Proposition \ref{getallfromiou}
translates directly to and from the corresponding results in
\cite{SchrSetRec}, all proofs are included to keep this presentation self-contained.

\begin{define}
\label{grouporb}

(Compare with Definition 8.1 in \cite{SchrSetRec}.)
Let $P$ be an ordered set, let $G$ be a subgroup of
${\rm Aut} (P)$ and let
$x\in P$.
Then the set
$G\cdot x :=\{ \Phi (x):\Phi \in G \} $
is called
the {\bf orbit of $x$ under the action of $G$} or the
{\bf $G$-orbit of $x$}.
Explicit mention of $G$ or $x$
can be dropped when there is only one group under consideration or when
specific knowledge of $x$ is not needed.
When no group $G$ is explicitly
mentioned at all, we assume by default that
$G={\rm Aut} (P)$.
The group generated by a single automorphism is denoted $\langle \Phi \rangle $.

\end{define}

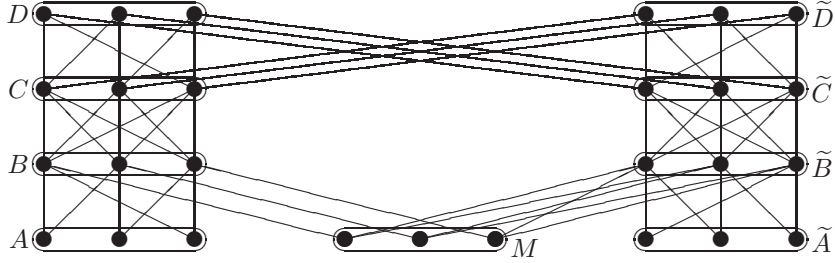
\begin{figure}

\centerline{
\unitlength 1mm 
\linethickness{0.4pt}
\ifx\plotpoint\undefined\newsavebox{\plotpoint}\fi 
\begin{picture}(107,36.5)(0,0)
\put(5,5){\circle*{2}}
\put(5,25){\circle*{2}}
\put(105,5){\circle*{2}}
\put(105,25){\circle*{2}}
\put(5,15){\circle*{2}}
\put(5,35){\circle*{2}}
\put(105,15){\circle*{2}}
\put(105,35){\circle*{2}}
\put(15,5){\circle*{2}}
\put(15,25){\circle*{2}}
\put(95,5){\circle*{2}}
\put(95,25){\circle*{2}}
\put(15,15){\circle*{2}}
\put(15,35){\circle*{2}}
\put(95,15){\circle*{2}}
\put(95,35){\circle*{2}}
\put(25,5){\circle*{2}}
\put(25,25){\circle*{2}}
\put(85,5){\circle*{2}}
\put(85,25){\circle*{2}}
\put(25,15){\circle*{2}}
\put(25,35){\circle*{2}}
\put(85,15){\circle*{2}}
\put(85,35){\circle*{2}}
\put(45,5){\circle*{2}}
\put(65,5){\circle*{2}}
\put(55,5){\circle*{2}}
\put(55,5){\circle*{2}}
\put(65,5){\circle*{2}}
\put(45,5){\circle*{2}}
\put(5,5){\line(0,1){10}}
\put(5,15){\line(0,1){10}}
\put(85,15){\line(0,1){10}}
\put(5,25){\line(0,1){10}}
\put(105,5){\line(0,1){10}}
\put(105,25){\line(0,1){10}}
\put(15,5){\line(0,1){10}}
\put(15,15){\line(0,1){10}}
\put(95,15){\line(0,1){10}}
\put(15,25){\line(0,1){10}}
\put(95,5){\line(0,1){10}}
\put(95,25){\line(0,1){10}}
\put(25,5){\line(0,1){10}}
\put(25,15){\line(0,1){10}}
\put(105,15){\line(0,1){10}}
\put(25,25){\line(0,1){10}}
\put(85,5){\line(0,1){10}}
\put(85,25){\line(0,1){10}}
\put(25,5){\line(-2,1){20}}
\put(25,15){\line(-2,1){20}}
\put(105,15){\line(-2,1){20}}
\put(5,15){\line(2,1){20}}
\put(85,15){\line(2,1){20}}
\put(25,25){\line(-2,1){20}}
\put(85,5){\line(2,1){20}}
\put(85,25){\line(2,1){20}}
\put(45,5){\line(-4,1){40}}
\put(65,5){\line(4,1){40}}
\put(55,5){\line(-4,1){40}}
\put(55,5){\line(4,1){40}}
\put(65,5){\line(-4,1){40}}
\put(45,5){\line(4,1){40}}
\put(5,5){\line(1,1){10}}
\put(5,15){\line(1,1){10}}
\put(85,15){\line(1,1){10}}
\put(25,15){\line(-1,1){10}}
\put(105,15){\line(-1,1){10}}
\put(5,25){\line(1,1){10}}
\put(105,5){\line(-1,1){10}}
\put(105,25){\line(-1,1){10}}
\put(15,5){\line(1,1){10}}
\put(15,15){\line(1,1){10}}
\put(95,15){\line(1,1){10}}
\put(15,15){\line(-1,1){10}}
\put(95,15){\line(-1,1){10}}
\put(15,25){\line(1,1){10}}
\put(95,5){\line(-1,1){10}}
\put(95,25){\line(-1,1){10}}
\put(3,5){\makebox(0,0)[rc]{\footnotesize $A$}}
\put(3,25){\makebox(0,0)[rc]{\footnotesize $C$}}
\put(107,5){\makebox(0,0)[lc]{\footnotesize $\widetilde{A} $}}
\put(107,25){\makebox(0,0)[lc]{\footnotesize $\widetilde{C} $}}
\put(67,5){\makebox(0,0)[lt]{\footnotesize $M$}}
\put(3,15){\makebox(0,0)[rc]{\footnotesize $B$}}
\put(3,35){\makebox(0,0)[rc]{\footnotesize $D$}}
\put(107,15){\makebox(0,0)[lc]{\footnotesize $\widetilde{B} $}}
\put(107,35){\makebox(0,0)[lc]{\footnotesize $\widetilde{D} $}}
\put(45,5){\line(5,1){50}}
\put(55,5){\line(5,1){50}}
\put(65,5){\line(2,1){20}}
\multiput(5,25)(.2693602694,.0336700337){297}{\line(1,0){.2693602694}}
\multiput(105,25)(-.2693602694,.0336700337){297}{\line(-1,0){.2693602694}}
\multiput(15,25)(.2693602694,.0336700337){297}{\line(1,0){.2693602694}}
\multiput(95,25)(-.2693602694,.0336700337){297}{\line(-1,0){.2693602694}}
\multiput(25,25)(.2693602694,.0336700337){297}{\line(1,0){.2693602694}}
\multiput(85,25)(-.2693602694,.0336700337){297}{\line(-1,0){.2693602694}}
\put(15,5){\oval(23,3)[]}
\put(95,5){\oval(23,3)[]}
\put(55,5){\oval(23,3)[]}
\put(15,15){\oval(23,3)[]}
\put(95,15){\oval(23,3)[]}
\put(15,25){\oval(23,3)[]}
\put(95,25){\oval(23,3)[]}
\put(15,35){\oval(23,3)[]}
\put(95,35){\oval(23,3)[]}
\end{picture}
}

\caption{An ordered set
$P$
with
orbits marked with ovals.
}
\label{transmit_drive}

\end{figure}

Note that,
if a strict subset $Q\subset P$ was obtained by
removing a union of ${\rm Aut} (P)$-orbits, then
${\rm Aut} (P)$-orbits
that are contained in $Q$
can be strictly contained in
${\rm Aut} (Q)$-orbits:
The ${\rm Aut} (P)$-orbits of the ordered set $P$ in
Figure \ref{transmit_drive} are marked by ovals.
We can see that, for $X\in \{ A,B,C,D\} $,
the ${\rm Aut} (P)$-orbits $X$ and $\widetilde{X} $
are strictly contained in
the ${\rm Aut} (P\setminus M)$-orbit $X\cup \widetilde{X} $.
For this reason (further elaborated later
in Remark \ref{orbstrictcont}),
dictated orbit structures will be
useful for the representation of ${\rm Aut} (P)$ in Proposition \ref{getallfromiou}
and they are vital for the induction
proof of
Theorem \ref{allowedconf}.

\begin{define}
\label{dictatedef}

(Compare with Definition 8.2 in \cite{SchrSetRec}.)
Let $P$ be an ordered set and let ${\cal D}$ be a partition of $P$
into antichains. Then ${\rm Aut} _{\cal D} (P)$ is the set of automorphisms
$\Phi :P\to P$ such that, for every $\Phi $-orbit
$O:=\langle \Phi \rangle \cdot x$ of $\Phi $, there is a
$D\in {\cal D}$ such that $O\subseteq D$.
In this context, the partition ${\cal D}$ is called a
{\bf dictated orbit structure} for $P$, and
the pair $(P,{\cal D})$ is called a {\bf structured ordered set}.
${\rm Aut} _{\cal D} (P)$-orbits will, more briefly, be called
{\bf ${\cal D}$-orbits}.

The partition of $P$ into its ${\rm Aut} (P)$-orbits is called the
{\bf natural orbit structure} of $P$, which will typically be denoted ${\cal N}$.
When working with the natural orbit structure,
explicit indications of the automorphism set, usually via subscripts
or prefixes
${\cal D}$, will often be
omitted.

\end{define}

Clearly,
${\rm Aut} _{\cal D} (P)$ is a subgroup of the
automorphism group
${\rm Aut} (P)$.
Moreover,
${\rm Aut} _{\cal D} (P)\not= {\rm Aut} (P)$
iff there are an orbit $O$ of $P$ and a $D\in {\cal D}$
such that $O\cap D\not= \emptyset $ and
$O\not\subseteq D$.

The proof of Lemma \ref{maxlockeddescr} has already shown
that, for an automorphism $\Phi $,
the values on a single orbit of $\Phi $
can completely determine
$\Phi $.
The relation of direct interdependence in Definition \ref{directdep}
below
provides a more detailed view of this observation as well as the
simple observation in
Lemma \ref{allbutoneac}.

\begin{define}
\label{directdep}

(Compare with Definition 8.3 in \cite{SchrSetRec}.)
Let
$(P,{\cal D})$ be a structured ordered set
and let $C,D$ be two ${\cal D}$-orbits
of $P$ such that there are a $c_1 \in C$
and a $d_1 \in D$ such that $c_1 <d_1 $.
We will write $C\upharpoonleft _{\cal D} D$
and $D\downharpoonright _{\cal D} C$
iff there are $c_2 \in C$ and a $d_2 \in D$
such that $c_2 \not\sim d_2 $.
In case $C\upharpoonleft _{\cal D} D$ or $C\downharpoonright _{\cal D} D$, we write
$C\updownharpoons _{\cal D} D$ and say that
$C$ and $D$ are
{\bf directly interdependent}.

\end{define}

Figure \ref{transmit_drive} shows
how
(connection through)
direct interdependence can
allow one orbit to determine the values of
automorphisms on many other orbits:
In the ordered set in
Figure \ref{transmit_drive}, we have
$A\upharpoonleft _{\cal N} B\downharpoonright _{\cal N} M
\upharpoonleft _{\cal N} \widetilde{B}
\downharpoonright _{\cal N} \widetilde{A}$ and the
values of any automorphism on
$A\cup B\cup M\cup \widetilde{B}
\cup \widetilde{A}$
are determined the by automorphism's values on $A$.
Note, however, that even direct interdependence $C\updownharpoons _{\cal N} D$
does not mean that
automorphisms are determined by the values on either
orbit
$C$ or $D$:
Consider a $6$-crown in which every element is replaced with an
order-autonomous $2$-antichain:
The minimal elements form an orbit, as do the maximal elements, the two
orbits are directly interdependent, but no automorphism is
completely determined solely
by its values on the minimal (or the maximal) elements.
This idea will be explored further in Section \ref{boundnra}.

We now turn our attention to interdependent orbit unions, which
are the unions of the connected components of the orbit graph
defined in Definition \ref{orbitgraph} below.
Orbit graphs themselves will take center stage starting in Section
\ref{orbgraphsec}.

\begin{define}
\label{orbitgraph}

Let
$(P, {\cal D} )$
be
a structured ordered set.
The {\bf orbit graph}
${\cal O} (P, {\cal D})$
of
$P$
is defined to be the graph whose vertex set is
the set of all ${\cal D}$-orbits
such that two ${\cal D}$-orbits
$D_1 $ and $D_2 $ are adjacent iff
$D_1 \updownharpoons _{\cal D} D_2 $.

\end{define}

Although our main focus will
ultimately be on orbit graphs
for structured ordered sets in which every $D\in {\cal D}$ is a ${\cal D} $-orbit,
also see part \ref{structset3}
of Definition \ref{structset} below,
until that time, we must keep in mind that the
$D\in {\cal D}$ need not be orbits themselves.

\begin{define}
\label{interdeporbundef}

(Compare with Definition 8.5 in \cite{SchrSetRec}.)
Let
$(P, {\cal D} )$
be
a structured ordered set and let ${\cal O} (P,{\cal D} )$ be its orbit graph.
If ${\cal E} $ is a connected component of
${\cal O} (P,{\cal D} )$,
then we will call the set $\bigcup {\cal E}\subseteq P$ an
{\bf interdependent ${\cal D}$-orbit union}.

\end{define}

\begin{define}

(Compare with Definition 8.6 in \cite{SchrSetRec}.)
Let
$(P,{\cal D})$ be a structured ordered set
and let
$Q\subseteq P$ such that, for all
${\cal D}$-orbits $D$,
we have $D\subseteq Q $ or $D\cap Q =\emptyset $.
The
dictated orbit structure for $Q$
{\bf induced by ${\cal D}$}, denoted ${\cal D} |Q $,
is defined to be the set of all ${\cal D}$-orbits
that are contained in $Q $.

For the natural orbit structure ${\cal N}$ for $P$,
the dictated orbit structure for $Q$ induced by ${\cal N} $
will be called the {\bf naturally required} automorphism structure
${\cal N}|Q$.

\end{define}

Note that, for induced orbit structures ${\cal D} |Q$, every
$D\in {\cal D} |Q$ is an orbit.
Proposition \ref{unionplacement}
below
now shows how
interdependent ${\cal D}$-orbit unions
reside in an ordered set.
It also lays the groundwork for representing
automorphisms through certain automorphisms on the
non-singleton interdependent orbit unions in
Proposition \ref{getallfromiou}.

\begin{prop}
\label{unionplacement}

(Compare with Proposition 8.7 in \cite{SchrSetRec}.)
Let
$(P,{\cal D})$ be a structured ordered set
and let
$U$ be an interdependent ${\cal D}$-orbit union.
Then the following hold.
\begin{enumerate}
\item
\label{unionplacement0}
For all $x\in P\setminus U$ and all
$C\in {\cal D}|U$, the following hold.
\begin{enumerate}
\item
\label{unionplacement1}
If there is a $c\in C$ such that $c<x$, then $C<x$.
\item
\label{unionplacement2}
If there is a $c\in C$ such that $c>x$, then $C>x$.
\end{enumerate}
\item
\label{unionplacement3a}
For every $\Phi \in {\rm Aut} _{\cal D} (P)$, we have that
$\Phi |_U \in {\rm Aut} _{{\cal D}|U} (U)$.
In particular, this means that the
${\cal D}|U$-orbits
are just the sets in
${\cal D}|U$.
\item
\label{unionplacement3b}
For
every $\Psi \in {\rm Aut} _{{\cal D}|U} (U)$, the function
$\Psi ^P (x) := \cases{
\Psi (x); & if $x\in U$, \cr
x; & if $x\in P\setminus U$,
} $
is an automorphism of $P$.
\end{enumerate}

\end{prop}

{\bf Proof.}
To prove part \ref{unionplacement1},
let $x\in P\setminus U$ and
$C\in {\cal D}|U$ be so that
there is a $c\in C$ such that $c<x$.
Let $X$ be the ${\cal D}$-orbit of $x$.
Because
$x\not\in U$, we must have
$X\not\updownharpoons _{\cal D} C$.
Because $C\ni c<x\in X$, we must have that
$C<X$ and hence $C<x$.

Part \ref{unionplacement2} is proved dually.

Part \ref{unionplacement3a}
follows directly from the definitions.

To prove part \ref{unionplacement3b},
let $\Psi \in {\rm Aut} _{{\cal D}|U} (U)$.
Clearly, $\Psi ^P $ is bijective.
To prove that $\Psi ^P $ is
order-preserving, let $x<y$.
If $x,y$ are both in $U$ or both in $P\setminus U$, we obtain
$\Psi ^P (x)<\Psi ^P (y)$.
In case $x\in P\setminus U$ and $y\in U$,
let $Y\in {\cal D}|U$ be so that $y\in Y$.
Then $\Psi (y)\in Y$.
By part \ref{unionplacement2},
we have that $x<Y$ and hence
$\Psi ^P (x)=x<Y\ni \Psi ^P (y)$.
The case in which $y\in P\setminus U$ and $x\in U$
is handled dually.
\qed

\begin{define}

(Compare with Definition 8.8 in \cite{SchrSetRec}.)
Let
$(P,{\cal D})$ be a structured ordered set
and let
$U$ be an interdependent ${\cal D}$-orbit union.
We define ${\rm Aut} _{{\cal D}|U} ^P (U)$ to be the set of
automorphisms $\Psi ^P \in {\rm Aut} (P)$ as in
part \ref{unionplacement3b} of Proposition \ref{unionplacement}.

\end{define}

Let $P$ be an ordered set, let ${\cal D}$ be a
dictated orbit structure for $P$ and let
$U, U'$ be disjoint interdependent ${\cal D}$-orbit unions.
Then, clearly, for
$\Psi ^P \in {\rm Aut} _{{\cal D}|U} ^P (U)$
and
$\Phi ^P \in {\rm Aut} _{{\cal D}|U'} ^P (U')$, we have
$\Psi ^P \circ \Phi ^P =\Phi ^P \circ \Psi ^P $.
Hence the following definition is sensible.

\begin{define}

(Compare with Definition 8.9 in \cite{SchrSetRec}.)
Let $P$ be an ordered set and let
${\cal A}_1 , \ldots , {\cal A} _z \subseteq {\rm Aut} (P)$ be
sets of automorphisms such that, for all pairs of distinct
$i,j\in \{ 1, \ldots , z\} $,
all $\Phi _i \in {\cal A} _i $ and all $\Phi _j \in {\cal A} _j $,
we have $\Phi _i \circ \Phi _j = \Phi _j \circ \Phi _i $.
We define $\bigcirc _{j=1} ^z
{\cal A} _j $
to be the set of
compositions $\Psi _1 \circ \cdots \circ \Psi _z $
such that, for $j=1, \ldots , z$, we have
$\Psi _j \in
{\cal A}_j $.

\end{define}

\begin{prop}
\label{getallfromiou}

(Compare with Proposition 8.10 in \cite{SchrSetRec}.)
Let $P$ be an
ordered set
with natural orbit structure ${\cal N}$
and let $U_1 , \ldots , U_z $
be the non-singleton interdependent orbit unions of $P$.
Then
${\rm Aut} (P)=\bigcirc _{j=1} ^z
{\rm Aut} _{{\cal N} |U_j } (U_j )$, and consequently
$
|{\rm Aut} (P)|=\prod _{j=1} ^z \left|{\rm Aut} _{{\cal N} |U_j } ^P (U_j )\right|
$.

\end{prop}

{\bf Proof.}
The containment ${\rm Aut} (P)\supseteq \bigcirc _{j=1} ^z
{\rm Aut} _{{\cal N} |U_j } ^P (U_j )$
follows from part \ref{unionplacement3b} of Lemma \ref{unionplacement}.

By part \ref{unionplacement3a} of Lemma \ref{unionplacement},
for every $\Phi \in {\rm Aut} (P)$ and every $j\in \{ 1, \ldots , z\} $,
we have that
$\Phi |_{U_j } \in {\rm Aut} _{{\cal N}|U_j } (U_j)$.
Because $\Phi $ fixes
all points in $P\setminus \bigcup _{j=1} ^z U_j $,
we have
$\Phi =\Phi |_{U_1 } ^P \circ \cdots \circ \Phi |_{U_z } ^P $.
Hence
${\rm Aut} (P)\subseteq \bigcirc _{j=1} ^z
{\rm Aut} _{{\cal N} |U_j } ^P (U_j )$.
\qed

\begin{remark}
\label{orbstrictcont}

(Compare with Remark 8.11 in \cite{SchrSetRec}.)
{\rm
Although the
naturally required automorphism structure
${\cal N}|U$
may look more technical than natural, it is indispensable for
the representation in Proposition \ref{getallfromiou}.
Consider the ordered set in Figure \ref{transmit_drive}.
The natural orbit unions in this ordered set are
$U_1 :=
A\cup B\cup M\cup \widetilde{B}\cup \widetilde{A} $
and
$U_2 :=
C\cup D\cup \widetilde{D}\cup \widetilde{C} $.
However, when considering $U_2 $
as an ordered set by itself, the
${\rm Aut} (U_2 )$-orbits
are $C\cup \widetilde{C}$ and
$D\cup \widetilde{D}$, whereas the
${\rm Aut} (P)$-orbits in $U_2 $
are
$C, D, \widetilde{C}$ and
$\widetilde{D}$.
Hence, we cannot use
the automorphism groups
${\rm Aut} (U_j )$
in place of their subgroups
${\rm Aut} _{{\cal N} |U_j } (U_j )$
in Proposition \ref{getallfromiou}.

The same effect was observed
in the introduction in the case that
$M$ is removed from the ordered set.
This is
why dictated orbit structures are also
crucial for the induction in the proof of
Theorem \ref{allowedconf}.
\qex
}

\end{remark}

Proposition \ref{getallfromiou} allows us to focus our
efforts to bound the number of automorphisms on ordered sets that
consist of a single interdependent orbit union.
Lemma \ref{w-1factexclude} below makes the,
admittedly simple, connection to max-locked ordered sets.

\begin{lem}
\label{w-1factexclude}

Let $P$ be an ordered set of width $\leq w$
that consists of a single ${\rm Aut} (P)$-orbit union
and which has
an
${\rm Aut} (P)$-orbit $X$
such that ${\rm Aut} (P)$ induces $N\geq (w-1)!$
permutations on $X$.
In case $w=4$, assume that
${\rm Aut} (P)$ induces $N>8$
permutations on $X$.
Then $P$ is max-locked.

\end{lem}

{\bf Proof.}
Because $P$ consists of a single
${\rm Aut} (P)$-orbit union, $P$ is
a coconnected ordered set
of width $w$
such that there is no
$p\in P$ that is fixed by all automorphisms of $P$.
With $k$ such that $X\subseteq R_k $,
the hypothesis guarantees that
$\alpha (R_k )\geq (w-1)!$ and, in case
$w=4$, $\alpha (R_k )>8$.
\qed

\vspace{.1in}

Lemma \ref{halfdictateslemmaxlock} below can be considered the prototype for
the key result in this paper, Theorem \ref{allowedconf}.
In most interdependent orbit unions in a non-max-locked ordered set,
the
upper bound on the number of automorphisms is
a product of factorials. By Lemma \ref{maxlockeddescr},
the factorials cannot exceed $(w(P)-1)!$ (for $w(P)=4$, some cases need to be
handled separately), and the sum of the numbers whose factorials are
multiplied will be bounded by half the size of the
interdependent orbit union minus the number of
(dictated) orbits in the interdependent orbit union.
We introduce some terminology to express this idea after
Lemma \ref{halfdictateslemmaxlock}.

\begin{lem}
\label{halfdictateslemmaxlock}

Let $P$ be a
connected
non-max-locked
ordered set of width $w\geq 3$
and let
$U\subset P$
be a max-locked
interdependent orbit union with $m\geq 3$ nonempty ranks and
width $a\leq w-1$ such that $(m,a)\not= (3,2)$, that is, such that
$P$ is not isomorphic to $2C_3 $.
Then
$|{\rm Aut } _{{\cal D}|U} (U)|
=a!
\leq
\min \left\{ (w-1)!,
\left\lfloor {1\over 2}
\left(
|U|-m \right)
\right\rfloor !\right\} $.

\end{lem}

{\bf Proof.}
By Lemma \ref{maxlockeddescr}, we have
$|{\rm Aut } _{{\cal D}|U} (U)|
=a!\leq (w-1)!$.
Because $U$ is max-locked, it consists of
$m$ natural orbits with $a$ elements each.
Hence
$\left\lfloor {1\over 2}
\left(
|U|-m \right)
\right\rfloor
=
\left\lfloor {1\over 2}
\left(
ma-m \right)
\right\rfloor
=
\left\lfloor {1\over 2}
\left(
m(a-1) \right)
\right\rfloor .
$
In case $m\geq 4$, we have
that the last term is at least
$2(a-1)\geq a$.
In case $m=3$ and $a\geq 3$, we have that the last term is
$
\left\lfloor {1\over 2}
\left(
3(a-1) \right)
\right\rfloor
=
\left\lfloor {1\over 2}
\left(
3a-3 \right)
\right\rfloor
\geq
\left\lfloor {1\over 2}
\left(
2a \right)
\right\rfloor
=
a
.
$
\qed

\vspace{.1in}

We conclude this section by introducing ideas and terminology that will be
fundamental for the remainder of this paper.

\begin{define}

We will say that the structured ordered set
$(P,{\cal D} _P )$ is
{\bf (dually) isomorphic}
to the structured ordered set
$(Q,{\cal D} _Q )$ iff there is a (dual) isomorphism
$\Phi :P\to Q$ such that, for all
$D\in {\cal D} _P $, we have
$\Phi [D]\in {\cal D} _Q $.

\end{define}

\begin{define}
\label{structset}

Let $(P, {\cal D} )$ be a structured ordered set.
\begin{enumerate}
\item
\label{structset1}
We call $(P,{\cal D})$ an {\bf interdependent orbit union}
iff
$\bigcup {\cal D} $
is an interdependent ${\cal D}$-orbit union.

\item
\label{structset2}
We say $(P,{\cal D})$ is
{\bf without slack}
iff
none of the sets in ${\cal D}$ contains a nontrivial order-autonomous
antichain.

\item
\label{structset3}
We call $(P,{\cal D})$ {\bf tight}
iff
it is without slack and
each $D \in {\cal D} $ is a
${\cal D}$-orbit, that is,
it is without slack and
${\rm Aut} _{\cal D} (P)$ acts transitively on
every $D \in {\cal D} $.

\item
\label{structset4}
If $(U,{\cal D})$ is an interdependent orbit union that is
not max-locked and such that $U$ is not a singleton,
we say $(U,{\cal D})$ is {\bf flexible}.
\end{enumerate}

\end{define}

By Lemma \ref{halfdictateslemmaxlock}, our main focus will be on
flexible interdependent orbit unions $(U,{\cal D})$.
(Max-locked interdependent orbit unions of height $1$
or of width $2$ and height $3$ will be handled separately.)
Because we will first focus on indecomposable ordered sets, we
will assume that $(U,{\cal D})$ is without slack.
Because we can
always refine the dictated orbit structure, we
are free to assume that $(U,{\cal D})$ is tight.

Similar to Lemma \ref{halfdictateslemmaxlock} we define the following.

\begin{define}
\label{adequatebound}

Let
$(P,{\cal D})$ be a
structured ordered set,
let
$w\geq 3$
and let $o\in {\mat Z} $.
Then $(P, {\cal D})$ is called {\bf
offset by $o$ from being
$w$-adequately bounded}
iff
there are numbers
$w_1 , \ldots , w_M \in \{ 0, \ldots , w-1\} $ such that
$|{\rm Aut } _{{\cal D}} (P)|
\leq
\prod _{j=1} ^M w_j !
$
and such that
$
\sum  _{j=1} ^M w_j
\leq
\left\lfloor {1\over 2}
\left(
|P|-|{\cal D} | +o\right)
\right\rfloor $.
When $o$ can be chosen to be zero,
$(P, {\cal D})$ is called {\bf
$w$-adequately bounded}.

\end{define}

Note that, although $w$ will often be the width of the ordered set,
in Definition \ref{adequatebound} above, it is just a number.
Lemma \ref{addonetomaxlock} and Theorem \ref{allowedconf}
show the ultimate use of this idea.
Similarly, for some specific constructions,
we will choose $M=m=|{\cal D} |$ and some $w_j !$ will be a bound on the number of
permutations induced by ${\rm Aut} _{\cal D} (U)$
on $D_j $, but
we should note that there is no
formal connection between the $w_j $ and the
$D_j $ other than the final inequality in
Definition \ref{adequatebound}.

By Lemma \ref{halfdictateslemmaxlock}
max-locked
interdependent orbit unions with $m\geq 3$ nonempty ranks and
width $a\leq w-1$ such that $(m,a)\not= (3,2)$
are $w$-adequately bounded.
The ordered set $2C_3 $
is offset by 1 from being $3$-adequately bounded.
We denote the natural dictated orbit structure of
$2C_3 $ by ${\cal D} (2C_3 )$.
By Lemma \ref{maxlock2lev},
max-locked interdependent orbit unions of height $1$
and width $w$
have $w!$ automorphisms, $2w$ elements and $2$ orbits,
so they are
are offset by $2$ from being $w$-adequately bounded.
From this point forward through the end of Section
\ref{boundnra}, our goal is to prove that most
non-max-locked interdependent
orbit unions $U$
are $\max \{ w(U), 5\} $-adequately bounded
and to treat the remaining cases separately.
Section \ref{orbgraphsec}
provides
the requisite structural insights and some
inequalities,
Section \ref{forbconfsec} provides the,
unfortunately a bit technical,
base step of the induction, which is then completed in
Section
\ref{boundnra}.

\section{The Orbit Graph of an Interdependent Orbit Union}
\label{orbgraphsec}

The idea for producing a bound on the number of
automorphisms for a flexible tight interdependent orbit
union is an induction on the number of dictated orbits.
Recall that, for tight interdependent orbit unions
$(U, {\cal D} )$, we have that
every $D\in {\cal D}$ is an orbit.
Hence the vertex set of
${\cal O} (U, {\cal D})$ is ${\cal D}$.
When a dictated orbit $D_n \in {\cal D}$ is removed,
the orbit graph can become disconnected and
the resulting structured ordered set may contain
nontrivial order-autonomous antichains.
We start with some standard notation that
assures that the orbits with indices smaller than $n$
form a component of the resulting
graph
${\cal O} (U, {\cal D})-D_{n} $, and
which allows easy reference to orbits with
certain properties.
This notation
will be used throughout this section.

\begin{nota}
\label{standardOGnotation}

Throughout this section,
$(U, {\cal D} )$
will be
a flexible tight
interdependent orbit union with $|{\cal D} |\geq 3$,
with the labeling of the elements of
${\cal D} =\{ D_1 , \ldots , D_m \} $
and $n\geq 3$ chosen
so that the following hold.
\begin{enumerate}
\item
$\{ D_1 , \ldots , D_{n-1} \} $
is a connected component of
the graph ${\cal O} (U, {\cal D})-D_{n} $.
\item
The orbits that are directly interdependent with
$D_n $ are $D_s , \ldots , D_{n-1} $ and $D_r , \ldots , D_m $.
\item
The orbits that contain nontrivial
$U\setminus \bigcup _{i=n} ^m D_i $-order-autonomous antichains are
$D_t , \ldots , D_{n-1} $.
\end{enumerate}

\end{nota}

Note that
Notation \ref{standardOGnotation} allows for $n=m$, that is, for
$D_n $ to not be a cutvertex of
${\cal O} (U, {\cal D})$.
Moreover, because
$(U, {\cal D} )$
is an
interdependent orbit union,
${\cal O} (U, {\cal D})$ is connected, so $s\leq n-1$.
For the other parameters,
$r> m$ or $t>n-1$ shall indicate that there are no orbits
as described via $r$ or $t$.

To easily refer to the (possibly trivial)
order-autonomous antichains in the orbits
$D_s , \ldots , D_{n-1} $,
we introduce the following notation.

\begin{nota}

For every $j\in \{ s, \ldots , n-1\} $,
we let $A_1 ^j , \ldots , A_{\ell _j} ^j $ be the
maximal
$U\setminus \bigcup _{i=n} ^m D_i $-order-autonomous antichains
that partition $D_j $. We set
${\cal A} ^j :=\{ A_1 ^j , \ldots , A_{\ell _j} ^j \} $.
Moreover,
for every $i\in \{ 1, \ldots , \ell _j \} $, we choose
a fixed element $a_i ^j \in A_i ^j $.
Note that, for $s\leq j\leq t-1$, the sets $A_i ^j $ are singletons.

\end{nota}

The fact that automorphisms map maximal order-autonomous antichains to
maximal order-autonomous antichains motivates the definition below, which
will frequently be used and which leads to our first insights.

\begin{define}

Let $S$ be a set, let ${\cal A} =\{ A_1 , \ldots , A_\ell \} $
be a partition of $S$ and let $\Phi :S\to S$ be a permutation of $S$.
We say that {\bf $\Phi $ respects the partition ${\cal A}$}
iff, for every $i\in \{ 1, \ldots , \ell \} $, there is
a $j\in \{ 1, \ldots , \ell \} $ such that $\Phi [A_i ]=A_j $.

\end{define}

\begin{lem}
\label{pruneorbit1}

For every $j\in \{ s, \ldots n-1\} $ we have that
$\ell _j >1$ and
every
$\Phi \in {\rm Aut} _{\cal D} (U)$
respects the partition
${\cal A} ^j $
of $D_j $.
Moreover,
for all $i,k\in \{ 1, \ldots , \ell _j \} $,
we have
$\left| A_i ^j \right| =\left| A_k ^j \right| $.

\end{lem}

{\bf Proof.}
Let
$j\in \{ s, \ldots , n-1\} $.
Because $n\geq 3$
and
$\{ D_1 , \ldots , D_{n-1} \} $
is a connected component of
the graph ${\cal O} (U, {\cal D})-D_{n} $, we have that
$D_j $ is directly interdependent with another
$D_{j'} $ with $j'\in \{ 1, \ldots , n-1\} \setminus \{ j\} $.
Consequently,
$D_j $ itself is not a
$U\setminus \bigcup _{i=n} ^m D_i $-order-autonomous antichain.
Hence $\ell _j >1$.

Let
$\Phi \in {\rm Aut} _{\cal D} (U)$.
Because
$\Phi |_{U\setminus \bigcup _{i=n} ^m D_i } \in
{\rm Aut} _{{\cal D} \setminus \{ D_n , \ldots , D_m \} }
(U\setminus \bigcup _{i=n} ^m D_i )$, we have
$\Phi [D_j ]=D_j $,
and because automorphisms map
maximal order-autonomous antichains to
maximal order-autonomous antichains,
$\Phi $ must map every $A_i ^j $ to another
$A_k ^j $, that is,
$\Phi $
respects
${\cal A} ^j $.

Because $(U, {\cal D})$ is tight,
for any $i,k\in \{ 1, \ldots , \ell _j \} $,
there is a
$\Phi \in {\rm Aut} _{\cal D} (U)$ such that
$\Phi \left( a_i ^j \right) =a_k ^j $.
Consequently, because
$\Phi $
respects
${\cal A} ^j $, we have
$\Phi \left[ A_i ^j \right] =A_k ^j $, and hence
$\left| A_i ^j \right| =\left| A_k ^j \right| $.
\qed

\begin{lem}
\label{pruneorbit3}

For all
$j\in \{ t, \ldots , n-1\} $,
$i\in \{ 1, \ldots , \ell _j \} $
and distinct $x,y\in A_i ^j $, there
is a $d\in D_n $ such that $d$ is comparable to one of
$x$ and $y$, but not the other.

\end{lem}

{\bf Proof.}
Let
$j\in \{ t, \ldots , n-1\} $ and
$i\in \{ 1, \ldots , \ell _j \} $.
Because $D_j $ is not directly interdependent with any
$D_k $ with $k>n$, we have that
$A_i ^j $ is order-autonomous in $U\setminus D_n $.
Because, for any distinct $x,y\in A_i ^j $, the set
$\{ x,y\} $ is not order-autonomous in $U$,
there
must be a $d\in D_n $ such that $d$ is comparable to one of
$x$ and $y$, but not the other.
\qed

\vspace{.1in}

In our analysis, the order-autonomous antichains
$A_i ^j $ will be collapsed into
singletons. Hence we introduce the following and again establish some natural
properties.

\begin{define}
\label{phindef}

We define the
{\bf $D_{n} $-pruned and compacted}
ordered set $U_n :=\bigcup _{i=1} ^{s-1} D_i \cup \bigcup _{j=s} ^{n-1}
\left\{ a_i ^j : i\in \{ 1, \ldots , \ell _j \} \right\} $
and we define
\begin{eqnarray*}
{\cal D}_n
& := &
\{ D_j \cap U_n :j \in \{ 1, \ldots , n-1\} \}
\\
& = &
\{ D_j : j\in \{ 1, \ldots , s-1\} \} \cup
\left\{ \left\{ a_1 ^j , \ldots , a_{\ell _j } ^j \right\} :  j\in \{ s, \ldots , n-1\} \right\} .
\end{eqnarray*}
For every
$\Phi \in {\rm Aut} _{\cal D} (U)$, we define
the function $\Phi _n :U_n \to U_n  $
by
$\Phi _n |_{\bigcup _{j=1} ^{s-1} D_j } :=\Phi |_{\bigcup _{j=1} ^{s-1} D_j } $,
and by, for any $j\in \{ s, \ldots , n-1\} $ and
$i\in \{ 1, \ldots , \ell _j \} $, setting $\Phi _n \left( a_i ^j \right)
$ to be the unique element of $\Phi \left[ A_i ^j \right] \cap U_n $.

\end{define}

\begin{lem}
\label{pruneorbit2}

For every $\Phi \in {\rm Aut} _{\cal D} (U)$, we have that
$\Phi _n \in {\rm Aut } _{{\cal D}_n} (U_n )$.
Moreover,
$(U_n , {\cal D}_n )$ is a tight interdependent
orbit union.

\end{lem}

{\bf Proof.}
Clearly,
${\cal D}_n$ is a dictated orbit structure for $U_n $.

Because $\{ D_1 , \ldots , D_{n-1} \} $
is a connected component of
the graph ${\cal O} (U, {\cal D})-D_{n} $,
by Lemma \ref{pruneorbit3},
$(U_n , {\cal D}_n )$ is a tight interdependent
orbit union.

Because $(U,{\cal D})$ is without slack and because we
choose exactly one element from each maximal
$U\setminus \bigcup _{i=n} ^m D_i $-order-autonomous antichain
to be in $U_n $, we obtain that $(U_n ,{\cal D}_n )$ is
an interdependent orbit union without slack.

For any
$\Phi \in {\rm Aut} _{\cal D} (U)$,
it follows from the definitions that
$\Phi _n \in {\rm Aut } _{{\cal D}_n} (U_n )$.
Because, for any $j\in \{ s, \ldots , n-1\} $ and
$x,y\in \{ 1, \ldots , \ell _j \} $, there is a
$\Phi \in {\rm Aut} _{\cal D} (U)$ with
$\Phi \left[ A_x ^j \right] =A_y ^j $, we conclude that
$(U_n , {\cal D}_n )$ is tight.
\qed

\vspace{.1in}

With the ``early orbits"
$D_1 , \ldots , D_{n-1} $ thus
discussed, we now turn to the remaining
orbits
$D_n , \ldots , D_m $ as well as the connection
between $D_n $ and the sets $A_i ^j $.

\begin{define}

Let $Q:=\bigcup _{j=s} ^{m} D_j $,
let ${\cal E}_Q
:=
\{ A_i ^j : j=s, \ldots , n-1; i=1, \ldots , \ell _j \}
\cup \{ D_j :j=n, \ldots ,m \} $
and
let ${\cal D}_Q $ be the set of all ${\cal E}_Q $-orbits.

\end{define}

\begin{lem}
\label{Qonlystruct}

$(Q , {\cal D}_Q )$ is a tight
structured ordered set.

\end{lem}

{\bf Proof.}
First note that,
for $j\geq n$, all $D_i $ that are directly interdependent with
$D_j $ are contained in $\bigcup _{k=s} ^m D_k = Q$.
Therefore,
any $Q$-order-autonomous antichain in a set $D_j $ with $j\geq n$
would be $U$-order-autonomous.
Hence
for $j\geq n$, no $D_j $ contains a nontrivial
$Q$-order-autonomous antichain.

By
Lemma \ref{pruneorbit3},
for $j\in \{ s, \ldots , n-1\} $ and
$i\in \{ 1, \ldots , \ell _j \} $,
no $A_i ^j $ contains a nontrivial
$Q$-order-autonomous antichain.
We conclude that $(Q,{\cal E}_Q )$ is
a structured ordered set
without slack.
Because
${\cal D}_Q $ is the set of all ${\cal E}_Q $-orbits, we conclude that
$(Q,{\cal D}_Q )$ is a
tight
structured ordered set.
\qed

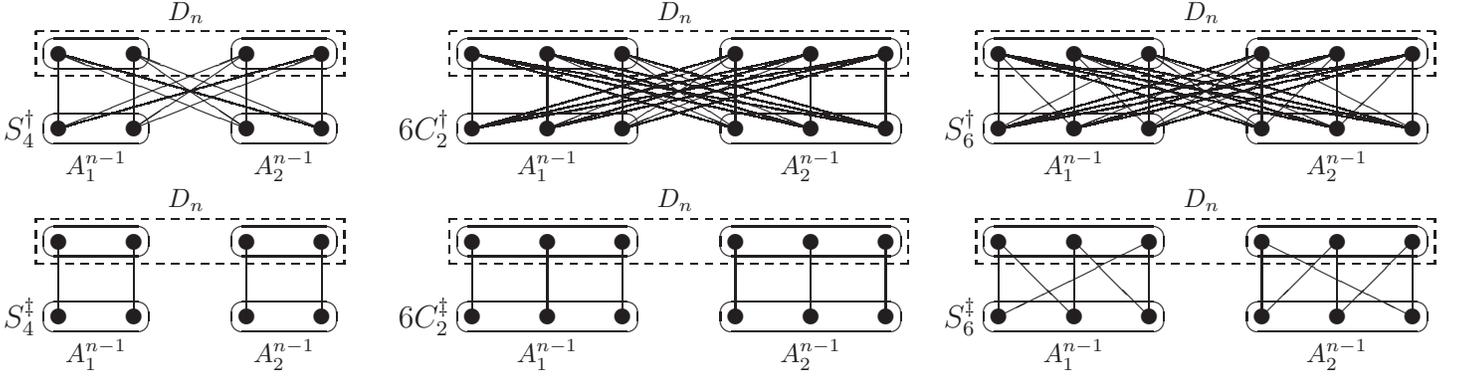
\begin{figure}

\centerline{
\unitlength 1mm 
\linethickness{0.4pt}
\ifx\plotpoint\undefined\newsavebox{\plotpoint}\fi 
\begin{picture}(188,44)(0,0)
\put(5,5){\line(0,1){10}}
\put(5,30){\line(0,1){10}}
\put(30,5){\line(0,1){10}}
\put(30,30){\line(0,1){10}}
\put(15,5){\line(0,1){10}}
\put(15,30){\line(0,1){10}}
\put(40,5){\line(0,1){10}}
\put(40,30){\line(0,1){10}}
\put(60,5){\line(0,1){10}}
\put(60,30){\line(0,1){10}}
\put(130,5){\line(0,1){10}}
\put(130,30){\line(0,1){10}}
\put(70,5){\line(0,1){10}}
\put(70,30){\line(0,1){10}}
\put(140,5){\line(0,1){10}}
\put(140,30){\line(0,1){10}}
\put(80,5){\line(0,1){10}}
\put(80,30){\line(0,1){10}}
\put(150,5){\line(0,1){10}}
\put(150,30){\line(0,1){10}}
\put(95,5){\line(0,1){10}}
\put(95,30){\line(0,1){10}}
\put(165,5){\line(0,1){10}}
\put(165,30){\line(0,1){10}}
\put(105,5){\line(0,1){10}}
\put(105,30){\line(0,1){10}}
\put(175,5){\line(0,1){10}}
\put(175,30){\line(0,1){10}}
\put(115,5){\line(0,1){10}}
\put(115,30){\line(0,1){10}}
\put(185,5){\line(0,1){10}}
\put(185,30){\line(0,1){10}}
\put(5,5){\circle*{2}}
\put(5,30){\circle*{2}}
\put(30,5){\circle*{2}}
\put(30,30){\circle*{2}}
\put(15,5){\circle*{2}}
\put(15,30){\circle*{2}}
\put(40,5){\circle*{2}}
\put(40,30){\circle*{2}}
\put(60,5){\circle*{2}}
\put(60,30){\circle*{2}}
\put(130,5){\circle*{2}}
\put(130,30){\circle*{2}}
\put(70,5){\circle*{2}}
\put(70,30){\circle*{2}}
\put(140,5){\circle*{2}}
\put(140,30){\circle*{2}}
\put(80,5){\circle*{2}}
\put(80,30){\circle*{2}}
\put(150,5){\circle*{2}}
\put(150,30){\circle*{2}}
\put(95,5){\circle*{2}}
\put(95,30){\circle*{2}}
\put(165,5){\circle*{2}}
\put(165,30){\circle*{2}}
\put(105,5){\circle*{2}}
\put(105,30){\circle*{2}}
\put(175,5){\circle*{2}}
\put(175,30){\circle*{2}}
\put(115,5){\circle*{2}}
\put(115,30){\circle*{2}}
\put(185,5){\circle*{2}}
\put(185,30){\circle*{2}}
\put(5,15){\circle*{2}}
\put(5,40){\circle*{2}}
\put(30,15){\circle*{2}}
\put(30,40){\circle*{2}}
\put(15,15){\circle*{2}}
\put(15,40){\circle*{2}}
\put(40,15){\circle*{2}}
\put(40,40){\circle*{2}}
\put(60,15){\circle*{2}}
\put(60,40){\circle*{2}}
\put(130,15){\circle*{2}}
\put(130,40){\circle*{2}}
\put(70,15){\circle*{2}}
\put(70,40){\circle*{2}}
\put(140,15){\circle*{2}}
\put(140,40){\circle*{2}}
\put(80,15){\circle*{2}}
\put(80,40){\circle*{2}}
\put(150,15){\circle*{2}}
\put(150,40){\circle*{2}}
\put(95,15){\circle*{2}}
\put(95,40){\circle*{2}}
\put(165,15){\circle*{2}}
\put(165,40){\circle*{2}}
\put(105,15){\circle*{2}}
\put(105,40){\circle*{2}}
\put(175,15){\circle*{2}}
\put(175,40){\circle*{2}}
\put(115,15){\circle*{2}}
\put(115,40){\circle*{2}}
\put(185,15){\circle*{2}}
\put(185,40){\circle*{2}}
\put(10,5){\oval(14,4)[]}
\put(10,30){\oval(14,4)[]}
\put(35,5){\oval(14,4)[]}
\put(35,30){\oval(14,4)[]}
\put(10,15){\oval(14,4)[]}
\put(10,40){\oval(14,4)[]}
\put(35,15){\oval(14,4)[]}
\put(35,40){\oval(14,4)[]}
\put(5,30){\line(5,2){25}}
\put(40,30){\line(-5,2){25}}
\put(30,40){\line(-3,-2){15}}
\put(15,40){\line(3,-2){15}}
\put(15,30){\line(5,2){25}}
\put(30,30){\line(-5,2){25}}
\multiput(40,40)(-.1178451178,-.0336700337){297}{\line(-1,0){.1178451178}}
\multiput(5,40)(.1178451178,-.0336700337){297}{\line(1,0){.1178451178}}
\put(70,5){\oval(24,4)[]}
\put(70,30){\oval(24,4)[]}
\put(140,5){\oval(24,4)[]}
\put(140,30){\oval(24,4)[]}
\put(105,5){\oval(24,4)[]}
\put(105,30){\oval(24,4)[]}
\put(175,5){\oval(24,4)[]}
\put(175,30){\oval(24,4)[]}
\put(70,15){\oval(24,4)[]}
\put(70,40){\oval(24,4)[]}
\put(140,15){\oval(24,4)[]}
\put(140,40){\oval(24,4)[]}
\put(105,15){\oval(24,4)[]}
\put(105,40){\oval(24,4)[]}
\put(175,15){\oval(24,4)[]}
\put(175,40){\oval(24,4)[]}
\multiput(60,30)(.1178451178,.0336700337){297}{\line(1,0){.1178451178}}
\multiput(130,30)(.1178451178,.0336700337){297}{\line(1,0){.1178451178}}
\multiput(115,30)(-.1178451178,.0336700337){297}{\line(-1,0){.1178451178}}
\multiput(185,30)(-.1178451178,.0336700337){297}{\line(-1,0){.1178451178}}
\put(95,40){\line(-5,-2){25}}
\put(165,40){\line(-5,-2){25}}
\put(80,40){\line(5,-2){25}}
\put(150,40){\line(5,-2){25}}
\multiput(70,30)(.1178451178,.0336700337){297}{\line(1,0){.1178451178}}
\multiput(140,30)(.1178451178,.0336700337){297}{\line(1,0){.1178451178}}
\multiput(105,30)(-.1178451178,.0336700337){297}{\line(-1,0){.1178451178}}
\multiput(175,30)(-.1178451178,.0336700337){297}{\line(-1,0){.1178451178}}
\put(105,40){\line(-5,-2){25}}
\put(175,40){\line(-5,-2){25}}
\put(70,40){\line(5,-2){25}}
\put(140,40){\line(5,-2){25}}
\multiput(80,30)(.1178451178,.0336700337){297}{\line(1,0){.1178451178}}
\multiput(150,30)(.1178451178,.0336700337){297}{\line(1,0){.1178451178}}
\multiput(95,30)(-.1178451178,.0336700337){297}{\line(-1,0){.1178451178}}
\multiput(165,30)(-.1178451178,.0336700337){297}{\line(-1,0){.1178451178}}
\multiput(115,40)(-.1515151515,-.0336700337){297}{\line(-1,0){.1515151515}}
\multiput(185,40)(-.1515151515,-.0336700337){297}{\line(-1,0){.1515151515}}
\multiput(60,40)(.1515151515,-.0336700337){297}{\line(1,0){.1515151515}}
\multiput(130,40)(.1515151515,-.0336700337){297}{\line(1,0){.1515151515}}
\multiput(105,40)(-.1515151515,-.0336700337){297}{\line(-1,0){.1515151515}}
\multiput(175,40)(-.1515151515,-.0336700337){297}{\line(-1,0){.1515151515}}
\multiput(70,40)(.1515151515,-.0336700337){297}{\line(1,0){.1515151515}}
\multiput(140,40)(.1515151515,-.0336700337){297}{\line(1,0){.1515151515}}
\multiput(60,30)(.1851851852,.0336700337){297}{\line(1,0){.1851851852}}
\multiput(130,30)(.1851851852,.0336700337){297}{\line(1,0){.1851851852}}
\multiput(115,30)(-.1851851852,.0336700337){297}{\line(-1,0){.1851851852}}
\multiput(185,30)(-.1851851852,.0336700337){297}{\line(-1,0){.1851851852}}
\put(80,30){\line(3,2){15}}
\put(150,30){\line(3,2){15}}
\put(95,30){\line(-3,2){15}}
\put(165,30){\line(-3,2){15}}
\put(2,12){\dashbox{1}(41,6)[cc]{}}
\put(2,37){\dashbox{1}(41,6)[cc]{}}
\put(57,12){\dashbox{1}(61,6)[cc]{}}
\put(57,37){\dashbox{1}(61,6)[cc]{}}
\put(127,12){\dashbox{1}(61,6)[cc]{}}
\put(127,37){\dashbox{1}(61,6)[cc]{}}
\put(22,19){\makebox(0,0)[cb]{\footnotesize $D_n $}}
\put(22,44){\makebox(0,0)[cb]{\footnotesize $D_n $}}
\put(87,19){\makebox(0,0)[cb]{\footnotesize $D_n $}}
\put(87,44){\makebox(0,0)[cb]{\footnotesize $D_n $}}
\put(157,19){\makebox(0,0)[cb]{\footnotesize $D_n $}}
\put(157,44){\makebox(0,0)[cb]{\footnotesize $D_n $}}
\put(10,2){\makebox(0,0)[ct]{\footnotesize $A_1 ^{n-1} $}}
\put(10,27){\makebox(0,0)[ct]{\footnotesize $A_1 ^{n-1} $}}
\put(70,2){\makebox(0,0)[ct]{\footnotesize $A_1 ^{n-1} $}}
\put(70,27){\makebox(0,0)[ct]{\footnotesize $A_1 ^{n-1} $}}
\put(140,2){\makebox(0,0)[ct]{\footnotesize $A_1 ^{n-1} $}}
\put(140,27){\makebox(0,0)[ct]{\footnotesize $A_1 ^{n-1} $}}
\put(35,2){\makebox(0,0)[ct]{\footnotesize $A_2 ^{n-1} $}}
\put(35,27){\makebox(0,0)[ct]{\footnotesize $A_2 ^{n-1} $}}
\put(105,2){\makebox(0,0)[ct]{\footnotesize $A_2 ^{n-1} $}}
\put(105,27){\makebox(0,0)[ct]{\footnotesize $A_2 ^{n-1} $}}
\put(175,2){\makebox(0,0)[ct]{\footnotesize $A_2 ^{n-1} $}}
\put(175,27){\makebox(0,0)[ct]{\footnotesize $A_2 ^{n-1} $}}
\put(150,15){\line(-2,-1){20}}
\put(150,40){\line(-2,-1){20}}
\put(165,15){\line(2,-1){20}}
\put(165,40){\line(2,-1){20}}
\put(130,15){\line(1,-1){10}}
\put(130,40){\line(1,-1){10}}
\put(185,15){\line(-1,-1){10}}
\put(185,40){\line(-1,-1){10}}
\put(140,15){\line(1,-1){10}}
\put(140,40){\line(1,-1){10}}
\put(175,15){\line(-1,-1){10}}
\put(175,40){\line(-1,-1){10}}
\put(2,5){\makebox(0,0)[rc]{$S_4 ^\ddag $}}
\put(2,30){\makebox(0,0)[rc]{$S_4 ^\dagger $}}
\put(127,5){\makebox(0,0)[rc]{$S_6 ^\ddag $}}
\put(127,30){\makebox(0,0)[rc]{$S_6 ^\dagger $}}
\put(57,5){\makebox(0,0)[rc]{$6C_2 ^\ddag $}}
\put(57,30){\makebox(0,0)[rc]{$6C_2 ^\dagger $}}
\end{picture}
}

\caption{
Structured ordered sets induced on $D_{n-1} \cup D_n $ with $n=m$
that are not interdependent orbit unions.
The orbits are marked with ovals,
$D_n $ is marked with a dashed box,
and the partition $A_1 ^{n-1} \cup A_2 ^{n-1} $
is labeled.
}
\label{inconvenient}

\end{figure}

\vspace{.1in}

Lemma \ref{Qonlystruct} may feel a little unsatisfying compared to
Lemma \ref{pruneorbit2} in that $(Q,{\cal D} _Q )$
need not be an interdependent orbit union.
However, Figure \ref{inconvenient} gives examples that this really need not be the case:
Any of the
structured ordered sets given there could be a
structured ordered set $(Q,{\cal D} _Q )$ in the case
in which $m=n$, that is,
$D_n $ is a pendant vertex, and $|U_n \cap D_{n-1} |=2$

\begin{define}

With
$\Phi _n $ as in Definition \ref{phindef},
we define
${\rm Aut} _{{\cal D} _Q } ^U (Q):=
\{ \Psi \in {\rm Aut} _{\cal D} (U): \Psi _n ={\rm id} _{U_n } \} $.
For every $\Delta \in
{\rm Aut } _{{\cal D}_Q } (Q)$, we define
$\Delta ^U $
by
$$\Delta ^U (x):=
\cases{
\Delta (x); & if $x\in Q$, \cr
x; & if $x\in U\setminus Q$.} $$

\end{define}

\begin{lem}
\label{AutDQUrepres}

${\rm Aut} _{{\cal D} _Q } ^U (Q)=
\{ \Psi \in {\rm Aut} _{\cal D} (U): \Psi _n ={\rm id} _{U_n } \}
=
\left\{ \Delta ^U :
\Delta \in {\rm Aut } _{{\cal D}_Q } (Q)
\right\} $.

\end{lem}

{\bf Proof.}
First note that,
for every $\Delta \in
{\rm Aut } _{{\cal D}_Q } (Q)$, the fact that
$\Delta $ fixes all sets $A_i ^j $ implies that
$\Delta ^U \in {\rm Aut } _{{\cal D}} (U)$ and that
$\Delta ^U _n ={\rm id} _{U_n } $.

Conversely, if
$\Psi \in {\rm Aut} _{\cal D} (U)$
satisfies $\Psi _n ={\rm id} _{U_n } $, then
$\Psi $ maps every $A_i ^j $ to itself.
Hence $\Psi |_Q \in {\rm Aut } _{{\cal D}_Q } (Q)$
and $\Psi =(\Psi |_Q )^U $.
\qed

\begin{theorem}
\label{pruneorbit6}

The set ${\rm Aut} _{{\cal D} _Q } ^U (Q)$ is a normal subgroup of
${\rm Aut} _{\cal D} (U)$ and
the factor group
${\rm Aut} _{\cal D} (U)/{\rm Aut} _{{\cal D} _Q } ^U (Q)$
is isomorphic to a subgroup of
${\rm Aut } _{{\cal D}_n} (U_n )$.
Consequently
$
|{\rm Aut } _{{\cal D}} (U)|
\leq
|{\rm Aut } _{{\cal D}_n} (U_n )|
|{\rm Aut } _{{\cal D}_Q } (Q)|$.

\end{theorem}

{\bf Proof.}
First note that, because
$\Phi _n \Psi _n =(\Phi \Psi )_n $ and
$(\Phi ^{-1} )_n =(\Phi _n )^{-1} $, we have that
${\rm Aut} _{{\cal D} _Q } ^U (Q)$ is a subgroup of
${\rm Aut} _{\cal D} (U)$.

Via Lemma \ref{AutDQUrepres}, let
$\Delta ^U \in {\rm Aut} _{{\cal D} _Q } ^U (Q)$
and let $\Phi \in {\rm Aut} _{\cal D} (U)$.
Then
$
\left( \Phi ^{-1} \Delta ^U \Phi \right) _n
=
\Phi ^{-1} _n \Delta ^U _n \Phi _n
=
\Phi ^{-1} _n \Phi _n
=
{\rm id} _{U_n } $.
Hence
$\Phi ^{-1} {\rm Aut} _{{\cal D} _Q } ^U (Q) \Phi
=
{\rm Aut} _{{\cal D} _Q } ^U (Q)$
and therefore ${\rm Aut} _{{\cal D} _Q } ^U (Q)$
is a normal subgroup of ${\rm Aut} _{\cal D} (U)$.

Moreover, for all
$\Phi, \Psi \in {\rm Aut} _{\cal D} (U)$, we have
$\Phi {\rm Aut} _{{\cal D} _Q } ^U (Q)=\Psi {\rm Aut} _{{\cal D} _Q } ^U (Q)$
iff $\Phi _n =\Psi _n $.
Hence the factor group
${\rm Aut} _{\cal D} (U)/{\rm Aut} _{{\cal D} _Q } ^U (Q)$
is, via $\Phi {\rm Aut} _{{\cal D} _Q } ^U (Q)\mapsto \Phi _n $,
isomorphic to a subgroup of
${\rm Aut } _{{\cal D}_n} (U_n )$.
[The subgroup may be proper, because it can happen that
not all $\Psi \in {\rm Aut } _{{\cal D}_n} (U_n )$ are equal to
a $\Phi _n $ with $\Phi \in {\rm Aut} _{\cal D} (U)$.]
The inequality now follows.
\qed

\begin{define}

The {\bf separation partition
${\cal S} (D_n )$
of $D_{n} $}
is the partition of $D_n $ that is contained in ${\cal D} _Q $.

\end{define}

\begin{lem}
\label{pruneorbit5}

Every $\Phi \in {\rm Aut} _{\cal D} (U)$ respects
${\cal D} _Q $.
Every nontrivial
${\cal D}_Q $-orbit
that is contained
in
a set $A_i ^j $
is directly interdependent
with a
${\cal D}_Q $-orbit
$S\in {\cal S} (D_n )$.
Every
${\cal D}_Q $-interdependent orbit union that is contained in
$\bigcup _{j=s} ^{n-1} D_j $
is a singleton.

\end{lem}

{\bf Proof.}
Suppose, for a contradiction, there are a
$\Phi \in {\rm Aut} _{\cal D} (U)$
and an $S\in {\cal D} _Q $ such that
$\Phi [S]\not\in {\cal D} _Q $.
Then $\Phi [S]$ intersects two distinct sets $B,C\in {\cal D}_Q $
or $\Phi [S]$ is strictly contained in a set $D\in {\cal D} _Q $.
Because we are free to work with the inverse, it suffices to consider
the case
in which
$\Phi [S]$ intersects two distinct sets $B,C\in {\cal D} _Q $.
Because $S\in {\cal D} _Q $, there is a
$\Delta ^U \in {\rm Aut} _{{\cal D}_Q } (Q)$
that maps a $b\in \Phi ^{-1} [B]\cap S$ to a $c\in \Phi ^{-1} [C]\cap S$.
Now $\Phi (b)\in B$,
$\Phi \Delta ^U \Phi ^{-1} \left( \Phi (b)\right) =\Phi (c)\in C$
and $\Phi \Delta ^U \Phi ^{-1} |_Q \in {\rm Aut} _{{\cal D}_Q } (Q)$,
a contradiction.
We thus conclude that
every $\Phi \in {\rm Aut} _{\cal D} (U)$ respects
${\cal D} _Q $.

Finally,
by
Lemma \ref{pruneorbit3},
any nontrivial
${\cal D}_Q $-orbit
that is contained
in
an $A_i ^j $
is directly interdependent
with an $S\in {\cal S} (D_n )$.
Because no two distinct $A_i ^j $ are directly interdependent and
no $A_i ^j $ is directly interdependent with a $D_k $ with $k>n$,
every
${\cal D}_Q $-interdependent orbit union that is contained in
$\bigcup _{j=s} ^{n-1} D_j $
must be a singleton.
\qed

\section{Induction Preparation: Estimates for
$(Q, {\cal D} _Q )$}
\label{indprep}

We now obtain the following
generalization of Lemma \ref{allbutoneac}.

\begin{lem}
\label{pruneorbit}

With notation as given so far, the following hold.
\begin{enumerate}
\item
\label{pruneorbit7}
If every $\Psi \in {\rm Aut } _{{\cal D}_Q } (Q)$ is uniquely determined by
$\Psi |_{D_n } $,
then
$
|{\rm Aut } _{{\cal D}_Q } (Q)|
\leq
\min \left\{
(w(U)-1)!,
\left( \left( {|D_n |\over |{\cal S} (D_n )|} \right) !\right)
^{|{\cal S} (D_n )|}
\right\}
$.

\item
\label{pruneorbit9}
If $n=m$, that is, if $D_n $ is not a cutvertex of ${\cal O} (U, {\cal D})$, then
$
|{\rm Aut } _{{\cal D}_Q } (Q)|
\leq
\min \left\{
(w(U)-1)!,
\left( \left( {|D_n |\over |{\cal S} (D_n )|} \right) !\right)
^{|{\cal S} (D_n )|}
,
\prod _{j=s} ^{n-1}
\prod _{i=1} ^{\ell _j } \left| A_i ^j \right| !
\right\}
$.

\item
\label{pruneorbit10}
If $n=m$, that is, if $D_n $ is not a cutvertex of ${\cal O} (U, {\cal D})$,
and if, for all $j\in \{ s, \ldots , m-1\} $ and $i\in \{ 1, \ldots , \ell _j \} $,
we have $|A_i ^j |=1$, then
$
|{\rm Aut } _{{\cal D}} (U)|
\leq
|{\rm Aut } _{{\cal D}_m} (U_m )|
$.

\end{enumerate}

\end{lem}

{\bf Proof.}
For
part \ref{pruneorbit7}, note that,
${\rm Aut } _{{\cal D}_Q } (Q)$
induces at most
$
\left( \left( {|D_n |\over |{\cal S} (D_n )|} \right) !\right)
^{|{\cal S} (D_n )|}
$
permutations on $D_n $.
By hypothesis for this part,
every $\Phi \in {\rm Aut } _{{\cal D}_Q } (Q)$ is determined by
its values on $D_n $, hence
$|{\rm Aut } _{{\cal D}_Q } (Q)|
\leq
\left( \left( {|D_n |\over |{\cal S} (D_n )|} \right) !\right)
^{|{\cal S} (D_n )|}
$.
Moreover, suppose, for a contradiction,
that
${\rm Aut } _{{\cal D}_Q } (Q)$
induces
$>(w(U)-1)!$
permutations on $D_n $.
Then, in case $w(U)\not= 4$,
$Q$ is max-locked.
By Lemma \ref{maxlockeddescr}
(or directly for $w(U)=4$), this
implies
that
all $|A_i ^j |$ are equal to $w(U)$ and then
$U$ would contain an antichain with
$2w(U)$ elements,
a contradiction.

For part \ref{pruneorbit9},
first note that,
by
Lemma \ref{pruneorbit3},
for all $\Phi \in {\rm Aut} _{{\cal D} _Q } (Q)$,
we have that
$\Phi |_{D_n } $ determines all $\Phi |_{A_i ^j } $,
which means we can apply
part \ref{pruneorbit7} to obtain the first two numbers in the minimum
on the right hand side.
For the last number in the minimum on the right hand side, note that, because any
two distinct $x,y\in D_m $ satisfy $\{ p\in P:p\sim x\} \not= \{ p\in P:p\sim y\} $
and all elements comparable to elements in $D_m $ are in
$\bigcup _{j=s} ^{n-1}
\bigcup _{i=1} ^{\ell _j } A_i ^j $,
we have that
$\Phi |_{
\bigcup _{j=s} ^{n-1}
\bigcup _{i=1} ^{\ell _j } A_i ^j
} $ determines $\Phi |_{D_m } $.

Part \ref{pruneorbit10} follows from part \ref{pruneorbit9}
and Theorem \ref{pruneorbit6}.
\qed

\vspace{.1in}

The idea for our induction proof is to add
$(Q,{\cal D} _Q )$
to $(U_n , {\cal D} _n )$
to obtain the interdependent orbit union
$(U, {\cal D})$.
We could (and will) work with
insights from the induction hypothesis regarding
$(Q,{\cal D} _Q )$ being $w(U)$-adequately bounded.
However,
one on hand, these insights may
not be available because $|{\cal D} _Q |$ is too large
for the induction hypothesis to apply.
On the other hand, such an estimate may be
more than what is needed:
$(U_n , {\cal D} _n )$ has
$\sum _{j=1} ^{s-1} |D_j |+\sum _{j=s} ^{n-1} \ell _j $ points
and $n-1$ orbits,
$(U, {\cal D})$ has
$\sum _{j=1} ^m |D_j |$ points and
$m$ orbits, and
$(Q,{\cal D} _Q )$ has
$\sum _{j=s} ^m |D_j |$ points,
but the number of orbits is, at least to this author,
not easy to determine.
Because we ultimately need an estimate for
$(U, {\cal D})$,
Definition \ref{adequateindbound}
below is the notion of boundedness we need for
$(Q,{\cal D} _Q )$, as
Lemma \ref{withandwithoutQ} will illustrate.

\begin{define}
\label{adequateindbound}

Let $q\in {\mat Z}$ and define
$\ell _Q :=\sum _{j=s} ^{n-1} \ell _j $ and
$$
T_q
:=
\left\lfloor
{1\over 2}
\left(
\left(
\sum _{j=n+1} ^m |D_j |
-(m-n)
\right)
+
\left( |D_n |-
{|{\cal S} (D_n )|}
\right)
+\left(
\sum _{j=s} ^{n-1}
|D_j|
-
\ell _Q \right)
+q
\right)
\right\rfloor
.$$

We will call
$(Q,{\cal D} _Q )$
{\bf induction $w$-adequately $q$-bounded} iff
there are numbers
$w_1 , \ldots , w_M \in \{ 0, \ldots , w-1\} $
such that
$|{\rm Aut } _{{\cal D}_Q } (Q)|
\leq
\prod _{j=1} ^M w_j !
$
and such that
$
\sum  _{j=1} ^M w_j
\leq
T_q $.
In cases in which
$q\leq |{\cal S} (D_n )|-1$
and
the specific value of $q$ is immaterial, we will
also say that $(Q,{\cal D} _Q )$ is {\bf induction $w$-adequately bounded}.

\end{define}

\begin{lem}
\label{withandwithoutQ}

With notation as given so far,
let
$(U_n , {\cal D} _n )$
be offset by $o$ from
being $w$-adequately bounded.
\begin{enumerate}
\item
\label{withoutQ}
If
$
|{\rm Aut } _{{\cal D}} (U)|
\leq
|{\rm Aut } _{{\cal D}_n} (U_n )|
$,
then $(U,{\cal D} )$ is
offset by $o-(|D_n |-1)$ from being
$w$-adequately bounded.

\item
\label{Qwadeq}
If $(Q, {\cal D} _Q )$ is
induction $w$-adequately $q$-bounded,
then
$(U, {\cal D} )$ is offset by $o+q-(|{\cal S} (D_n )| -1)$ from
being $w$-adequately bounded.

In particular, if
$(U_n , {\cal D} _n )$
is $w$-adequately bounded
and $(Q, {\cal D} _Q )$ is
induction $w$-adequately bounded,
then
$(U, {\cal D} )$ is
$w$-adequately bounded.

Moreover, if
$(U_n , {\cal D} _n )$
is offset by $1$ from being $w$-adequately bounded
and $(Q, {\cal D} _Q )$ is
induction $w$-adequately $(|{\cal S} (D_n )| -2)$-bounded,
then
$(U, {\cal D} )$ is
$w$-adequately
bounded.

\item
\label{Qwadeqoddiff}
If $o=1$, $(Q, {\cal D} _Q )$ is
induction
$w$-adequately $(|{\cal S} (D_n )|-1)$-bounded,
and
$\left(
\sum _{j=n+1} ^m |D_j |
-(m-n)
\right)
+
\left( |D_n |
-1
\right)
+\left(
\sum _{j=s} ^{n-1}
|D_j|
-
\ell _Q \right)
$
is odd,
then
$(U, {\cal D} )$ is
$w$-adequately bounded.

\end{enumerate}

\end{lem}

{\bf Proof.}
Because
$(U_n , {\cal D} _n )$
is offset by $o$ from
being
$w$-adequately bounded, there are numbers
$w_i \in \{ 0, \ldots , w-1\} $,
$i=1, \ldots , M$ such that
$|{\rm Aut } _{{\cal D}_n} (U_n )|
\leq
\prod _{i=1} ^{M} w_i !$
and
$
\sum  _{i=1} ^{M} w_i
\leq
\left\lfloor {1\over 2}
\left(
|U_n |-|{\cal D} _n | +o\right)
\right\rfloor $.
Now part \ref{withoutQ} follows easily, because
the size of the set increases by $\sum _{i=n} ^m |D_i | $ and the number of
orbits increases by $m-n+1$ as we go from
$(U_n , {\cal D} _n )$ to
$(U , {\cal D} )$.

For part \ref{Qwadeq},
because
$(Q, {\cal D} _Q )$ is
induction $w$-adequately $q$-bounded,
there are numbers $w_{M+1} , \ldots , w_{M+K} $ such that
$|{\rm Aut } _{{\cal D}_Q } (Q)|
\leq \prod _{i=M+1} ^{M+K} w_i !$ with
$w_i \in \{ 0,\ldots , w-1\} $ and
$\sum _{j=M+1} ^{M+K} w_j \leq T_q $.

By
Theorem \ref{pruneorbit6},
we have
$$
|{\rm Aut } _{{\cal D}} (U)|
\leq
|{\rm Aut } _{{\cal D}_n} (U_n )|
|{\rm Aut } _{{\cal D}_Q } (Q)|
\leq
\prod _{i=1} ^{M} w_i !
\prod _{i=M+1} ^{M+K} w_i !
=\prod _{i=1} ^{M+K} w_i ! ,$$

Because $|{\cal D} _n |=n-1$
and $|{\cal D}  |=m$,
the sum of the $w_i $ is bounded as follows.

\vspace*{-.2in}
\hspace*{-1in}
\begin{minipage}[t]{3in}
\begin{eqnarray*}
\lefteqn{\sum  _{i=1} ^{M+K} w_i
=
\sum  _{i=1} ^{M} w_i + \sum _{i=M+1 } ^{M+K} w_i
}
\\
& \leq &
\left\lfloor {1\over 2}
\left( |U_n | -|{\cal D} _n | +o\right)
\right\rfloor
+
\left\lfloor {1\over 2}
\left(
\sum_{i=n} ^m |D_i |
-(m-n+1)
+
\sum _{j=s} ^{n-1}
|D_j |
-
\ell _Q
+q-(|{\cal S} (D_n )|-1)
\right)
\right\rfloor
\\
& \leq &
\left\lfloor {1\over 2}
\left( |U_n |+
\sum _{j=s} ^{n-1}
|D_j |
-
\ell _Q
+\sum_{i=n} ^m |D_i |
-(m-n+1)-(n-1)
+o+q-(|{\cal S} (D_n )|-1)
\right)
\right\rfloor
\\
& = &
\left\lfloor {1\over 2}
\left( |U|
-|{\cal D}|
+o+q-(|{\cal S} (D_n )|-1)
\right)
\right\rfloor
\end{eqnarray*}
\end{minipage}

~\\

For part \ref{Qwadeqoddiff},
$
|{\rm Aut } _{{\cal D}} (U)|
$ is bounded by the same product.
Because,
for odd numbers $2k+1$, we have
$\left\lfloor {1\over 2} (2k+1)\right\rfloor
=
\left\lfloor {1\over 2} (2k+1-1)\right\rfloor $,
we obtain the following.

\vspace*{-.2in}
\hspace*{-1in}
\begin{minipage}[t]{3in}
\begin{eqnarray*}
\lefteqn{
\sum  _{i=1} ^{M+K} w_i
=
\sum  _{i=1} ^{M} w_i + \sum _{i=M+1 } ^{M+K} w_i
}
\\
& \leq &
\left\lfloor {1\over 2}
\left( |U_n | -|{\cal D} _n |+1 \right)
\right\rfloor
+
\left\lfloor {1\over 2}
\left(
\left(
\sum _{j=n+1} ^m |D_j |
-(m-n)
\right)
+
\left( |D_n |
-1
\right)
+\left(
\sum _{j=s} ^{n-1}
|D_j|
-
\ell _Q \right)
-1 \right)
\right\rfloor
\\
& \leq &
\left\lfloor {1\over 2}
\left( |U_n |
+
\left(
\sum_{i=n} ^m |D_i |
+
\sum _{j=s} ^{m-1}
|D_j |
-
\ell _Q
\right)
-|{\cal D} _n |
-
(m-n+1)
\right)
\right\rfloor
\\
& = &
\left\lfloor {1\over 2}
\left( |U | -|{\cal D} |\right)
\right\rfloor .
\end{eqnarray*}
\end{minipage}

\qed

We conclude this section with some sufficient conditions for
$(Q, {\cal D} _Q )$ being
induction $w(U)$-adequately bounded.

\begin{lem}
\label{nontrivgcd}

Let $(P , {\cal D} )$ be a
tight structured ordered set
and let $C,D\in {\cal D}$ be so that
$C\upharpoonleft _{\cal D} D$.
Then $\gcd (|C|,|D|)>1$.

\end{lem}

{\bf Proof.}
Because $(P , {\cal D} )$ is tight, the automorphism group acts transitively on
$C$ and on $D$. Hence any two elements of $C$ have the same number
$1\leq u<|D|$ of upper
bounds in $D$ and
any two elements of $D$ have the same number $1\leq \ell <|C|$ of lower
bounds in $C$.
Moreover,
$|C|u=|D|\ell $.

Suppose, for a contradiction, that
$\gcd (|C|,|D|)=1$.
Because $|C|$ is a factor of $|D|\ell $, we have that
$|C|$ must divide $\ell $,
contradicting
$\ell < |C|$.
\qed

\vspace{.1in}

Similar to Definition \ref{inducedonrankdef}, we define the following.

\begin{define}

Let $(P,{\cal D} )$ be a structured ordered set and let
$D\in {\cal D}$.
We define
$\alpha _{\cal D} (D):=|\{ \Phi |_D :\Phi \in {\rm Aut} _{\cal D} (P)\} |$.

\end{define}

The goal throughout the following proofs, in which $n=m$,
is to
bound
$|{\rm Aut } _{{\cal D}_Q } (Q)|$
with a product of factorials such that the
sum of the individual numbers is less than or equal to
$T_q =
{1\over 2} \left\lfloor
\left( |D_n |-
{|{\cal S} (D_n )|}
\right)
+\left(
\sum _{j=s} ^{n-1}
|D_j |
-
\ell _Q
\right)
+q
\right\rfloor
.$
In many cases, $q$ will be equal to zero.
Because $|D_j |=\sum _{i=1} ^{\ell _j } \left| A_i ^j \right| $,
we will feel free to switch between
the two quantities, and the same holds true for
$\ell _Q =\sum _{j=s} ^{n-1} \ell _j $.

\begin{lem}
\label{seppart1}

With notation as given so far,
let $D_n $ not be a cutvertex of ${\cal O} (U, {\cal D})$,
$\sum _{j=t} ^{n-1}
\ell _j \geq 3$
and $|{\cal S}(D_n )|=1$.
Then
$(Q, {\cal D} _Q )$ is
induction $w(U)$-adequately bounded.

\end{lem}

{\bf Proof.}
By Lemma \ref{pruneorbit3}, we have
$|{\rm Aut} _{{\cal D}_Q} (Q) |\leq
\alpha _{{\cal D}_Q } (D_n )$. By part \ref{pruneorbit9}
of Lemma \ref{pruneorbit}, we have
$|{\rm Aut} _{{\cal D}_Q} (Q) |\leq
\prod _{j=s} ^{n-1}
\prod _{i=1} ^{\ell _j } \left| A_i ^j \right| !$.
Because $|{\cal S} (D_n )|=1$ there is no
$\left| A_i ^j \right| $ that is equal to $1$.

{\em Case 1:
$
\alpha _{{\cal D}_Q } (D_n )
>
(|D_n| -1)!$.}
Clearly
$|{\rm Aut} _{{\cal D}_Q} (Q) |\leq
\alpha _{{\cal D}_Q } (D_n )
\leq
|D_n|!$.
By Lemma \ref{maxlockeddescr}
(or directly for $|D_n |=4$),
there are
a $j$ and distinct $i,k$ such that
$|A_i ^j |\geq |D_n |$
and
$|A_k ^j |\geq |D_n |$.
Hence $w(U)\geq w(Q)\geq 2|D_n |$.
Moreover, because $\sum _{j=t} ^{n-1}
\ell _j \geq 3$,
there is another set $A_{i'} ^{j'} $ not equal to either
of $A_i ^j $ and $A_k ^j $ such that
$|A_{i'} ^{j'} |>1$.
Now
$T_0
\geq
\left\lfloor {1\over 2} \left(
(|A_i ^j | + |A_k ^j |-2)+
(|D_n| -1)+(|A_{i'} ^{j'} |-1)
\right) \right\rfloor
\geq
\left\lfloor {1\over 2} \left(
3|D_n| -3 +1
\right) \right\rfloor
=
\left\lfloor {1\over 2} \left(
2|D_n| +(|D_n |-2)
\right) \right\rfloor
\geq |D_n |
$,
proves that
$(Q,{\cal D} _Q )$ is induction $w(U)$-adequately
bounded.

{\em Case 2:
$
\alpha _{{\cal D}_Q } (D_n )
\leq
(|D_n| -1)!$.}

{\em Case 2.1:
$
(|D_n| -1)\leq
\left(
\sum _{j=s} ^{n-1}
\sum _{i=1} ^{\ell _j } \left| A_i ^j \right|
-
\sum _{j=s} ^{n-1}
\ell _j \right) $.}
In this case $T_0 \leq |D_n |-1$.
Via
$|{\rm Aut} _{{\cal D} _Q } (Q)|\leq
\alpha _{{\cal D}_Q } (D_n )
\leq
(|D_n| -1)!$,
we obtain that
$(Q,{\cal D}_Q )$
is induction
$w(U)$-adequately bounded.

{\em Case 2.2:
$
\left(
\sum _{j=s} ^{n-1}
\sum _{i=1} ^{\ell _j } \left| A_i ^j \right|
-
\sum _{j=s} ^{n-1}
\ell _j \right) < (|D_n| -1)$, and
there is an $|A_i ^j |$ that is not equal to $2$, or,
if all $|A_i ^j |=2$, then
$\sum _{j=s} ^{n-1}
\ell _j >3$.}
Because all cases
in which
Lemma \ref{factprodcombine} cannot be applied are excluded,
we obtain
$|{\rm Aut} _{{\cal D} _Q } (Q)|\leq
\prod _{j=s} ^{n-1}
\prod _{i=1} ^{\ell _j } \left| A_i ^j \right| !
\leq
\left( \sum _{j=s} ^{n-1}
\sum _{i=1} ^{\ell _j } \left| A_i ^j \right|
-
\sum _{j=s} ^{n-1}
\ell _j \right) !
$.
Because
$
\left(
\sum _{j=s} ^{n-1}
\sum _{i=1} ^{\ell _j } \left| A_i ^j \right|
-
\sum _{j=s} ^{n-1}
\ell _j \right) $ is smaller than $(|D_n| -1)\leq w(U)-1$
and less than or equal to $T_0 $,
we obtain that
$(Q,{\cal D}_Q )$
is induction
$w(U)$-adequately bounded.

{\em Case 2.3:
$
\left(
\sum _{j=s} ^{n-1}
\sum _{i=1} ^{\ell _j } \left| A_i ^j \right|
-
\sum _{j=s} ^{n-1}
\ell _j \right) < (|D_n| -1)$, and
all $|A_i ^j |$ are equal to $2$ and
$\sum _{j=s} ^{n-1}
\ell _j =3$.}
In this case, we have
$|A_1 ^{n-1} |=|A_2 ^{n-1} |=|A_3 ^{n-1} |=2$,
$|D_{n-1} |=6$
and
$
\left(
\sum _{j=s} ^{n-1}
\sum _{i=1} ^{\ell _j } \left| A_i ^j \right|
-
\sum _{j=s} ^{n-1}
\ell _j \right)
=3
$. Consequently $|D_n |>4$.
By Lemma \ref{nontrivgcd},
the case $|D_n|=5$ does not occur,
Clearly,
$|{\rm Aut} _{{\cal D} _Q } (Q)|\leq
\prod _{j=s} ^{n-1}
\prod _{i=1} ^{\ell _j } \left| A_i ^j \right| !
=8<4!$, and, because $|D_n |\geq 6$, we have
$T_0 \geq 4$, so
$(Q,{\cal D}_Q )$
is induction
$w(U)$-adequately bounded.
\qed

\begin{lem}
\label{boundsforQ}

With notation as given so far,
let $D_n $ not be a cutvertex of ${\cal O} (U, {\cal D})$
and let
$\sum _{j=t} ^{n-1}
\ell _j \geq 3$.
Then
$(Q, {\cal D} _Q )$ is
induction $w(U)$-adequately bounded.

\end{lem}

{\bf Proof.}
By Lemma \ref{seppart1}, we can assume that $|{\cal S} (D_n )|>1$.
In case
$
T_0 \geq w(U) - 1$,
using $|{\rm Aut } _{{\cal D}_Q } (Q)|
\leq
(w(U)-1)!$, we have that
$(Q,{\cal D} _Q )$ is
induction $w(U)$-adequately bounded.
Hence, for the remainder of this proof, we
can assume that
$
T_0 < w(U)-1$.

Recall that, by part \ref{pruneorbit9} of Lemma \ref{pruneorbit},
we have
$$
|{\rm Aut } _{{\cal D}_Q } (Q)|
\leq
\min
\left\{
(w(U)-1)!,
\left( \left( {|D_n |\over |{\cal S} (D_n )|} \right) !\right)
^{|{\cal S} (D_n )|}
,
\prod _{j=s} ^{n-1}
\prod _{i=1} ^{\ell _j } \left| A_i ^j \right| !\right\}
.$$
Further, recall that $|A_i ^j |>1$ iff $j\in \{ t, \ldots , n-1\} $.
Hence $\sum _{j=s} ^{t-1}
\sum _{i=1} ^{\ell _j } \left| A_i ^j \right|
-
\sum _{j=s} ^{t-1}
\ell _j =0$ and $\prod _{j=s} ^{t-1}
\prod _{i=1} ^{\ell _j } \left| A_i ^j \right| !=1
$, that is, this sum and product have no effect on the argument.
To stay connected to the
definition of
induction $w$-adequate boundedness, we will continue to
start all sums and products at $s$, but the estimates
are all derived and needed for the corresponding sums
and products starting at $t$.

{\em Case 1:
$\sum _{j=t} ^{n-1}
\ell _j \geq 4$ or
$\sum _{j=t} ^{n-1}
\ell _j =3$ and one of the
$\left| A_i ^{n-1} \right| $ is not equal to $2$.}

{\em Case 1.1:
$\left(
\sum _{j=s} ^{n-1}
\sum _{i=1} ^{\ell _j } \left| A_i ^j \right|
-
\sum _{j=s} ^{n-1}
\ell _j \right)
\leq
\left( |D_n |- {|{\cal S} (D_n )|}
\right)
$.}
In case $\sum _{j=t} ^{n-1}
\ell _j =3$,
we have $t=n-1$ and
then $|A_1 ^{n-1} |=|A_2 ^{n-1} |=|A_3 ^{n-1} |\not= 2$.
By Lemma \ref{factprodcombine},
we obtain
$
\prod _{j=s} ^{n-1}
\prod _{i=1} ^{\ell _j } \left| A_i ^j \right| !
\leq
\left(
\sum _{j=s} ^{n-1}
\sum _{i=1} ^{\ell _j } \left| A_i ^j \right|
-
\sum _{j=s} ^{n-1}
\ell _j \right) !
$.
Now
$
|{\rm Aut } _{{\cal D}_Q } (Q)|
\leq
\prod _{j=s} ^{n-1}
\prod _{i=1} ^{\ell _j } \left| A_i ^j \right| !
\leq
\left(
\sum _{j=s} ^{n-1}
\sum _{i=1} ^{\ell _j } \left| A_i ^j \right|
-
\sum _{j=s} ^{n-1}
\ell _j \right) !
$,
and, by assumption for this case,
$\left(
\sum _{j=s} ^{n-1}
\sum _{i=1} ^{\ell _j } \left| A_i ^j \right|
-
\sum _{j=s} ^{n-1}
\ell _j \right) $ is less than or equal to $T_0 $
as well as less than or equal to
$\left( |D_n |- {|{\cal S} (D_n )|}
\right) \leq w(Q)-1$.
Hence $(Q, {\cal D} _Q )$ is
induction $w(Q)\leq w(U)$-adequately bounded.

In the remaining subcases, 1.2 and 1.3,
because
$
\left( \left( {|D_n |\over |{\cal S} (D_n )|} \right) !\right)
^{|{\cal S} (D_n )|}
\leq
\left( |D_n |-{|{\cal S} (D_n )|}
+1\right) !$, we set
$H:=\min \left\{ (w(U)-1), \left( |D_n |-{|{\cal S} (D_n )|}
+1\right)  \right\} $, and we have
the inequality
$
|{\rm Aut } _{{\cal D}_Q } (Q)|
\leq
\min \left\{ (w(U)-1)!, \left( |D_n |-{|{\cal S} (D_n )|}
+1\right) ! \right\}
=H!$.

{\em Case 1.2:
$\left(
\sum _{j=s} ^{n-1}
\sum _{i=1} ^{\ell _j } \left| A_i ^j \right|
-
\sum _{j=s} ^{n-1}
\ell _j \right)
\geq \left( |D_n |-{|{\cal S} (D_n )|}
+2\right)
$.}
Because
$
H\leq
\left( |D_n |-{|{\cal S} (D_n )|}
+1\right)
\leq
T_0 $,
we obtain that
$(Q, {\cal D} _Q )$ is induction
$w(U)$-adequately bounded.

{\em Case 1.3:
$\left(
\sum _{j=s} ^{n-1}
\sum _{i=1} ^{\ell _j } \left| A_i ^j \right|
-
\sum _{j=s} ^{n-1}
\ell _j \right)
= \left( |D_n |-{|{\cal S} (D_n )|}
+1\right)
$.}
Because
$
|{\rm Aut } _{{\cal D}_Q } (Q)|
\leq
H!$
and $H\leq T_1 $,
we have that
$(Q, {\cal D} _Q )$ is
induction $w(U)$-adequately $1$-bounded.
Because $|{\cal S} (D_n )|>1$, this means that
$(Q, {\cal D} _Q )$ is
induction
$w(U)$-adequately bounded.

{\em Case 2:
$|A_1 ^{n-1} |=|A_2 ^{n-1} |=|A_3 ^{n-1} |=2$.}
In this case, we have
$|D_{n-1} |=6$.
and
$
\left(
\sum _{j=s} ^{n-1}
\sum _{i=1} ^{\ell _j } \left| A_i ^j \right|
-
\sum _{j=s} ^{n-1}
\ell _j \right)
=3
$.

We first consider small sizes $|D_n |\leq 5$.
Because
$|{\cal S} (D_n )|>1$, the cases $|D_n |=3$ and $|D_n |=5$
do not occur, and the case
$|D_n |=2$ leads to $|{\rm Aut} _{{\cal D}_Q } (Q)|=1$,
which is trivial.
This leaves the case
$|D_n |=4$, and there the only nontrivial case is
$|{\cal S} (D_n )|=2$.
Now
$|{\rm Aut} _{{\cal D} _Q } (Q)|\leq
\alpha _{{\cal D}_Q } (D_n )=
4<3!$.
Because $T_1 =3$,
we conclude that $(Q, {\cal D} _Q )$ is
induction $w(U)$-adequately bounded.

Finally, we consider $|D_n |\geq 6$.
First note that
$|{\rm Aut } _{{\cal D}_Q } (Q)|
\leq
\prod _{i=1} ^3 \left| A_i ^{n-1} \right|
=
(2!)^3 =8<4!$
and $4< w(Q)-1\leq w(U)-1$.
In case
$\left( |D_n |-{|{\cal S} (D_n )|}
\right) \geq 5$, we have $T_0 \geq 4$ and
$(Q, {\cal D} _Q )$ is
induction
$w(U)$-adequately bounded.
In case
$\left( |D_n |-{|{\cal S} (D_n )|}
\right) \leq 4$,
because the sets in
${\cal S} (D_n )$ all are of the same size,
we obtain $|D_n |\in \{ 6,8\} $.

For $|D_n |=8$, we obtain $|{\cal S} (D_n ) |=4$,
and,
for $|D_n |=6$, by Lemma \ref{nontrivgcd},
we obtain $|{\cal S} (D_n ) |=3$.
In either case, we have that $T_2 =4$, which proves that
$(Q, {\cal D} _Q )$ is
induction $w(U)$-adequately bounded.
\qed

\vspace{.1in}

By Lemma \ref{pruneorbit1}, the case
$\sum _{j=t} ^{n-1}
\ell _j = 1$ does not occur.
This leaves the case
$\sum _{j=t} ^{n-1}
\ell _j = 2$.

\begin{lem}
\label{seppart3ormore}

With notation as given so far,
let $D_n $ not be a cutvertex of ${\cal O} (U, {\cal D})$,
$\sum _{j=t} ^{n-1}
\ell _j = 2$
and $|{\cal S}(D_n )|\geq 3$.
Then
$(Q, {\cal D} _Q )$ is
induction $w(U)$-adequately bounded.

\end{lem}

{\bf Proof.}
Because $\sum _{j=t} ^{n-1}
\ell _j = 2$, we have $\ell _{n-1} =2$, that is,
${\cal A} ^{n-1} =\{ A_1 ^{n-1} , A_2 ^{n-1} \} $.

{\em Case 1:
$
\left| A_1 ^{n-1} \right|
+
\left| A_2 ^{n-1} \right|
-
2
\geq
\left( |D_n |- {|{\cal S} (D_n )|}+1
\right)
$.}
In this case,
$|{\rm Aut } _{{\cal D}_Q } (Q)|
\leq
\alpha _{\cal D} (D_n )
\leq
\left( \left( {|D_n |\over |{\cal S} (D_n )|} \right) !\right)
^{|{\cal S} (D_n )|}
\leq
(|D_n |-{|{\cal S} (D_n )|}+1)!$
and
$|D_n |-{|{\cal S} (D_n )|}+1$ is less than or equal to
$w(Q)-2\leq w(U)-1$
as well as less than or equal to
$T_1 $.
Hence
$(Q,{\cal D}_Q )$
is induction
$w(U)$-adequately bounded.

{\em Case 2:
$
\left| A_1 ^{n-1} \right|
+
\left| A_2 ^{n-1} \right|
-
2
\leq
\left( |D_n |- {|{\cal S} (D_n )|}
\right)
$.}
In this case,
$|{\rm Aut } _{{\cal D}_Q } (Q)|
\leq
\left| A_1 ^{n-1} \right| !
\left| A_2 ^{n-1} \right| !
\leq
\left(
\left| A_1 ^{n-1} \right|+
\left| A_2 ^{n-1} \right| -1\right) !
$
and
$\left| A_1 ^{n-1} \right|+
\left| A_2 ^{n-1} \right| -1
\leq T_0 +1\leq T_2 $.
Hence
$(Q,{\cal D}_Q )$
is induction
$w(U)$-adequately bounded.
\qed

\vspace{.1in}
The only cases for which we have not yet proved that
$(Q,{\cal D} _Q )$ is induction $w(U)$-adequately bounded
are $\sum _{j=t} ^{n-1}
\ell _j =2=\ell _{n-1} $ and
$|{\cal S}(D_n )|\leq 2$.
The
case
$\sum _{j=t} ^{n-1}
\ell _j = 2=\ell _{n-1} $,
$|{\cal S}(D_n )|=1$, will
follow from the base step for $|{\cal D} |=3$
in Section \ref{forbconfsec}.
This leaves the case
$\sum _{j=t} ^{n-1}
\ell _j = 2=\ell _{n-1} =
|{\cal S}(D_n )|$.

\begin{lem}
\label{seppart2and2Aij}

With notation as given so far,
let $D_n $ not be a cutvertex of ${\cal O} (U, {\cal D})$,
$\sum _{j=t} ^{n-1}
\ell _j = 2 =
|{\cal S}(D_n )|$ and $|D_n |\not= |D_{n-1} |$.
Then
$(Q, {\cal D} _Q )$ is
induction $w(U)$-adequately bounded.

\end{lem}

{\bf Proof.}
Without loss of generality, assume that
$|D_n |<|D_{n-1} |$.
Because $\sum _{j=t} ^{n-1}
\ell _j = 2 =
|{\cal S}(D_n )|$, we obtain
$|D_n |\leq |D_{n-1} |-2$.
Now
$|{\rm Aut } _{{\cal D}_Q } (Q)|
\leq
\left( {D_n \over 2} \right) !
\left( {D_n \over 2} \right) !
\leq
(|D_n |-1)!$
and $|D_n |-1
=
\left\lfloor {1\over 2} \left(
|D_n |-2+|D_{n} |
\right) \right\rfloor
\leq
\left\lfloor {1\over 2} \left(
|D_n |-2+|D_{n-1} |-2
\right) \right\rfloor
=
T_0 $.
Thus
$(Q, {\cal D} _Q )$ is
induction $w(U)$-adequately bounded.
\qed

\begin{lem}
\label{seppart2and2Aijequalsize}

With notation as given so far,
let $D_n $ not be a cutvertex of ${\cal O} (U, {\cal D})$,
$\sum _{j=t} ^{n-1}
\ell _j = 2 =
|{\cal S}(D_n )|$ and $|D_n |= |D_{n-1} |\geq 8$.
Then
$(Q, {\cal D} _Q )$ is
induction $w(U)$-adequately bounded.

\end{lem}

{\bf Proof.}
By Lemma \ref{factprodcombine},
$|{\rm Aut } _{{\cal D}_Q } (Q)|
\leq
\left( {D_n \over 2} \right) !
\left( {D_n \over 2} \right) !
\leq
(|D_n |-2)!$
and $|D_n |-2
=
\left\lfloor {1\over 2} \left(
|D_n |-2+|D_{n-1} |-2
\right) \right\rfloor
=
T_0 $.
Thus
$(Q, {\cal D} _Q )$ is
induction $w(U)$-adequately bounded.
\qed

\begin{lem}
\label{inconvenientSDn2}

With notation as given so far,
let $D_n $ not be a cutvertex of ${\cal O} (U, {\cal D})$,
$\sum _{j=t} ^{n-1}
\ell _j \geq 2 $,
$|{\cal S}(D_n )|\geq 2$,
let $(Q, {\cal D} _Q )$ be not
induction $w(U)$-adequately bounded
and let
$(V,{\cal D} _V )$ be the structured ordered subset that
consists of all the non-singleton
interdependent orbit unions of $(Q, {\cal D} _Q )$.
Then
$\ell_{n-1} =|{\cal S}(D_n )|= 2$,
$V=D_{n-1} \cup D_n $,
and $(V,{\cal D} _V )$
is isomorphic or dually isomorphic to
one of
the structured ordered sets
$\left( S_4 ^\dagger , {\cal D} _{S_4 ^\dagger } \right) $,
$\left( 6C_2 ^\dagger , {\cal D} _{6C_2 ^\dagger } \right) $,
$\left( S_6 ^\dagger , {\cal D} _{S_6 ^\dagger } \right) $,
$\left( S_4 ^\ddag , {\cal D} _{S_4 ^\ddag } \right) $,
$\left( 6C_2 ^\ddag , {\cal D} _{6C_2 ^\ddag } \right) $,
$\left( S_6 ^\ddag , {\cal D} _{S_6 ^\ddag } \right) $
in Figure \ref{inconvenient}.
In particular, $V=D_{n-1} \cup D_n $ does not contain any order-autonomous antichains.

\end{lem}

{\bf Proof.}
Let
$(Q, {\cal D} _Q )$ be not
induction $w(U)$-adequately bounded
such that
$\sum _{j=t} ^{n-1}
\ell _j \geq 2 $ and
$|{\cal S}(D_n )|\geq 2$.
By Lemma \ref{boundsforQ}, we have
$\sum _{j=t} ^{n-1}
\ell _j = 2 $, that is, ${\cal A} ^{n-1} =\{ A_1 ^{n-1} , A_2 ^{n-1} \} $.
By Lemma \ref{seppart3ormore}, we have
$|{\cal S}(D_n )|= 2$, that is, ${\cal S} (D_n )=\left\{ B_1 ^n , B_2 ^n \right\} $,
and we also conclude
$(V,{\cal D} _V )=( D_{n-1} \cup D_n ,
\{ A_1 ^{n-1} , A_2 ^{n-1} , B_1 ^n , B_2 ^n \} )$.
By Lemmas \ref{seppart2and2Aij} and \ref{seppart2and2Aijequalsize},
we have
$|D_{n-1} |=|D_n |<8$
and hence, because $|{\cal S}(D_n )|= 2$
and ${\cal S}(D_n )$ must be nontrivial,
$|D_{n-1} |=|D_n |\in \{ 4,6\} $.

Because no nontrivial ${\cal D} _Q $-orbit is order-autonomous in $Q$,
without loss of generality, we can assume that
$A_1 ^{n-1} \updownharpoons _{\cal D} B_1 ^n $
and
$A_2 ^{n-1} \updownharpoons _{\cal D} B_2 ^n $.
The only possibilities for
two orbits with
$2$ or $3$ elements
each to be directly interdependent
are isomorphic to $2C_2 $, $3C_2 $ and the six crown $C_6 =S_3 $.
Hence, each
$A_i ^{n-1} \cup B_i ^n $
must be isomorphic to one of these sets.
Because $(U, {\cal D} )$ is tight and, by Lemma \ref{pruneorbit5},
the ${\cal D} $-automorphisms respect ${\cal S} (D_n )$,
$A_1 ^{n-1} \cup B_1 ^n $
must be isomorphic to
$A_2 ^{n-1} \cup B_2 ^n $

Any further direct interdependence
$A_i ^{n-1} \updownharpoons _{\cal D} B_j ^n $
with $i\not= j$
would induce
$|{\rm Aut} _{{\cal D} _Q } (Q)|\leq 2!$
in case
$|D_{n-1} |=|D_n |=4$, or,
$|{\rm Aut} _{{\cal D} _Q } (Q)|\leq 3!$
in case
$|D_{n-1} |=|D_n |=6$.
In either case,
$(Q, {\cal D} _Q )$ would be
induction $w(U)$-adequately bounded.
Thus there are no further direct interdependences
in $(V,{\cal D} _V )$ and
all possible combinations are listed
in Figure \ref{inconvenient}.
\qed

\section{Base Step: Forbidden Configurations}
\label{forbconfsec}

Clearly, Lemma \ref{withandwithoutQ} is the key building block
to
show, via an induction on $|U|$, that ``many"
flexible tight interdependent
orbit unions
$(U, {\cal D})$ are
$w(U)$-adequately bounded.
Unlike in most induction proofs, for this particular situation,
the base step is a bit complex: There
are a number of
flexible tight interdependent
orbit unions
$(U,{\cal D})$ with two or three orbits
which are not $w(U)$-adequately bounded.
On the positive side, the proofs in this section will
show
how these ``forbidden configurations," see Definition \ref{forbconfdef}
below for details, arise naturally,
though annoyingly, via
max-locked ordered sets
for which Lemma \ref{halfdictateslemmaxlock}
does not provide the requisite bound, and
through the
merging of an $8$-crown or an ordered set
$S_4 $ or $4C_2 $
or $2V$ (see Definition \ref{forbconffirst} below)
with
the ordered set
$2V$.
The $8$-crown arises naturally though Lemma \ref{halfdictates2clique}
below.

\begin{lem}
\label{halfdictates2clique}

Let $(U, \{ B,T\} )$ be a flexible tight interdependent orbit union
such that
the elements of $B$ are minimal and
the elements of $T $ are maximal.
If
$U$ is not an 8-crown, then
$(U, \{ B,T\} )$ is $\max \{ |B|,|T|\} $-adequately bounded.

\end{lem}

{\bf Proof.}
Without loss of generality, assume that $|B|\leq |T|$.
Because $U$ is without slack,
by Lemma \ref{allbutoneac},
$|{\rm Aut } _{\{ B,T\} } (U)|
\leq
\alpha _{\{ B,T\} } (B)$.

Because
$(U,\{ B,T\} )$ is tight,
by Lemma \ref{nontrivgcd}, we have $\gcd (|B|,|T|)>1$.
Therefore, in
case $|B|<|T|$, we have
$|B|\leq |T|-2$.
Now
$|{\rm Aut } _{\{ B,T\} } (U)|
\leq
\alpha _{\{ B,T\} } (B)
\leq
|B|!$ and
$|B|
=
\left\lfloor {1\over 2} \left(
|B |+ |B|\right) \right\rfloor
\leq
\left\lfloor {1\over 2} \left(
|B |+ |T|-2\right) \right\rfloor $.

This leaves the
case $|B|=|T|$.
First suppose, for a contradiction, that
$|B|=|T|=4$.
Because
$(U,\{ B,T\} )$ is tight, and because
every element of $T$ has at least one lower cover, there is a $c\in \{ 1,2,3\} $
such that every element of $T$ has exactly $c$ lower covers and
every element of $B$ has exactly $c$ upper covers.
If every element of $T$ had exactly one lower cover,
then
$U$ would be a copy of $4C_2 $, contradicting that
$(U,\{ B,T\} )$ is flexible.
If every element of $T$ had exactly three lower covers,
then,
because $(U,\{ B,T\} )$ is without slack,
$U$ would be a copy of $S_4 $, again contradicting that
$(U,\{ B,T\} )$ is flexible.
Thus
every element of $T$ has two lower covers and
every element of $B$ has two upper covers.
Because $U$
is not an $8$-crown,
$U$ could only
be the disjoint union of two $4$-crowns,
contradicting that $(U,\{ B,T\} )$ is without slack.
Thus the case $|B|=|T|=4$
does not occur.

This leaves the case $|B|=|T|\not= 4$.
Because $U$ is flexible,
we have
$\alpha _{\{ B,T\} } (B)
\leq
(|B|-1)!$.
Hence
$|{\rm Aut } _{\{ B,T\} } (U)|
\leq
(|B|-1)!$ and
$|B|-1=
\left\lfloor {1\over 2} \left(
|B |+ |T|-2\right) \right\rfloor $.
\qed

\vspace{.1in}

The $8$-crown $C_8 $
has $8$ automorphisms, $8$ elements and $2$ natural orbits, so it is
offset by $4$ from being $w(C_8 )$-adequately bounded.
By Lemma \ref{halfdictates2clique},
the $8$-crown $C_8 $
with its natural dictated orbit structure
${\cal D} (C_8 )$
is
the only
flexible tight interdependent orbit union
$(U, {\cal D} )$
with exactly $2$ dictated orbits
that is not
$w(U)$-adequately bounded.
For exactly three orbits, we need to define a few ordered sets.

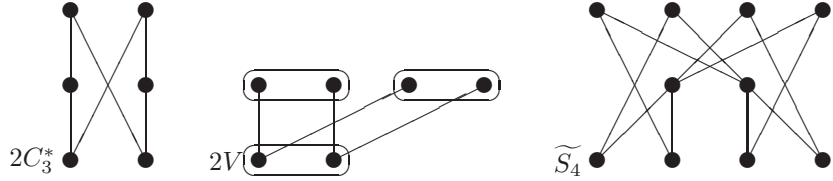
\begin{figure}

\centerline{
\unitlength 1mm 
\linethickness{0.4pt}
\ifx\plotpoint\undefined\newsavebox{\plotpoint}\fi 
\begin{picture}(111,26)(0,0)
\put(10,5){\line(0,1){10}}
\put(20,5){\line(0,1){10}}
\put(10,15){\line(0,1){10}}
\put(20,15){\line(0,1){10}}
\put(10,25){\line(1,-2){10}}
\put(20,25){\line(-1,-2){10}}
\put(10,5){\circle*{2}}
\put(35,5){\circle*{2}}
\put(45,5){\circle*{2}}
\put(20,5){\circle*{2}}
\put(10,15){\circle*{2}}
\put(35,15){\circle*{2}}
\put(80,5){\circle*{2}}
\put(80,25){\circle*{2}}
\put(45,15){\circle*{2}}
\put(90,5){\circle*{2}}
\put(90,25){\circle*{2}}
\put(20,15){\circle*{2}}
\put(10,25){\circle*{2}}
\put(20,25){\circle*{2}}
\put(8,5){\makebox(0,0)[rc]{\footnotesize $2C_3 ^* $}}
\put(35,15){\line(0,-1){10}}
\put(45,15){\line(0,-1){10}}
\put(35,5){\line(2,1){20}}
\put(45,5){\line(2,1){20}}
\put(55,15){\circle*{2}}
\put(100,5){\circle*{2}}
\put(90,15){\circle*{2}}
\put(100,25){\circle*{2}}
\put(65,15){\circle*{2}}
\put(110,5){\circle*{2}}
\put(100,15){\circle*{2}}
\put(110,25){\circle*{2}}
\put(33,5){\makebox(0,0)[rc]{\footnotesize $2V$}}
\put(78,5){\makebox(0,0)[rc]{\footnotesize $\widetilde{S_4 } $}}
\put(80,5){\line(1,1){10}}
\put(110,5){\line(-1,1){10}}
\put(90,15){\line(0,-1){10}}
\put(100,15){\line(0,-1){10}}
\put(100,25){\line(-1,-1){10}}
\put(90,25){\line(1,-1){10}}
\put(90,15){\line(2,1){20}}
\put(100,15){\line(-2,1){20}}
\put(80,5){\line(1,2){10}}
\put(110,5){\line(-1,2){10}}
\put(90,5){\line(-1,2){10}}
\put(100,5){\line(1,2){10}}
\put(40,15){\oval(14,4)[]}
\put(60,15){\oval(14,4)[]}
\put(40,5){\oval(14,4)[]}
\end{picture}
}

\caption{The forbidden configurations from Definition \ref{forbconffirst}.
Dictated orbit structures that are not equal to the natural one
are indicated with ovals.}
\label{forbconf_sporadic}

\end{figure}

\begin{define}
\label{forbconffirst}

(See Figure \ref{forbconf_sporadic}.)
We define the following ordered sets with dictated orbit structures.
\begin{enumerate}
\item
\label{forbconffirst1}
We define the ordered set
$2C_3 ^* $ to be the set
$\{ b_1 , b_2 , m_1 , m_2 , t_1 , t_2 \} $
with the comparabilities
$b_1 <m_1 <t_1 $,
$b_2 <m_2 <t_2 $, and
$b_1 , b_2 < t_1 , t_2 $.
We define
the dictated orbit structure ${\cal D} \left( 2C_3 ^* \right)
:=\big\{ \{ b_1 , b_2 \} ,
\{ m_1 , m_2 \} ,
\{ t_1 , t_2 \} \big\} $.

\item
We define the ordered set $2V$ to be the set
$\left\{ b_1 , b_2 , t_1 ^1 , t_2 ^1 , t_1 ^2 , t_2 ^2 \right\} $
with the comparabilities
$b_1 < t_1 ^1 , t_2 ^1 $ and
$b_2 < t_1 ^2 , t_2 ^2 $.
We define the dictated orbit structure
${\cal D} (2V):=
\left\{ \{ b_1 , b_2 \} , \left\{ t_1 ^1 , t_1 ^2 \right\} ,
\left\{ t_2 ^1 , t_2 ^2 \right\} \right\} $.

\item
We define the ordered set
$\widetilde{S_4 } $ to be the set
$\{ b_1 , b_2 , b_3 , b_4 , m_1 , m_2 , t_1 , t_2 , t_3 , t_4 \} $
with the comparabilities
$b_1 , b_2 <m_2 < t_3 , t_4 $,
$b_3 , b_4 <m_1 < t_1 , t_2 $,
$b_1 <t_2$,
$b_2 <t_1$,
$b_3 <t_4$,
$b_4 <t_3$,
plus the comparabilities dictated by transitivity.
Moreover, we define
the dictated orbit structure ${\cal D} \left( \widetilde{S_4 } \right)
:=\big\{ \{ b_1 , b_2 , b_3 , b_4 \} ,
\{ m_1 , m_2 \} ,
\{ t_1 , t_2 , t_3 , t_4 \} \big\} $.

\end{enumerate}

\end{define}

Next, we define a way to obtain ordered sets from smaller ordered sets by
attaching 2 points to generate an ordered set $2V$ within the new set.

\begin{figure}

\centerline{
\unitlength 1mm 
\linethickness{0.4pt}
\ifx\plotpoint\undefined\newsavebox{\plotpoint}\fi 
\begin{picture}(171,66)(0,0)
\put(115,25){\circle*{2}}
\put(95,45){\circle*{2}}
\put(5,35){\circle*{2}}
\put(35,55){\circle*{2}}
\put(115,35){\circle*{2}}
\put(95,55){\circle*{2}}
\put(125,25){\circle*{2}}
\put(105,45){\circle*{2}}
\put(15,35){\circle*{2}}
\put(25,55){\circle*{2}}
\put(125,35){\circle*{2}}
\put(105,55){\circle*{2}}
\put(25,25){\circle*{2}}
\put(15,25){\circle*{2}}
\put(135,25){\circle*{2}}
\put(115,45){\circle*{2}}
\put(25,35){\circle*{2}}
\put(15,55){\circle*{2}}
\put(135,35){\circle*{2}}
\put(115,55){\circle*{2}}
\put(60,25){\circle*{2}}
\put(60,45){\circle*{2}}
\put(60,5){\circle*{2}}
\put(150,45){\circle*{2}}
\put(80,25){\circle*{2}}
\put(100,25){\circle*{2}}
\put(80,45){\circle*{2}}
\put(80,5){\circle*{2}}
\put(170,45){\circle*{2}}
\put(35,25){\circle*{2}}
\put(5,25){\circle*{2}}
\put(145,25){\circle*{2}}
\put(125,45){\circle*{2}}
\put(60,35){\circle*{2}}
\put(60,55){\circle*{2}}
\put(60,15){\circle*{2}}
\put(150,55){\circle*{2}}
\put(80,35){\circle*{2}}
\put(80,55){\circle*{2}}
\put(80,15){\circle*{2}}
\put(170,55){\circle*{2}}
\put(35,35){\circle*{2}}
\put(5,55){\circle*{2}}
\put(145,35){\circle*{2}}
\put(125,55){\circle*{2}}
\put(70,25){\circle*{2}}
\put(90,25){\circle*{2}}
\put(70,45){\circle*{2}}
\put(70,5){\circle*{2}}
\put(160,45){\circle*{2}}
\put(50,25){\circle*{2}}
\put(50,45){\circle*{2}}
\put(50,5){\circle*{2}}
\put(140,45){\circle*{2}}
\put(70,35){\circle*{2}}
\put(70,55){\circle*{2}}
\put(70,15){\circle*{2}}
\put(160,55){\circle*{2}}
\put(50,35){\circle*{2}}
\put(50,55){\circle*{2}}
\put(50,15){\circle*{2}}
\put(140,55){\circle*{2}}
\put(60,25){\line(0,1){10}}
\put(60,45){\line(0,1){10}}
\put(60,5){\line(0,1){10}}
\put(150,45){\line(0,1){10}}
\put(25,25){\line(0,1){10}}
\put(80,25){\line(0,1){10}}
\put(80,45){\line(0,1){10}}
\put(80,5){\line(0,1){10}}
\put(170,45){\line(0,1){10}}
\put(70,25){\line(0,1){10}}
\put(70,45){\line(0,1){10}}
\put(70,5){\line(0,1){10}}
\put(160,45){\line(0,1){10}}
\put(50,25){\line(0,1){10}}
\put(50,45){\line(0,1){10}}
\put(50,5){\line(0,1){10}}
\put(140,45){\line(0,1){10}}
\put(35,25){\line(0,1){10}}
\put(15,35){\line(2,-1){20}}
\put(125,35){\line(2,-1){20}}
\put(105,55){\line(2,-1){20}}
\put(5,35){\line(2,-1){20}}
\put(115,35){\line(2,-1){20}}
\put(95,55){\line(2,-1){20}}
\put(125,25){\line(2,1){20}}
\put(105,45){\line(2,1){20}}
\put(115,25){\line(2,1){20}}
\put(95,45){\line(2,1){20}}
\put(145,25){\line(-1,1){10}}
\put(125,45){\line(-1,1){10}}
\put(135,25){\line(-1,1){10}}
\put(115,45){\line(-1,1){10}}
\put(125,25){\line(-1,1){10}}
\put(105,45){\line(-1,1){10}}
\put(135,25){\line(1,1){10}}
\put(115,45){\line(1,1){10}}
\put(70,5){\line(1,1){10}}
\put(160,45){\line(1,1){10}}
\put(125,25){\line(1,1){10}}
\put(105,45){\line(1,1){10}}
\put(60,5){\line(1,1){10}}
\put(150,45){\line(1,1){10}}
\put(115,25){\line(1,1){10}}
\put(95,45){\line(1,1){10}}
\put(50,5){\line(1,1){10}}
\put(140,45){\line(1,1){10}}
\put(145,25){\line(-3,1){30}}
\put(125,45){\line(-3,1){30}}
\put(80,5){\line(-3,1){30}}
\put(170,45){\line(-3,1){30}}
\put(145,35){\line(-3,-1){30}}
\put(125,55){\line(-3,-1){30}}
\put(2,25){\makebox(0,0)[rc]{\footnotesize $\overline{2V} $}}
\put(3,45){\makebox(0,0)[rc]{\footnotesize $\widehat{2V} $}}
\put(48,25){\makebox(0,0)[rc]{\footnotesize $\overline{4C_2 } $}}
\put(48,45){\makebox(0,0)[rc]{\footnotesize $\widehat{4C_2 } $}}
\put(113,25){\makebox(0,0)[rc]{\footnotesize $\overline{S_4 } $}}
\put(93,45){\makebox(0,0)[rc]{\footnotesize $\widehat{S_4 } $}}
\put(47,5){\makebox(0,0)[rc]{\footnotesize $\overline{C_8 } $}}
\put(138,45){\makebox(0,0)[rc]{\footnotesize $\widehat{C_8 } $}}
\put(5,55){\line(1,-2){5}}
\put(25,55){\line(1,-2){5}}
\put(10,45){\line(1,2){5}}
\put(30,45){\line(1,2){5}}
\put(10,45){\circle*{2}}
\put(30,45){\circle*{2}}
\put(5,35){\line(0,-1){10}}
\put(5,25){\line(1,1){10}}
\put(25,35){\line(-1,-1){10}}
\put(15,25){\line(2,1){20}}
\put(10,25){\oval(14,4)[]}
\put(30,25){\oval(14,4)[]}
\put(5,55){\line(1,1){10}}
\put(50,55){\line(1,1){10}}
\put(95,55){\line(1,1){10}}
\put(140,55){\line(1,1){10}}
\put(35,55){\line(-1,1){10}}
\put(80,55){\line(-1,1){10}}
\put(125,55){\line(-1,1){10}}
\put(170,55){\line(-1,1){10}}
\put(15,65){\line(1,-1){10}}
\put(60,65){\line(1,-1){10}}
\put(105,65){\line(1,-1){10}}
\put(150,65){\line(1,-1){10}}
\put(25,65){\line(-1,-1){10}}
\put(70,65){\line(-1,-1){10}}
\put(115,65){\line(-1,-1){10}}
\put(160,65){\line(-1,-1){10}}
\put(15,65){\circle*{2}}
\put(60,65){\circle*{2}}
\put(105,65){\circle*{2}}
\put(150,65){\circle*{2}}
\put(25,65){\circle*{2}}
\put(70,65){\circle*{2}}
\put(115,65){\circle*{2}}
\put(160,65){\circle*{2}}
\put(115,35){\line(4,-1){40}}
\put(50,35){\line(4,-1){40}}
\put(50,15){\line(4,-1){40}}
\put(125,35){\line(4,-1){40}}
\put(60,35){\line(4,-1){40}}
\put(60,15){\line(4,-1){40}}
\put(155,25){\line(-2,1){20}}
\put(90,25){\line(-2,1){20}}
\put(90,5){\line(-2,1){20}}
\put(165,25){\line(-2,1){20}}
\put(100,25){\line(-2,1){20}}
\put(100,5){\line(-2,1){20}}
\put(155,25){\circle*{2}}
\put(90,5){\circle*{2}}
\put(165,25){\circle*{2}}
\put(100,5){\circle*{2}}
\put(20,35){\oval(34,4)[]}
\put(65,5){\oval(34,4)[]}
\put(95,5){\oval(14,4)[]}
\put(65,15){\oval(34,4)[]}
\end{picture}
}

\caption{The forbidden configurations constructed from
$2V$, $4C_2 $, $S_4 $ and $C_8 $ via
Definition \ref{forbconfadd2}.
Dictated orbit structures that are not equal to the natural one
are indicated with ovals. As a mnemonic aid: Hats indicate a
peaked roof, bars a flat one.}
\label{forbconf_hats_bars}

\end{figure}
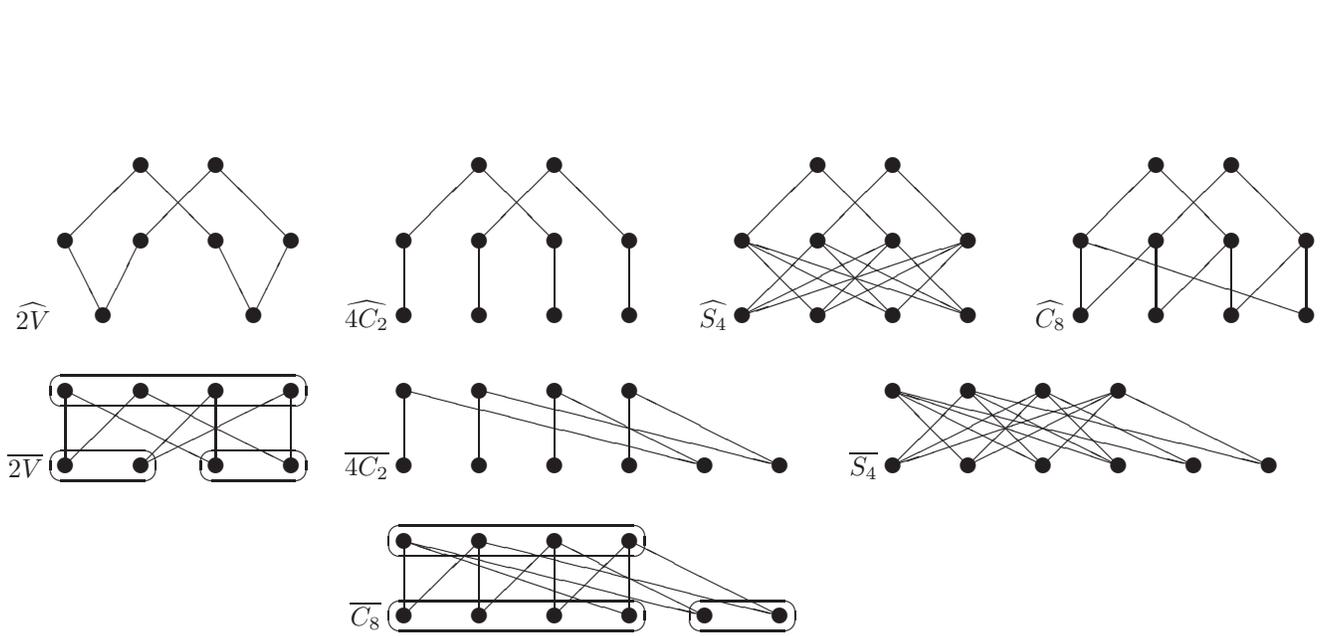

\begin{define}
\label{forbconfadd2}

(Also see Figure \ref{forbconf_hats_bars}.)
Let $Q\in \{ 2V, 4C_2 , S_4 , C_8 \} $,
let
$B$
be the set of minimal
elements of $Q$, let
$T=\{ t_1 , t_2 , t_3 , t_4 \} $
be the set of maximal elements of $Q$, labeled such that, in case
$Q\in \{ 2V, C_8 \} $, $t_1 $ and $t_3 $ do not have a common lower bound
and $t_2 $ and $t_4 $ do not have a common lower bound.
Let $A=\{ a_1 , a_2 \} $ be a set of two additional points
that are not in $Q$.
\begin{enumerate}
\item
We define the
ordered set $\widehat{Q} $ to be the ordered set
obtained by adding $a_1 $ and $a_2 $ to $Q$ such that
$a_1 $ is an upper bound of $t_1 $ and $t_3 $
and such that
$a_2 $ is an upper bound of $t_2 $ and $t_4 $.
We define ${\cal D} \left( \widehat{Q} \right) :=\{ B,T,A\} $.

\item
We define the
ordered set $\overline{Q} $ to be the ordered set
obtained by adding $a_1 $ and $a_2 $ to $Q$ such that
$a_1 $ is a lower bound of $t_1 $ and $t_3 $
and such that
$a_2 $ is a lower bound of $t_2 $ and $t_4 $.
We define ${\cal D} \left( \overline{Q} \right) :=\{ B,T,A\} $.

\end{enumerate}

\end{define}

\begin{define}
\label{forbconfdef}

An interdependent orbit union
$(U, {\cal D})$ will be called a
{\bf forbidden configuration} iff
$U$ is max-locked of height $1$ or
if
$(U, {\cal D})$ is isomorphic
or dually isomorphic to
any of
$\left( C_8 , {\cal D} \left( C_8 \right) \right)
$,
$\left( 2C_3 , {\cal D} \left( 2C_3 \right) \right)
$,
(see Figures \ref{forbconf_sporadic} and \ref{forbconf_hats_bars} for the
remaining sets)
$\left( 2C_3 ^* , {\cal D} \left( 2C_3 ^* \right) \right)
$,
$\left( 2V , {\cal D} \left( 2V \right) \right)
$,
$\left( \widehat{2V } , {\cal D} \left( \widehat{2V } \right) \right)
$,
$\left( \overline{2V } , {\cal D} \left( \overline{2V } \right) \right)
$,
$\left( \widehat{C_8 } , {\cal D} \left( \widehat{C_8 } \right) \right)
$,
$\left( \overline{C_8 } , {\cal D} \left( \overline{C_8 } \right) \right)
$,
$\left( \widehat{4C_2 } , {\cal D} \left( \widehat{4C_2 } \right) \right)
$,
$\left( \overline{4C_2 } , {\cal D} \left( \overline{4C_2 } \right) \right)
$,
$\left( \widehat{S_4 } , {\cal D} \left( \widehat{S_4 } \right) \right)
$,
$\left( \overline{S_4 } , {\cal D} \left( \overline{S_4 } \right) \right)
$,
$\left( \widetilde{S_4 } , {\cal D} \left( \widetilde{S_4 } \right) \right)
$.

\end{define}

\begin{lem}
\label{forboff1from5adeq}

Every forbidden configuration with $3$
dictated orbits
is
offset by $1$ from being
$5$-adequately bounded.

\end{lem}

{\bf Proof.}
For $i\in \{ 1,2,3\} $, we say that
$(U, {\cal D} )$ is of {\bf type $i$}
iff $(U, {\cal D} )$
has $2^i $ automorphisms and $4+2i$ points.

The ordered sets
$\left( 2C_3 , {\cal D} \left( 2C_3 \right) \right)
$,
$\left( 2C_3 ^* , {\cal D} \left( 2C_3 ^* \right) \right)
$ and
$(2V, {\cal D} (2V))$
have
$2$ automorphisms and $6$ elements;
$\left( \widetilde{S_4 } , {\cal D} \left( \widetilde{S_4 } \right) \right) $
has
$8$ automorphisms and $10$ elements;
$\left( \widehat{2V } , {\cal D} \left( \widehat{2V } \right) \right)
$ and
$\left( \overline{2V } , {\cal D} \left( \overline{2V } \right) \right)
$ each
have
$4$ automorphisms and $8$ points.

For the remaining
forbidden configurations
$(U, {\cal D} )$
with 3 dictated orbits,
note that each has $10$ points.
Let the
orbit with 2 elements be $D_3 $
and assume $D_3 \updownharpoons _{\cal D} D_2 $.
Then $D_3 \cup D_2 $ is isomorphic to
$2V$ or its dual, which, with its natural
automorphism structure, has $8$ automorphisms.
Moreover,
$D_1 \cup D_2 $
is a set $C_8 $, $S_4 $ or $4C_2 $,
and in each case\footnote{In the case of $C_8 $,
this is because of the
way the two-element orbit was attached.},
we obtain
$|{\rm Aut} _{\cal D} (U)|
=
8$.

Thus, every
forbidden configuration with $3$
dictated orbits
is of type $i$ with $i\in \{ 1,2,3\} $.
Now note that $2^i \leq (i+1)!$ and
$\left\lfloor {1\over 2} (4+2i-3+1)\right\rfloor
=i+1$.
\qed

\vspace{.1in}

Our induction's
base step with ${\cal D} $ having 3 orbits,
is Lemma \ref{addonetomaxlock} below.
Within its proof we will introduce further
useful items that will be used later.

\begin{lem}
\label{addonetomaxlock}

Let $(U, {\cal D})$
with ${\cal D} =\{ D_1 , D_2 , D_3 \} $
be a
flexible
tight
interdependent orbit union
that is not a
forbidden
configuration.
Then
$(U, {\cal D} )$ is
$\max\{ w(U), 5\} $-adequately bounded.

\end{lem}

{\bf Proof.}
This proof consists of separately discussing the cases that there are
$\{ i,j,k\} =\{ 1,2,3\} $ such that
$|D_i |
>|D_j |+|D_k |$, or,
such that
$|D_i |
=|D_j |+|D_k |$, or,
$|D_1 |=|D_2 |=|D_3 |=2$, or,
$|D_i |
<|D_j |+|D_k |$ and not all three are equal to $2$.
Cases and subcases will be indicated in {\em italics}.

{\em Case 1:
There are $\{ i,j,k\} =\{ 1,2,3\} $ such that
$|D_i |
>|D_j |+|D_k |$.}
In this case, by Lemma \ref{allbutoneac},
$|{\rm Aut} _{\cal D} (U)|
\leq
|D_j |!|D_k |!
\leq
(|D_j |+|D_k |-1)!$.
Moreover,
$\left\lfloor {1\over 2} (|D_i |+|D_j |+|D_k |-3) ) \right\rfloor
\geq
\left\lfloor {1\over 2} (2|D_j |+2|D_k |-2) ) \right\rfloor
=
|D_j |+|D_k |-1$.
Hence, in this case,
$(U,{\cal D} )$ is $w(U)$-adequately bounded.

{\em Case 2.
There are $\{ i,j,k\} =\{ 1,2,3\} $ such that
$|D_i |=|D_j |+|D_k |$.}
In case
$|{\rm Aut} _{\cal D} (U)|\leq (|D_i |-2)!$, because
$\left\lfloor {1\over 2} (|D_1 |+|D_2 |+|D_3 |-3) ) \right\rfloor
=
\left\lfloor {1\over 2} (2|D_i |-3) ) \right\rfloor
=
|D_i |-2$, we obtain that
$(U,{\cal D} )$ is $w(U)$-adequately bounded.
Thus, for the remainder of this case,
we can assume that
$|{\rm Aut} _{\cal D} (U)|> (|D_i |-2)!$.
Because
$|{\rm Aut} _{\cal D} (U)|
\leq
|D_j |!|D_k |!$,
by Lemma \ref{factprodcombine} and by symmetry,
we only need to consider the cases
$(|D_j |, |D_k |)\in \{ (3,3), (3,4)\} \cup \{ (2,n):n\geq 2\} $.
By Lemma \ref{nontrivgcd}, the case $(|D_j |, |D_k |)= (3,4)$, which
leads to $|D_i |=7$, does not occur, and, in the case $(|D_j |, |D_k |)= (2,n)$,
$n$ must be even.

Because, for two
directly interdependent tight orbits,
any two elements of
the first orbit are comparable to the same number of elements of
the second orbit,
the following definition is safe.

\begin{define}

Let $(U, {\cal D} )$ be a tight interdependent
orbit union and let
\linebreak
$D_a \updownharpoons _{\cal D} D_b $.
We define $z_{a\to b} $ to be the number of elements of
$D_b $ that are comparable to an individual element of $D_a $.

\end{define}

We will frequently use that $z_{a\to b} |D_a |=z_{b\to a} |D_b |$.

{\em Case 2.1. $|D_j |=|D_k |=2$, $|D_i |=4$.}
First, we note that, if one of $D_j $ or $D_k $ were not
directly interdependent with $D_k $, then
$D_k $ would contain two order-autonomous $2$-antichains,
which cannot be. Thus $D_j , D_k \updownharpoons _{\cal D} D_i $.
Therefore,
we obtain that
$z_{j\to i}=z_{k\to i} =2$ and
$z_{i\to j}=z_{i\to k} =1$ and then, because $D_i $ cannot contain any
nontrivial order-autonomous antichains,
$(U,{\cal D} )$ is isomorphic or dually isomorphic to
one of
$\left( \widehat{2V } , {\cal D} \left( \widehat{2V } \right) \right)
$
and
$\left( \overline{2V } , {\cal D} \left( \overline{2V } \right) \right)
$.

{\em Case 2.2. $|D_j |=2$, $|D_k |=n$, $n\geq 4$, $|D_i |=2+n$.}
Because the case
$|{\rm Aut} _{\cal D} (U)|\leq (|D_i |-2)!=n!$
has already been discussed,
and because
$|{\rm Aut} _{\cal D} (U)|
\leq
\alpha _{\cal D} (D_j )
\alpha _{\cal D} (D_k )
\leq 2!n!$,
the
only case
that needs to be considered is
$\alpha _{\cal D} (D_j ) =2$ and
$\alpha _{\cal D} (D_k ) =n!$.
In this case,
$D_k \not\updownharpoons _{\cal D} D_j $, and hence
$D_k \updownharpoons _{\cal D} D_i $.
Now
$z_{i\to k} \in \{ 2, \ldots , |D_k | -2\} $ and, because
$\alpha _{\cal D} (D_k ) =n!$, we conclude that
$|D_i |\geq \pmatrix{ n\cr z_{i\to k} \cr } \geq
\pmatrix{ n\cr 2\cr } ={1\over 2} n(n-1)$.
For $n>4$, we thus obtain that
$\pmatrix{ n\cr z_{i\to k} \cr } >n+2$, which cannot be.
For $n=4$, we have
$z_{i\to k} =2$, and,
for every $\Phi \in {\rm Aut} _{\cal D} (U)$,
$\Phi |_{D_k } $ determines first $\Phi |_{D_i} $ and then $\Phi $ itself.
Hence, in this case
$|{\rm Aut} _{\cal D} (U)|\leq
|D_k |!=n!$, which was already discussed.

{\em Case 2.3. $|D_j |=|D_k |=3$, $|D_i |=6$.}
Assume without loss of generality that
$D_i \upharpoonleft _{\cal D} D_j $.
Independent
of whether
$z_{i\to j} =1$
or
$z_{i\to j} =2$,
we have that
$D_i $ consists of three
$D_i \cup D_j $-order-autonomous antichains
$D_i ^1, D_i ^2 , D_i ^3 $
with two elements each.
Because $D_i $ does not contain any
nontrivial
$D_i \cup D_j \cup D_k $-order-autonomous subsets,
we conclude that
$D_i \updownharpoons _{\cal D} D_k $.
Similar
to the preceding argument,
$D_i $ consists of three
$D_i \cup D_k $-order-autonomous antichains with two elements each, such that
each intersects exactly two $D_i ^\ell $.
Therefore, for every $\Phi \in {\rm Aut} _{\cal D} (U)$, the
$3!$ ways that $\Phi $ could permute the $D_i ^\ell $ as sets
determine
$\Phi |_{D_k } $, which then determines how
$\Phi $ maps every $D_i ^\ell $ to its image.
Hence
$|{\rm Aut} _{\cal D} (U)|
\leq
3!$ and
$(U,{\cal D} )$ is $w(U)$-adequately bounded.

{\em Case 3:
$|D_1 |=|D_2 |=|D_3 |=2$.}
In this case, $(U,{\cal D})$ is one of
$\left( 2C_3 , {\cal D} \left( 2C_3 \right) \right)
$,
$\left( 2C_3 ^* , {\cal D} \left( 2C_3 ^* \right) \right)
$, and
$\left( 2V , {\cal D} \left( 2V \right) \right)
$.

{\em
Case 4: For all
$\{ i,j,k\} =\{ 1,2,3\} $, we have
$|D_i |< |D_j |+|D_k |$, and there is a
$z$ with $|D_z |>2$.}
First note that, for all $i\in \{ 1,2,3\} $, we have
$\left\lfloor {1\over 2} (|D_1 |+|D_2 |+|D_3 |-3) ) \right\rfloor
\geq
\left\lfloor {1\over 2} (2|D_i |-2) ) \right\rfloor
=
|D_i |-1$.
Hence, if
$|{\rm Aut} _{\cal D} (U)|\leq (|D_i |-1)!$,
then $(U,{\cal D} )$ is $w(U)$-adequately bounded.
We can assume
$|{\rm Aut} _{\cal D} (U)|> (|D_i |-1)!$
for the remainder of this case.

Via Lemma \ref{nontrivgcd}, it is not possible for
two orbits to have two elements and the third to have $3<2+2$ elements.
Therefore, at most one
orbit has
$2$ elements.

{\em Case 4.1:
For all $\{ i,j,k\} =\{ 1,2,3\} $ such that
$D_i \updownharpoons _{\cal D} D_j , D_k $, we have that
$D_i $ does not contain any nontrivial
$D_i \cup D_k $-
or
$D_i \cup D_k $-order-autonomous antichains.}
In this case, for any $\Phi \in {\rm Aut} _{\cal D} (U)$
and $D_i \updownharpoons _{\cal D} D_j $, the
restriction $\Phi | _{D_i } $ determines
$\Phi |_{D_i \cup D_j } $ which determines $\Phi $.
Hence
${\rm Aut} _{\cal D} (U)|\leq
\min \{ \alpha _{\cal D} (D_1 ), \alpha _{\cal D} (D_2 ), \alpha _{\cal D} (D_3 ) \} $.
The case that there is an $i$ such that
$\alpha _{\cal D} (D_i )\leq (|D_i |-1)!$
has already been discussed.
This leaves the case that, for all $i\in \{ 1,2,3\} $, we have that
$\alpha _{\cal D} (D_i )> (|D_i |-1)!$.
In case $\max \{ |D_1 |, |D_2 |, |D_3 |\} \not= 4$, this implies that
$(U,{\cal D} )$ is max-locked, which was excluded.
This leaves the case that
$\max \{ |D_1 |, |D_2 |, |D_3 |\} = 4$.
In this case,
by Lemma \ref{nontrivgcd},
we first conclude that two interdependent orbits have size $2$ or $4$,
and then
that the third size
$|D_i |$
cannot be equal to
$3$ either.
Because any $D_i \cup D_j $ with
$D_i \updownharpoons _{\cal D} D_j $,
$|D_i |=2$ and $|D_j |=4$
contains
nontrivial
$D_i \cup D_j $-order-autonomous antichains,
the last case to consider is
$|D_1 |= |D_2 |= |D_3 | = 4$.
In this case,
${\rm Aut} _{\cal D} (U)|\leq
4!$ and
$
\left\lfloor {1\over 2} (|D_1 |+|D_2 |+|D_3 |-3) ) \right\rfloor
=
\left\lfloor {1\over 2} (12-3) ) \right\rfloor
=4$, and
$(U,{\cal D} )$ is $5$-adequately bounded.

\begin{define}

Let $(P, {\cal D})$ be
a structured
ordered set
and let $D_j , D_k
\in {\cal D}$
be so that
$D_j \updownharpoons _{\cal D} D_k $.
For every $x\in D_j $, define
$B_x ^k :=\{ z\in D_j :\updownarrow z\setminus (D_k \cup \{ z\} )
=\updownarrow x\setminus (D_k  \cup \{ x\} ) \} $.
We define the
{\bf partition of $D_j $ induced by removal of $D_k $}
to be
${\cal B} ^{j\setminus k} :=\left\{ B_x ^k :x\in D_j \right\} $.
This partition is called {\bf nontrivial}
iff it contains at least one set that is not the
underlying set and not a singleton.

\end{define}

Note that, clearly, every $\Phi \in {\rm Aut} _{\cal D} (P)$ respects
every partition
${\cal B} ^{j\setminus k} $.

{\em Case 4.2:
There are $\{ i,j,k\} =\{ 1,2,3\} $ such that
$D_i \updownharpoons _{\cal D} D_j , D_k $,
such that
$D_i $ contains nontrivial
$D_i \cup D_k $-order-autonomous antichains
and
there is a $q\in \{ j,k\} $ such that
$D_q $ does not contain nontrivial
$D_i \cup D_q $-order-autonomous antichains.}
In this case, because $D_i \updownharpoons _{\cal D} D_q $,
for any $\Phi \in {\rm Aut} _{\cal D} (U)$,
the
restriction $\Phi | _{D_i } $ determines
$\Phi |_{D_i \cup D_q } $ which determines $\Phi $,
and hence we have that
$|{\rm Aut} _{\cal D} (U)|\leq \alpha _{\cal D} (D_i )$.
By Lemma \ref{noverklem}, because every
$\Phi \in {\rm Aut} _{\cal D} (U)$
respects ${\cal B} ^{i\setminus j} $,
for $|D_i |\geq 6$,
we obtain
$|{\rm Aut} _{\cal D} (U)|\leq \alpha _{\cal D} (D_i )\leq (|D_i |-1)!$,
which has already been discussed.
Hence $|D_i |<6$.
Because $D_i $ can be partitioned
nontrivially into sets of pairwise equal sizes, we obtain that
$|D_i |$ cannot be prime.
Hence $|D_i |=4$.

Because the case
$|{\rm Aut} _{\cal D} (U)|\leq (|D_i |-1)!$
has already been treated, we can assume that
$3!<|{\rm Aut} _{\cal D} (U)|\leq \alpha _{\cal D} (D_i )\leq 4!$.
In case $|U|\geq 11$, we obtain
$\left\lfloor {1\over 2} (|D_1 |+|D_2 |+|D_3 |-3) ) \right\rfloor
\geq
\left\lfloor {1\over 2} \cdot 8 \right\rfloor
=4$ and
$(U,{\cal D} )$ is $5$-adequately bounded.
This leaves the case $|U|\leq 10$.
In this case, by Lemma \ref{nontrivgcd},
both $|D_j |$ and $|D_k |$ are even and because at most one of
them could be equal to $2$, we can assume without loss of generality that
$|D_j |=4$ and $|D_k |=2$.
Now $D_i \cup D_k $ must be isomorphic to
$2V$ or its dual.
In case $z_{i\to j} =z_{j\to i} =1$,
$D_i \cup D_j $ must be isomorphic to
$4C_2 $.
Because $D_j $ does not contain any
nontrivial $D_i \cup D_j $-order-autonomous antichains,
in case $z_{i\to j} =z_{j\to i} =2$,
$D_i \cup D_j $ must be isomorphic to
$C_8 $.
In case $z_{i\to j} =z_{j\to i} =3$,
$D_i \cup D_j $ must be isomorphic to
$S_4 $.
Therefore,
$D_i \cup D_j $ must be isomorphic to
$C_8 $, $S_4 $ or $4C_2 $.
Now $(U,{\cal D})$ must be
isomorphic to
$\left( \widehat{D_i \cup D_j } , {\cal D} \left(
\widehat{D_i \cup D_j } \right) \right) $,
$\left( \overline{D_i \cup D_j } , {\cal D} \left(
\overline{D_i \cup D_j } \right) \right) $, one of their
duals, or, in case
$D_i \cup D_j $ is isomorphic to
$S_4 $, to
$\left( \widetilde{S_4 } , {\cal D} \left( \widetilde{S_4 } \right) \right)
$, all of which are forbidden configurations.

{\em Case 4.3.
There are $\{ i,j,k\} =\{ 1,2,3\} $ such that
$D_i \updownharpoons _{\cal D} D_j , D_k $,
we have that $D_i $
and $D_k $
contain nontrivial
$D_i \cup D_k $-order-autonomous antichains
and $D_j $ contains nontrivial
$D_i \cup D_j $-order-autonomous antichains.}
In this case, we must have that
$D_j \updownharpoons _{\cal D} D_k $,
and then, by symmetry and Case 4.2,
we can assume that, for all
distinct $i,j\in \{ 1,2,3\} $, we have that
both $D_i $ and $D_j $ contain nontrivial
$D_i \cup D_j $-order-autonomous antichains.

The notation in
Section \ref{orbgraphsec},
with
the removed orbits
being the ones with the highest indices,
is convenient to state
the results so far.
However, in symmetric situations, as we encounter
here and, later,
in Lemma \ref{QiswUadequateUmforbidden}, we may need to
remove a noncutvertex $D_k $ of ${\cal O} (U,{\cal D} )$
that does not have the highest index.
To handle such situations, we introduce the following notation.

\begin{define}
\label{relabeling4.4}

With notation as in
Section \ref{orbgraphsec},
let $k\leq m$ and let
$D_k $ be a noncutvertex of ${\cal O} (U,{\cal D} )$.
We define $(U^k , {\cal D} ^k )$ and
$(Q^k , {\cal D} _{Q^k } )$ to be the ordered sets
as constructed in Section \ref{orbgraphsec}
with $D_k $ in the place of $D_m $ and vice versa.
All corresponding quantities will carry a superscript $k$.

\end{define}

For $\{ i,j,k\} =\{ 1,2,3\} $,
because $\ell _i ^k +\ell _j ^k \geq 4$,
by Lemma \ref{boundsforQ}, we have that
$(Q^k , {\cal D} _{Q^k } )$ is
induction
$w(U)$-adequately bounded.
Thus, by Lemma \ref{halfdictates2clique}
and by part \ref{Qwadeq} of Lemma \ref{withandwithoutQ},
the only case left to consider is if
every $(U^k , {\cal D} ^k )$ is an $8$-crown
or max-locked, that is, a set $pC_2 $ or
a set $S_p $.
This means that, for any $\{ i,j,k\} =\{ 1,2,3\} $, we
have $|{\cal B} ^{i\setminus k} |=|{\cal B} ^{j\setminus k} |$
because each is the number of minimal or maximal elements of $U^k $.

Without loss of generality, we can
assume that $|D_3 |\leq |D_2 |,|D_1 |$.

Suppose,
for a contradiction, that $|D_3 |=4$.
Then
$|{\cal B} ^{2\setminus 1} |=
|{\cal B} ^{3\setminus 1} |=
2$
and
$|{\cal B} ^{1\setminus 2} |=
|{\cal B} ^{3\setminus  2} |=
2$.
Because any two of our three sets are
directly interdependent, one of
$D_1 $ and $D_2 $ is the set of
maximal elements or the set of minimal elements.
By duality and because the
sizes of
$D_1 $ and $D_2 $
are immaterial for the argument,
without loss of generality, we can assume that
$D_1 $ is the set of maximal elements of $U$.
Because $(U, {\cal D} )$ is tight, and because
${\cal B} ^{2\setminus 1} $ and ${\cal B} ^{3\setminus 1} $
are nontrivial, no element
$x\in D_1 $ is above all minimal elements or above all but
one minimal element.
Because ${\cal B} ^{3\setminus 2} $ is not trivial, every
$x\in D_1 $ is above at least two elements of $D_3 $.
Because ${\cal B} ^{3\setminus 2} \not= {\cal B} ^{3\setminus 1} $,
every $x\in D_1 $ is above two elements in different sets of
${\cal B} ^{3\setminus 1} $.

Let $x\in D_1 $.
In case $D_2 \upharpoonleft _{\cal D} D_3 \upharpoonleft _{\cal D} D_1 $,
because
$x$ is above two elements in different sets of
${\cal B} ^{3\setminus 1} $ and
$|{\cal B} ^{3\setminus 1} |=2=|{\cal B} ^{2\setminus 1} |$,
we obtain $x>D_2 $, a contradiction.
Thus
$D_3 \upharpoonleft _{\cal D} D_2 \upharpoonleft _{\cal D} D_1 $
and $D_3 $ is the set of minimal elements of $U$.
Because $x$ is above some element in
$D_2 $, it is also above at least one whole set
in
${\cal B} ^{3\setminus 1} $.
Because $x$ is above at least one element in another set of
${\cal B} ^{3\setminus 1} $,
we conclude that $x$ is above at least
$3$ of the $4$
minimal elements of $U$, a contradiction.

Because $D_3 $ has a nontrivial partition,
$|D_3 |$ cannot be a prime number.
Thus $|D_3 |\geq 6$.

Note that, for $j\in \{ 1,2\} $, we have
$|{\rm Aut} _{\cal D} (U)|
\leq
\alpha _{\cal D} (D_3 )
|{\cal B} ^{j\setminus 3} |!
$.
Because $|D_3 |\geq 6$, by Lemma \ref{noverklem}, because
${\cal B} ^{3\setminus 1} $ is not trivial, we have
$\alpha _{\cal D} (D_3 )
\leq (|D_3 |-1)!$.
Moreover,
$|{\cal B} ^{j\setminus 3} |\leq {|D_j |\over 2} $.
Hence, trivially,
$|{\cal B} ^{j\setminus 3} |\leq {|D_j |\over a} $
with $a\in \{ 2,3\} $,
and $\alpha _{\cal D} (D_3 )
\leq (|D_3 |-b)!$ with $b\in \{ 1,2\} $.
For
$a=3$ or $b=2$ or $|D_3 |<|D_1 |$ (and $j=2$, $k=1$), we have the following
for $\{ j,k\} =\{ 1,2\} $, $a\in \{ 2,3\} $ and $b\in \{ 1,2\} $.
\begin{eqnarray*}
|D_3|-2b
& \leq &
|D_k|-3+\left( 1-{2\over a}\right)  |D_j |
\\
{1\over 2}
|D_3|-b
& \leq &
{1\over 2}
(|D_k|-3 )
+\left( {1\over 2}-{1\over a}\right) |D_j |
\\
{|D_j |\over a} +|D_3|-b
& \leq &
{1\over 2}
(|D_1|+|D_2|+|D_3|-3 )
\\
\left\lceil {|D_j |\over a} \right\rceil +|D_3|-b
& \leq &
\left\lfloor
{1\over 2}
(|D_1|+|D_2|+|D_3|-3 )
\right\rfloor
\end{eqnarray*}

Hence, in all these cases,
$(U,{\cal D} )$ is $w(U)$-adequately bounded, and we
are left with
the case that
$|D_1 |=|D_2 |=|D_3 |\geq 6$
and
$|{\cal B} ^{j\setminus 3} |> {|D_j |\over 3} $, that is,
$|{\cal B} ^{j\setminus 3} |= {|D_j |\over 2} $.
By symmetry, we can assume that all
$|{\cal B} ^{j\setminus k} |= {|D_j |\over 2} $.

In case $|D_j |\geq 8$, we obtain
$\alpha _{\cal D} (D_3 )
\leq
\left( {|D_3 |\over 2}\right)!
2^{|D_3 |\over 2}
\leq
(|D_3 |-2)!$
and,
by the above with $b=2$,
$(U,{\cal D} )$ is $w(U)$-adequately bounded.

Finally, suppose, for a contradiction, that
$|D_1 |=|D_2 |=|D_3 |=6$
and all ${\cal B} ^{i\setminus j} $ partition $D_i $ into doubletons.
Without loss of generality, assume that
$D_3 \upharpoonleft _{\cal D} D_2 \upharpoonleft _{\cal D} D_1 $.
Because ${\cal B} ^{3\setminus 2} $ is nontrivial, no element $x\in D_1 $
is above $5$ or more elements of $D_3 $.
Let $x\in D_1 $.
Because every set in ${\cal B} ^{2\setminus 3} $ is a doubleton,
$x$ is above at least two elements of
$D_2 $. Because no set of ${\cal B} ^{2\setminus 3} $
is contained in a set of ${\cal B} ^{2\setminus 1} $,
these two elements must be in different sets of ${\cal B} ^{2\setminus 1} $.
Hence $x$ is above at least 2 different sets of ${\cal B} ^{3\setminus 1} $,
which means $x$ is above at least 4 elements of $D_3 $.
Because
${\cal B} ^{3\setminus 2} \not= {\cal B} ^{3\setminus 1} $,
$x$ must be above exactly one
more element of $D_1 $. However, then
$x$ is above $5$ elements of $D_1 $,
a contradiction.

This concludes the proof of Lemma \ref{addonetomaxlock}.
\qed

\section{Induction Step: Bounding the Number of Automorphisms in
Non-Max-Locked
Interdependent Orbit Unions}
\label{boundnra}

As in Section \ref{orbgraphsec}, $(U, {\cal D} )$ is a
flexible tight interdependent orbit
union, but with $|{\cal D} |\geq 4$.
The idea for the proof of Theorem \ref{allowedconf} below is
an induction in which we
obtain estimates of
$|{\rm Aut} _{\cal D} (U)|$ by
applying
Theorem \ref{pruneorbit6},
Lemma \ref{withandwithoutQ}
and
Lemma \ref{inconvenientSDn2} to
a noncutvertex $D_m $ of the orbit graph ${\cal O} (U,{\cal D} )$
and the resulting
ordered sets $(U_m , {\cal D} _m )$ and
$(Q, {\cal D} _Q )$.

By part \ref{Qwadeq} of Lemma \ref{withandwithoutQ}, when
$(Q, {\cal D} _Q )$ is
induction $w(U)$-adequately bounded
and
$(U_n , {\cal D}_n )$ is $\max \{ w(U),5\} $-adequately
bounded, then
$(U , {\cal D} )$ is $\max \{ w(U),5\} $-adequately
bounded.
Hence, our last obstacle for the induction proof of
Theorem \ref{allowedconf} below, are
sets
$(Q, {\cal D} _Q )$
that are not
induction $w(U)$-adequately bounded
and forbidden configurations
$(U_n , {\cal D}_n )$.
Lemma \ref{inconvenientSDn2} and
Lemma \ref{QnotwUadequate} below show that
there are very few
sets
$(Q, {\cal D} _Q )$
that are not
induction $w(U)$-adequately bounded.

\begin{lem}
\label{QnotwUadequate}

With notation as in Section \ref{orbgraphsec},
if $D_m $ is not a cutvertex of ${\cal O} (U, {\cal D})$
and $(Q, {\cal D} _Q )$ is not
induction $w(U)$-adequately bounded
with
$\sum _{j=t} ^{m-1}
\ell _j = 2 $ and
$|{\cal S}(D_m )|=1$,
then
$|D_{m-1} |=4$, and
$(Q, {\cal D} _Q )$ is isomorphic to
$(2V, {\cal D} (2V))$, or
$\left( \overline{2V } , {\cal D} \left( \overline{2V } \right) \right)
$,
or one of their duals.
Moreover, $D_m $ is a pendant vertex of ${\cal O} (U, {\cal D})$.

\end{lem}

{\bf Proof.}
By Lemma \ref{addonetomaxlock} and because
$(Q, {\cal D}_Q )$ being $w(U)$-adequately bounded implies that
$(Q, {\cal D}_Q )$ is induction $w(U)$-adequately bounded,
$(Q, {\cal D}_Q )$ must be a forbidden configuration.
Because the
${\cal D} _Q $-orbits
$A_1 ^{m-1} $ and $A_2 ^{m-1} $
are not directly interdependent,
the orbit graph
${\cal O} (Q, {\cal D}_Q )$
must be a path with $3$ vertices, and,
because $A_1 ^{m-1} , A_2 ^{m-1} \subseteq D_{m-1} $,
the endvertices must have the same number of elements and their union
must be an antichain, too.
Now
$(Q, {\cal D}_Q )$ is isomorphic to
$(2V, {\cal D} (2V))$,
$\left( \overline{2V } , {\cal D} \left( \overline{2V } \right) \right)
$
or one of their duals, because these are the only
forbidden configurations whose orbit graph is a path
whose endvertices are antichains of the same size
whose union is an antichain, too.

Finally,
because ${\cal S} (D_m )=1$, there are no singleton sets $A_i ^j $, which means that
$D_{m-1} $ is the only orbit that is directly interdependent with
$D_m $, that is,
$D_m $ is a pendant vertex of ${\cal O} (U,{\cal D} )$.
\qed

\begin{lem}
\label{QiswUadequateUmforbidden}

With notation as in Section \ref{orbgraphsec},
if
$D_m $ is not a cutvertex of ${\cal O} (U, {\cal D})$,
$(Q, {\cal D} _Q )$ is
induction $w(U)$-adequately bounded,
and $(U_m , {\cal D} _m )$
is a forbidden configuration with $|{\cal D} _m |=3$,
then
$(U, {\cal D} )$ is
$\max \{ w(U),5\} $-adequately bounded.

\end{lem}

{\bf Proof.}
Clearly, $m=4$.
By Lemma \ref{forboff1from5adeq},
every forbidden configuration
with $3$ orbits is offset by
$1$ from being $5$-adequately bounded.
Therefore, by Lemma \ref{withandwithoutQ},
we obtain that $(U,{\cal D} )$ is
$\max \{ w(U),5\} $-adequately bounded
in case
$
|{\rm Aut } _{{\cal D}} (U)|
\leq
|{\rm Aut } _{{\cal D}_4} (U_4 )|
$,
in case
$(Q, {\cal D} _Q )$ is induction
$w(U)$-adequately $({\cal S} (D_4 )-2)$-bounded,
and in case
$(Q, {\cal D} _Q )$ is induction
$w(U)$-adequately $({\cal S} (D_4 )-1)$-bounded and
$
\left( |D_4 |
-1
\right)
+\left(
\sum _{j=s} ^{3}
|D_j|
-
\ell _Q \right)
$
is odd.

This leaves us the case that
$
|{\rm Aut } _{{\cal D}} (U)|
>
|{\rm Aut } _{{\cal D}_4} (U_4 )|
$,
$(Q, {\cal D} _Q )$ is not induction
$w(U)$-adequately $({\cal S} (D_4 )-2)$-bounded,
which means that $(Q, {\cal D} _Q )$ is induction
$w(U)$-adequately $({\cal S} (D_4 )-1)$-bounded,
and
$
\left( |D_4 |
-1
\right)
+\left(
\sum _{j=s} ^{3}
|D_j|
-
\ell _Q \right)
$
is even.
Because every orbit
in a forbidden configuration with $3$ orbits
has an even number of at most $4$ elements,
all $\ell _j $ are even.
Moreover, because $|U_4 |\leq 10$, we have
$\ell _Q=\sum _{j=s} ^{n-1}
\ell _j \in \{ 2,4,6,8,10\} $.
Finally, because
$
\left( |D_4 |
-1
\right)
+\left(
\sum _{j=s} ^{3}
|D_j|
-
\ell _Q \right)
$
is even,
we obtain that $|D_4 |$ is odd.




Without loss of generality, we can assume that
$D_1 $ is not a cutvertex of ${\cal O} (U, {\cal D} )$.
Let
$(U^1 , {\cal D} ^1 )$
be the ordered set
we obtain through removal of $D_1 $
as in Theorem \ref{pruneorbit6}.
Then $|{\cal D} ^1 |=3$.
Because $|D_4 |$ is odd and can therefore only split into
an odd number of orbits,
${\cal D} ^1 $ contains an orbit with an
odd number of elements, which means that
$(U^1 , {\cal D} ^1 )$
is not a forbidden configuration.
Hence, by Lemma \ref{addonetomaxlock},
$(U^1 , {\cal D} ^1 )$ is $\max \{ w(U), 5\} $-adequately bounded.
If the
ordered set
$(Q^1 ,{\cal D}_{Q^1 } )$,
obtained through removal of $D_1 $
as in Theorem \ref{pruneorbit6}, is
induction $w(U)$-adequately bounded, then, by part \ref{Qwadeq}
of Lemma \ref{withandwithoutQ},
$(U, {\cal D} )$ is
$\max \{ w(U),5\} $-adequately bounded.

We are left to consider the case that
$(Q^1 ,{\cal D}_{Q^1 } )$ is not induction $w(U)$-adequately bounded.
By Lemmas \ref{inconvenientSDn2} and \ref{QnotwUadequate},
we have that $|D_1 |\in \{ 2,4,6\} $.
Because $|D_4 |$ is odd,
by Lemma \ref{nontrivgcd}, $D_1 $ is not directly interdependent with $D_4 $.
Let $k\in \{ 2,3\} $ be so that $D_1 \updownharpoons _{\cal D} D_k $
and $D_k $ contains nontrivial
order-autonomous antichains induced by the removal of $D_1 $.
In case $|D_1 |=6$, by Lemma \ref{inconvenientSDn2},
we would obtain
$|D_k |=6$ and that
no subset of $D_1 \cup D_k $ would be
an order-autonomous antichain in
$U\setminus D_4 $. This would mean that
$(U_4 , {\cal D} _4 )$ contains orbits of size
$6$, contradicting the fact that
$(U_4 , {\cal D} _4 )$ is a forbidden configuration.
Hence, $|D_1 |\in \{ 2,4\} $.

By Lemmas \ref{inconvenientSDn2} and \ref{QnotwUadequate},
the orbit $D_k $
has $4$ elements.
Because $|D_4 |$ is odd,
by Lemma \ref{nontrivgcd},
$D_4 $ is not directly interdependent with
$D_k $. Because
(standard assumption for Section \ref{orbgraphsec})
$D_4 \updownharpoons _{\cal D} D_3 $, we obtain
$k=2$, $|D_2 |=4$, and
$D_4 $ is not directly interdependent with
either of $D_1 $ and $D_2 $.
Because
removal of $D_1 $ also produces an
interdependent orbit union, we conclude that
$D_2  \updownharpoons _{\cal D} D_3 $.
In particular,
by Lemma \ref{nontrivgcd},
$|D_3 |$ is even.

The only way removal of
$D_1 $ can induce nontrivial $U\setminus D_1 $-order-autonomous
antichains
in $D_2 $ is when
$D_2 $ contains
$D_2 \cup D_3 $-order-autonomous
antichains.
Because the existence of
$D_2 \cup D_3 $-order-autonomous
antichains in $D_2 $ only depends on $D_3 \cap U_4 $,
and because
$U_4 $ is a forbidden configuration, by Lemma \ref{addonetomaxlock}
and inspection of Figures \ref{forbconf_sporadic}
and \ref{forbconf_hats_bars}, we see that this
is only possible when $|D_3 \cap U_4 |=2$.
Hence $\sum _{j=t} ^3 \ell _j =2$.
Note that this means that,
to establish that $(Q, {\cal D} _Q )$ is
induction $w(U)$-adequately
$({\cal S} (D_4 )-2)$-bounded,
the sum of the bounding factorials would need to be
bounded by $\left\lfloor {1\over 2} (|D_3 |+|D_4 |-4)\right\rfloor $.

Again by Lemma \ref{nontrivgcd},
there are odd numbers
$b,c$ with $b>1$ and an even number $e$ such that $|D_3 |=eb$
and $|D_4 |=cb$.
Because every $\Phi \in {\rm Aut} _{{\cal D} _Q } (Q)$ is
determined by its restriction to $D_3 $ and respects the
partition into $A_1 ^3 $ and $A_2 ^3 $, $c\geq e+1$ would imply
$|{\rm Aut} _{{\cal D} _Q} (Q)|\leq
\left( {e\over 2} b\right) ! \left( {e\over 2} b\right) !<(eb-1)!$ and
$
\left\lfloor
{1\over 2} (eb+cb-4)
\right\rfloor
\geq
\left\lfloor
{1\over 2} ((e+c)b-4)
\right\rfloor
\geq
\left\lfloor
eb+{b\over 2} -2
\right\rfloor
\geq
eb-1$,
which would mean that
$(Q, {\cal D} _Q )$ is
induction $w(U)$-adequately
$({\cal S} (D_4 )-2)$-bounded,
which was excluded.
Hence $c\leq e-1$.
Therefore
$|{\rm Aut} _{{\cal D} _Q} (Q)|\leq \alpha _{{\cal D} _Q} (D_n ) \leq (c b)!$.
For
$b>3$,
$
\left\lfloor
{1\over 2} (cb+(c+1)b-4)
\right\rfloor
=
\left\lfloor
{1\over 2} (2cb+b-4)
\right\rfloor
\geq
cb$,
which would mean that
$(Q, {\cal D} _Q )$ is
induction $w(U)$-adequately
$({\cal S} (D_4 )-2)$-bounded, which was excluded.
Hence $c\leq e-1$ and $b=3$.
Now $|{\rm Aut} _{\cal D} (U)|
\leq
|{\rm Aut } _{{\cal D}_4} (U_4 )|
|{\rm Aut } _{{\cal D}_Q } (Q)|
\leq
8\cdot (3c)!$
and
$
\left\lfloor
{1\over 2} (|U|-|{\cal D} |)
\right\rfloor
=
\left\lfloor
{1\over 2} (3c+3e+|D_1 |+|D_2 |-4)
\right\rfloor
\geq
\left\lfloor
{1\over 2} (6c+3+|D_1 |+|D_2 |-4)
\right\rfloor
\geq
\left\lfloor
{1\over 2} (6c+5)
\right\rfloor
\geq 3c+2
$.
Because $3c+2<3e\leq w(U)$ and
$8\cdot (3c)!<(3c+2)!$,
we obtain that
$(U, {\cal D} )$ is $w(U)$-adequately bounded.
\qed



%

\begin{lem}
\label{twobad}

With notation
as in Section \ref{orbgraphsec},
let $(U, {\cal D} )$
be a flexible tight interdependent orbit union
with $|{\cal D} |\geq 4$,
let $D^1 $ and $D^2 $ be two distinct noncutvertices of
${\cal O} (U,{\cal D} )$ whose removal
leads to
sets $(Q^i , {\cal D} _{Q^i } )$ that are
not induction $w(U)$-adequately bounded, and
let
$(U_{m-2} , {\cal D} _{m-2} )$
be the
interdependent orbit union
obtained
by consecutively removing $D^1 $ and $D^2 $ by iterating the
process in Section \ref{orbgraphsec}.
Then, independent of whether
$(U_{m-2} , {\cal D} _{m-2} )$
is a
forbidden configuration or $\max \{ w(U),5\} $-adequately bounded,
$(U, {\cal D} )$ is
$\max \{ w(U),5\} $-adequately bounded.

\end{lem}

{\bf Proof.}
First suppose, for a contradiction, that
there is a $D\in {\cal D} $
such that
$
D^1 \updownharpoons _{\cal D}
D \updownharpoons _{\cal D}
D^2 $ and $D$ is nontrivially partitioned in
each ${\cal D} _{Q^i } $.
Because none of the sets from Lemma \ref{inconvenientSDn2}
contain any nontrivial order-autonomous antichains,
both sets $D^j $ must be as in Lemma \ref{QnotwUadequate}.
Thus, by Lemma \ref{QnotwUadequate},
$D^1 $ and $D^2 $ are pendant vertices
and $|D|=4$.

Let the partitions of $D$ induced by removal of $D^1 $ or $D^2 $, respectively, be
$D=A_1 ^1 \cup A_2 ^1 =A_1 ^2 \cup A_2 ^2 $.
If $A_1 ^1 \in \{ A_1 ^2 , A_2 ^2\} $, then
$A_1 ^1 $ is order-autonomous in $U$, which cannot be.
Thus $A_1 ^1 $ intersects each of
$A_1 ^2 $ and $A_2 ^2 $.
Because all $A_i ^j $ are order-autonomous
in $U\setminus (D^1 \cup D^2 )$ we obtain that
$D$ is order-autonomous in
$U\setminus (D^1 \cup D^2 )$.
Because both $D^1 $ and $D^2 $ are pendant vertices
and we have just shown that $D$ is not adjacent to any other orbits,
we conclude that $\{ D, D^1 , D^2 \} $ is a connected component of
${\cal O} (U,{\cal D} ) $. Hence it is the whole vertex set
of the orbit graph,
contradicting that
$|{\cal D} |\geq 4$.
Thus, if
$
D^1 \updownharpoons _{\cal D}
E^1 $ and $D^2 \updownharpoons _{\cal D}
E^2 $
and each $E^i $ is nontrivially partitioned in
each ${\cal D} _{Q^i } $, then
$E_1 \not= E_2 $.

Now let
$\alpha _{n-2}
:=\left| \left\{ \Phi |_{U_{m-2}} :\Phi \in {\rm Aut} (U)\right\} \right|
\leq |{\rm Aut} _{D_{m-2} } (U_{m-2} )|$
and note that
$|{\rm Aut} _{\cal D} (U)|\leq \alpha _{n-2} |{\rm Aut} (Q^1 , {\cal D} _{Q^1 } )|
|{\rm Aut} (Q^2 , {\cal D} _{Q^2 } )|$.
Also note that, by
Lemmas \ref{inconvenientSDn2} and \ref{QnotwUadequate},
we have that
$|{\rm Aut} (Q^i , {\cal D} _{Q^1 } )|\in \{ 2,4,36\} $.

In case
$(U_{m-2} , {\cal D} _{m-2} )$
is $\max \{ w(U),5\} $-adequately bounded,
an estimate similar to part \ref{Qwadeq} of Lemma \ref{withandwithoutQ}
proves that, if
$|{\rm Aut} (Q^1 , {\cal D} _{Q^1 } )|
|{\rm Aut} (Q^2 , {\cal D} _{Q^2 } )|$
is bounded by a product of factorials whose arguments are
individually at most $w(U)-1$ and
such that the sum of these arguments is bounded by
$\left\lfloor {1\over 2} (|Q_1 |+|Q_2 |-6)\right\rfloor $,
then
$(U , {\cal D} )$
is $\max \{ w(U),5\} $-adequately bounded:
We use the numbers $w_i $ that we have for 
$|{\rm Aut} _{D_{m-2} } (U_{m-2} )|$
and for 
$|{\rm Aut} (Q^1 , {\cal D} _{Q^1 } )|
|{\rm Aut} (Q^2 , {\cal D} _{Q^2 } )|$
and then note that their sum is bounded by 
$
\left\lfloor {1\over 2} (|U_{m-2} |-|{\cal D} _{m-2} |)\right\rfloor
+
\left\lfloor {1\over 2} (|Q_1 |+|Q_2 |-6)\right\rfloor 
=
\left\lfloor {1\over 2} ((|U |-(|Q_1 |+|Q_2 |)+4)-(|{\cal D} | -2 ))\right\rfloor
+
\left\lfloor {1\over 2} (|Q_1 |+|Q_2 |-6)\right\rfloor
\leq 
\left\lfloor {1\over 2} (|U |-|{\cal D} |) \right\rfloor
$.
The details
for bounding $|{\rm Aut} (Q^1 , {\cal D} _{Q^1 } )|
|{\rm Aut} (Q^2 , {\cal D} _{Q^2 } )|$
with a product of factorials 
whose arguments are bounded by $\max \{ w(U), 5\} -1$
and such that 
the sum of the arguments is bounded by 
$\left\lfloor {1\over 2} (|Q_1 |+|Q_2 |-6)\right\rfloor $
are
listed in Table \ref{lastcases1}.

\begin{table}

\centerline{
\begin{tabular}{|c|c|c|c|c|c|}
\hline
$|{\rm Aut} _{{\cal D} _{Q^1 }} (Q^1 )|
\rule{0in}{.2in}$ &
$|Q_1 |$ &
$|{\rm Aut} _{{\cal D} _{Q^2 }} (Q^2 )|$ &
$|Q_2 |$ &
$\left\lfloor {1\over 2} (|Q_1 |+|Q_2 |-6)\right\rfloor $ &
bound
\\
\hline
36 & 12 & 36 & 12 & 9 & $36\cdot 36< 4!5!$
\\
\hline
36 & 12 & 4 & 8 & 7 & $36\cdot 4= 4!3!$
\\
\hline
36 & 12 & 2 & 6 & 6 & $36\cdot 2< 5!$
\\
\hline
4 & 8 & 4 & 8 & 5 & $4\cdot 4< 4!$
\\
\hline
4 & 8 & 2 & 6 & 4 & $4\cdot 2< 4!$
\\
\hline
2 & 6 & 2 & 6 & 3 & $2\cdot 2< 3!$
\\
\hline
\end{tabular}
}

\caption{Bounding $|{\rm Aut} _{{\cal D} _{Q^1 } }(Q^1 )|
\cdot |{\rm Aut} _{{\cal D} _{Q^2 } } (Q^2 )|$
when $(U_{n-2}, {\cal D} _{n-2} )$
is $\max \{ w(U),5\} $-adequately bounded.
The offset we must subtract from $|Q_1 |+|Q_2 |$ is $s=6$ and the sum of the
arguments of the factorials in the bound must be bounded by
$\left\lfloor {1\over 2} |Q_1 |+|Q_2 |-6\right\rfloor $.}
\label{lastcases1}

\end{table}

In case
$(U_{m-2} , {\cal D} _{m-2} )$
is a forbidden configuration, we first note that,
because $E^1 \not= E^2 $, we have that
$(U_{m-2} , {\cal D} _{m-2} )$ contains $2$
distinct orbits, namely the $E^i \cap U_{m-2} $,
that have exactly $2$ elements.
Thus, by Lemma \ref{addonetomaxlock}
and inspection of Figures \ref{forbconf_sporadic}
and \ref{forbconf_hats_bars},
$(U_{m-2} , {\cal D} _{m-2} )$ is one of
$2C_2 $, $2C_3 $, $2C_3 ^* $, $2V$; $\widehat{ 2V} $, $\overline{2V} $
or their dual.

In case
$(U_{m-2} , {\cal D} _{m-2} )$ is one of
$2C_2 $, $2C_3 $, $2C_3 ^* $, $2V$
or their dual,
we have that $|U|= |Q_1 |+|Q_2 |$ and ${\cal D} =4$, or
$|U|= |Q_1 |+|Q_2 |+2$ and ${\cal D} =5$.
Thus
$
\left\lfloor {1\over 2} (|U|-|{\cal D} |)\right\rfloor
\geq
\left\lfloor {1\over 2} (|Q_1 |+|Q_2 |-4 )\right\rfloor
$.
Moreover,
$
|{\rm Aut} _{{\cal D} }(U )|
\leq
2|{\rm Aut} _{{\cal D} _{Q^1 } }(Q^1 )|
|{\rm Aut} _{{\cal D} _{Q^2 } } (Q^2 )|$
and
the details for bounding $|{\rm Aut} _{\cal D} (U)|$
with the right product of
the right factorials are
listed in Table \ref{lastcases2}.

\begin{table}

\centerline{
\begin{tabular}{|c|c|c|c|c|c|c|}
\hline
$|{\rm Aut} _{{\cal D} _{Q^1 } }(Q^1 )|
\rule{0in}{.2in}$ &
$|Q_1 |$ &
$|{\rm Aut} _{{\cal D} _{Q^2 } } (Q^2 )|$ &
$|Q_2 |$ &
$\left\lfloor {1\over 2} (|Q_1 |+|Q_2 |-4)\right\rfloor $ &
bound
\\
\hline
36 & 12 & 36 & 12 &
10 &
$2\cdot 36\cdot 36< 4!5!$
\\
\hline
36 & 12 & 4 & 8 &
8 &
$2\cdot 36\cdot 4< 4!4!$
\\
\hline
36 & 12 & 2 & 6 &
7 &
$2\cdot 36\cdot 2< 5!2!$
\\
\hline
4 & 8 & 4 & 8 &
6 &
$2\cdot 4\cdot 4< 4!2!$
\\
\hline
4 & 8 & 2 & 6 &
5 &
$2\cdot 4\cdot 2< 4!$
\\
\hline
2 & 6 & 2 & 6 &
4 &
$2\cdot 2\cdot 2< 4!$
\\
\hline
\end{tabular}
}

\caption{Bounding $|{\rm Aut} _{\cal D} (U)|$
when $(U_{n-2}, {\cal D} _{n-2} )$
is one of $2C_2 $, $2C_3 $, $2C_3 ^* $, $2V$ or their dual.
The sum of the
arguments of the factorials in the bound must be bounded by
$\left\lfloor {1\over 2} (|U |-|{\cal D} |)\right\rfloor $, for which a
lower bound is provided in the second to last column.
}
\label{lastcases2}

\end{table}

Finally,
in case
$(U_{m-2} , {\cal D} _{m-2} )$ is one of
$\widehat{ 2V} $, $\overline{2V} $
or their dual,
we have that $|U|= |Q_1 |+|Q_2 |+4$ and ${\cal D} =5$.
Thus
$
\left\lfloor {1\over 2} (|U|-|{\cal D} |)\right\rfloor
=
\left\lfloor {1\over 2} (|Q_1 |+|Q_2 |+4-5 )\right\rfloor
$.
Moreover,
$
|{\rm Aut} _{{\cal D} }(U )|
\leq
4|{\rm Aut} _{{\cal D} _{Q^1 } }(Q^1 )|
|{\rm Aut} _{{\cal D} _{Q^2 } } (Q^2 )|$
and
the details for bounding $|{\rm Aut} _{\cal D} (U)|$
with the right product of
the right factorials are
listed in Table \ref{lastcases3}.
\qed

\begin{table}

\centerline{
\begin{tabular}{|c|c|c|c|c|c|c|}
\hline
$|{\rm Aut} _{{\cal D} _{Q^1 }} (Q^1 )|
\rule{0in}{.2in}$ &
$|Q_1 |$ &
$|{\rm Aut} _{{\cal D} _{Q^2 }} (Q^2 )|$ &
$|Q_2 |$ &
$|U|$ &
$\left\lfloor {1\over 2} (|U |-|{\cal D} |)\right\rfloor $ &
bound
\\
\hline
36 & 12 & 36 & 12 &
28 & 11 &
$4\cdot 36\cdot 36< 5!5!$
\\
\hline
36 & 12 & 4 & 8 &
24 & 9 &
$4\cdot 36\cdot 4= 4!4!$
\\
\hline
36 & 12 & 2 & 6 &
22 & 8 &
$4\cdot 36\cdot 2< 4!4!$
\\
\hline
4 & 8 & 4 & 8 &
20 & 7 &
$4\cdot 4\cdot 4< 4!3!$
\\
\hline
4 & 8 & 2 & 6 &
18 & 6 &
$4\cdot 4\cdot 2< 4!2!$
\\
\hline
2 & 6 & 2 & 6 &
16 & 5 &
$4\cdot 2\cdot 2< 4!$
\\
\hline
\end{tabular}
}

\caption{Bounding $|{\rm Aut} _{\cal D} (U)|$
when $(U_{n-2}, {\cal D} _{n-2} )$
is one of $\widehat{ 2V} $, $\overline{2V} $ or their dual.
The sum of the
arguments of the factorials in the bound must be bounded by
$\left\lfloor {1\over 2} (|U |-|{\cal D} |)\right\rfloor
=
\left\lfloor {1\over 2} (|Q_1 |+|Q_2 |-1 )\right\rfloor
$,
which is provided in the second to last column.}
\label{lastcases3}

\end{table}

\vspace{.1in}

We can now prove that flexible tight interdependent orbit unions
$(U, {\cal D})$ that are not
forbidden are
$\max \{ w(U), 5\} $-adequately bounded.

\begin{theorem}
\label{allowedconf}

Let $(U, {\cal D})$
be a flexible tight interdependent orbit union.
If $(U, {\cal D})$ is not a forbidden configuration, then
$(U, {\cal D})$ is
$\max \{ w(U), 5\} $-adequately bounded.

\end{theorem}

{\bf Proof.}
The proof is an induction on
$|{\cal D}|$.
For the base case, we note that, by Lemmas \ref{halfdictates2clique}
and \ref{addonetomaxlock},
the result holds for
$|{\cal D} |\in \{ 2,3\} $.

{\em Induction Step, $|{\cal D} |\geq 4$.}
Because every graph has at least $2$ noncutvertices,
let $D^1 , D^2 \in {\cal D} $
be noncutvertices of ${\cal O} (U, {\cal D} )$.

In case one of the resulting sets $(Q ^i , {\cal D} _{Q^i } )$ is induction
$w(U)$-adequately bounded, then,
by Lemma \ref{QiswUadequateUmforbidden}
and by part \ref{Qwadeq} of Lemma \ref{withandwithoutQ},
$(U, {\cal D})$ is $\max \{ w(U), 5\} $-adequately bounded.
Otherwise,
by Lemma \ref{twobad},
$(U, {\cal D})$ is $\max \{ w(U), 5\} $-adequately bounded.
\qed

\section{The Automorphism Conjecture for Ordered Sets of Small Width}
\label{ACsmwidthproof}

To finally prove the Automorphism Conjecture for
ordered sets of width up to $11$ in Theorem \ref{ACw<=12},
we first establish a
lower bound on the number of endomorphisms in Lemma \ref{2tonforbddwidth}.
Then,
from
Proposition \ref{usingEndD}
through
Proposition \ref{manyinmaxlock}, we show
that a relative abundance of
forbidden configurations
guaranteees that the Automorphism Conjecture holds.
Lemma \ref{wid12AClem} combines the
work so far into an upper bound for the number of automorphisms
when there are no nontrivial order-autonomous antichains, and
Proposition \ref{onelexiter} shows that
a single execution of the lexicographic sum construction,
as occurs for example when nontrivial order-autonomous antichains are
inserted, does not affect the status of the Automorphism Conjecture.

\begin{lem}
\label{2tonforbddwidth}

Let $w\in {\mat N}$ and $\varepsilon >0$.
There is an $N\in {\mat N}$ such that
every
ordered set $P$ of width $\leq w$ with $n:=|P|\geq N$ elements
has at least
$2^{\left( 1
-\varepsilon \right) n} $
endomorphisms.

\end{lem}

{\bf Proof.}
Let $P$ be an ordered set of height $h$ with $n$ elements.
The proof of Theorem 1 in \cite{DRSW} (on page 20 of \cite{DRSW})
shows that $P$ has
at least $2^{{h\over h+1} n} $
endomorphisms that are surjective onto a chain
of length $h$.

Let $N\in {\mat N}$ be so that
${(N/w)-1\over [(N/w)-1]+1}>1-\varepsilon $ and let
$P$ be an ordered set of width $w$ with $n\geq N$ elements.
Then the height $h$ of $P$ satisfies
$h\geq {n\over w}-1 \geq {N\over w} -1$.
Hence
$P$ has
at least $2^{{h\over h+1} n} \geq 2^{(1-\varepsilon )n} $
endomorphisms.
\qed

\vspace{.1in}

Similar to ${\rm Aut} _{\cal D} (P)$,
we define
${\rm End} _{\cal D} (P)$.

\begin{define}
\label{EndDdef}

Let $(P, {\cal D} )$ be
a structured
ordered set.
We define ${\rm End} _{\cal D} (P)$ to be the set of
order-preserving maps $f:P\to P$ such that,
for all $D\in {\cal D}$, we have
$f[D]\subseteq D$.

\end{define}

\begin{prop}
\label{usingEndD}

Let $P$ be an ordered set and let
${\cal U}$ denote the set of
nontrivial natural
interdependent orbit unions of $P$. Then
$${
|{\rm Aut} (P)|
\over
|{\rm End} (P)|
}
=
{
|{\rm Aut}_{\cal N} (P)|
\over
|{\rm End}_{\cal N} (P)|
}
\leq
\prod _{U\in {\cal U} }
{
|{\rm Aut}_{{\cal N}|U} (U)|
\over
|{\rm End}_{{\cal N}|U} (U)|
}
.$$

\end{prop}

{\bf Proof.}
The functions in each
set ${\rm End}_{{\cal N}|U} (U)$ can be combined into
endomorphisms of $P$ in the same way as the functions in
${\rm Aut}_{{\cal N}|U} (U)$ are combined in
Proposition \ref{getallfromiou}.
\qed

\begin{lem}
\label{forbconfbound}

Let $(U,{\cal D} )$ be a forbidden configuration
that is not a max-locked interdependent
orbit union of height $1$.
Then
${
|{\rm Aut}_{\cal D} (U)|
\over
|{\rm End}_{\cal D} (U)|
}
\leq
{3\over 4} .$

\end{lem}

{\bf Proof.}
Giving the details for each of the
$13$ possibilities that must be considered after accounting for
duality would be tedious at best.
Instead, simply note
that each of these forbidden configurations has 2 or 3 orbits
and that there
are a sufficient number of ways to collapse
each orbit $D_j $ to a point $x_j \in D_j $ in an order-preserving fashion.
\qed

\begin{lem}
\label{maxlockh1presrank}

$|{\rm End}_{\cal N} (wC_2 )|\geq w^w $
and
$|{\rm End}_{\cal N} (S_w )|\geq (w-1)^w $.

\end{lem}

{\bf Proof.}
The result about
$|{\rm End}_{\cal N} (wC_2 )|$ follows from the fact that every self map of
the minimal elements can be extended to a function in
${\rm End}_{\cal N} (wC_2 )$.
For the result about
$|{\rm End}_{\cal N} (S_w )|$, let $t\in S_w $ be maximal in $S_w $.
Now we can map the maximal elements of $S_w $
to $t$, we can map
every one of
the $w$ minimal elements of $S_w $
to any of the $w-1$ minimal elements in
$\downarrow t$, and, through each such choice,
we obtain a function in
${\rm End}_{\cal N} (S_w )$.
\qed

\begin{prop}
\label{manyinmaxlock}

The Automorphism Conjecture is true for the class
of ordered sets $P$
such that at least
$\lg (|P|)$ elements of $P$ are contained in
forbidden configurations.

\end{prop}

{\bf Proof.}
Let $P$ be an ordered set
with $|P|=n$ elements
such that at least
$\lg (n)$ elements of $P$ are contained in
forbidden configurations.
Recall that all
forbidden configurations that are not max-locked of height $1$
have at most $10$ elements.

In case
$\geq \sqrt{\lg (n)} $ elements of $P$ are contained in
one max-locked interdependent orbit union
$U$ of height
$1$ and width $w$, we have
$w\geq {1\over 2}
\sqrt{\lg (n)} $.
By
Lemma \ref{maxlockh1presrank},
we obtain
$|{\rm End}_{{\cal N}|U} (U)|
\geq (w-1)^w $.
Moreover, $U$ has
exactly $w!$
automorphisms. By Proposition \ref{usingEndD}, we obtain the following.
$${ |{\rm Aut} (P)|
\over
|{\rm End} (P)|}
\leq
{w!\over (w-1)^w }
\leq
{\left( {1\over 2}
\sqrt{\lg (n)}\right) !\over
\left( {1\over 2}
\sqrt{\lg (n)}-1\right) ^{{1\over 2} \sqrt{\lg (n)}} }
,$$
and the right hand side goes to zero as $n\to \infty $.

In case no max-locked interdependent orbit union
$U\subseteq P$
of height
$1$
contains
$\geq \sqrt{\lg (n)} $ elements, we have that
$P$ contains
$\geq \sqrt{\lg (n)} $
interdependent orbit unions
$(U, {\cal N}|U )$
that are
forbidden configurations.
By Lemma \ref{maxlockh1presrank},
if
$(U, {\cal N}|U )$
is max-locked of height $1$, we have
${
|{\rm Aut}_{{\cal N}|U} (U)|
\over
|{\rm End}_{{\cal N}|U} (U)|
}
\leq
{w!\over (w-1)^w }
, $
which, for
$w\geq 3$,
is $\leq {3\over 4}$.
For
$(U, {\cal N}|U )$
max-locked of height $1$ and width $2$, we have
${
|{\rm Aut}_{{\cal N}|U} (U)|
\over
|{\rm End}_{{\cal N}|U} (U)|
}
\leq
{2\over 4}<{3\over 4} $.
For any
forbidden configuration $(U, {\cal N} |U)$
that is not a max-locked interdependent
orbit union of height $1$,
by Lemma \ref{forbconfbound}, we have that
${
|{\rm Aut}_{{\cal N} |U} (U)|
\over
|{\rm End}_{{\cal N} |U} (U)|
}
\leq
{3\over 4} .$
By Proposition \ref{usingEndD},
we obtain
$${ |{\rm Aut} (P)|
\over
|{\rm End} (P)|}
\leq
\left( {3\over 4} \right) ^{\sqrt{\lg (n)} },
$$
and the right hand side goes to zero as $n\to \infty $.

Therefore, the Automorphism Conjecture is true for this
class of ordered sets.
\qed

\begin{lem}
\label{wid12AClem}

There is an $N\in {\mat N} $ such that, for every
coconnected
ordered set $P$
of width $w\leq 11$ with
$n\geq N$ elements,
no nontrivial order-autonomous antichains,
and at most $\lg (n) $
elements
contained in
forbidden configurations, we have that
${|{\rm Aut} (P)|
\over
|{\rm End} (P)|}
\leq 2^{-0.005 n} $.

\end{lem}

{\bf Proof.}
Because the claim is trivial for max-locked ordered sets of width $w$
and because $P$ is coconnected, we can focus on
flexible interdependent orbit unions and on
max-locked subsets of width less than $w$.

By
Lemma \ref{halfdictateslemmaxlock} and
Theorem \ref{allowedconf},
for all
interdependent orbit unions
$(U,{\cal D})$,
that are not forbidden configurations,
there are numbers
$w_1 , \ldots , w_M \in \{ 0, \ldots , \max \{ w-1 , 4\} \} $ such that
$|{\rm Aut } _{{\cal D}} (U)|
\leq
\prod _{j=1} ^M w_j !
$
and such that
$
\sum  _{j=1} ^M w_j
\leq
\left\lfloor {1\over 2}
\left(
|U|-|{\cal D} | \right)
\right\rfloor $.

Let
$Z\subseteq P$ be the set of all points in $P$
that are
contained in (possibly trivial) interdependent orbit
unions that are
not forbidden configurations.
Let $K:=|Z|$
and let $D$ be the number of orbits of $P$
contained in $Z$.
By assumption $n-K\leq \lg (n) $.
By the above,
there are numbers
$w_1 , \ldots , w_M \in \{ 0, \ldots , \max\{ w-1, 4\} \} $ such that
the number of
restrictions of automorphisms of $P$
to
$Z$ is
bounded by
$\prod _{j=1} ^M w_j !
$
and such that
$
\sum  _{j=1} ^M w_j
\leq
\left\lfloor {1\over 2}
\left(
K-D \right)
\right\rfloor $.
The number of
restrictions of automorphisms of $P$
to $P\setminus Z$
is bounded by
$|P\setminus Z|! = (n-K)!\leq \lg (n)!\leq
\lg (n)^{\lg (n)} =2^{\lg (n) \lg (\lg (n))}$.

Moreover, we have
$D\geq {K\over w} $
and, for $w\geq 5$ we obtain
\begin{eqnarray*}
|{\rm Aut} _{{\cal N}_P } (Z)|\leq
\prod _{j=1} ^M w_j !
& \leq &
((w-1)!)^{\left\lceil {1\over w-1} \sum  _{j=1} ^M w_j
\right\rceil }
\\
& \leq &
((w-1)!)^{{1\over w-1} \left\lfloor {1\over 2}
\left(
K-D \right)
\right\rfloor +1}
\\
& \leq &
((w-1)!)^{{1\over w-1} \left\lfloor {1\over 2}
\left(
K-{K\over w} \right)
\right\rfloor +1}
\\
& = &
((w-1)!)^{{1\over 2(w-1)}
\left(
1-{1\over w} \right)
K+1 }
\\
& \leq &
2^{\lg ((w-1)!)
\left(
{1\over 2w}
n+1 \right) }
.
\end{eqnarray*}

For $5\leq w\leq 11$, we have
$\lg ((w-1)!)
{1\over 2w}
\leq
0.991$.
For $w=4$, the
$w-1$ in the above estimate must be replaced with
$4$, and then
$
\lg (4!)
{1\over 2\cdot 4}
<0.574$.
For $w=3$, we can use the crude estimate
$|{\rm Aut} (P)|\leq (3!)^{{n\over 3} +1} =
2^{\lg (3!)\left( {n\over 3} +1\right) } $
and ${\lg (3!)\over 3} <0.87 $.

Hence, in any case,
the number of automorphisms of
$P$ is bounded by
$
2^{0.991n+1}
2^{\lg (n) \lg (\lg (n))}
$,
which, for large enough $n$, is smaller than
$
2^{0.992n}
$.
The claim now follows from
Lemma \ref{2tonforbddwidth}
with $\varepsilon :=0.003$.
\qed

\vspace{.1in}

We are left with the task to consider
ordered sets with nontrivial order-autonomous antichains.
Proposition \ref{onelexiter} below focuses on the
more general lexicographic sum construction.

\begin{define}
\label{lexsumdef}

(See, for example, \cite{HH1,MR}.)
Let $T$ be a nonempty ordered set considered as an index set.
Let
$\{ P_t \} _{t\in T} $ be a family of
pairwise
disjoint nonempty ordered sets that
are all disjoint from $T$. We define
the {\bf lexicographic sum} $L\{ P_t \mid t\in T\} $
{\bf (of the $P_t $ over $T$)}
to be
the union
$\bigcup _{t\in T} P_t $
ordered by letting $p_1 \leq p_2 $ iff
either
\begin{enumerate}
\item
There are distinct
$t_1 , t_2 \in T$ with
$t_1 <t_2 $,
such that
$p_i\in P_{t_i } $,
or
\item
There is a $t\in T$ such that
$p_1 , p_2 \in P_t $ and $p_1 \leq  p_2 $ in $P_t $.
\end{enumerate}
The ordered sets
$P_t $ are the {\bf pieces} of the lexicographic sum and
$T$ is the {\bf index set}.

\end{define}

\begin{prop}
\label{onelexiter}

Let ${\cal C}$ be a class of finite
ordered sets for which the Automorphism
Conjecture holds and let
${\cal L} _{\cal C}$ be the
class of lexicographic sums $L\{ P_t |t\in T\} $ such that all $P_t \in {\cal C}$
and
$T\in {\cal C} $ is indecomposable or a chain or
$T\in {\cal C} $ does not contain any nontrivial order-autonomous antichains and
all $P_t $ are antichains.
Then
the Automorphism
Conjecture holds for
${\cal L} _{\cal C}$.

\end{prop}

{\bf Proof.}
For every $n\in {\mat N}$, let
$$h(n):=
\max _{P\in {\cal C} , |P|\geq n}
{ |{\rm Aut} (P)|\over |{\rm End} (P)| } .
$$
Because
$\lim _{n\to \infty } h(n)=0$,
for $n\geq 2$, we have that
$h(n)<1$.

For every
$L\{ P_t |t\in T\} \in {\cal L} _{\cal C} $,
let
$A(P_t:t\in T):=
\sum _{|P_t|>1} |P_t |$.

Let
$\varepsilon >0$.
Fix
$K\in {\mat N}$ such that,
for all $k\geq K$, we have that
$(h(2))^k
<\varepsilon $ and
$h(k)
<\varepsilon $.
Choose $N\in {\mat N}$ such that
$N\geq K^2 $ and such that, for all
$n\geq N$, we have
$h\left( n-
K^2 \right) \leq \varepsilon ^2 $
and $\left( K^2 \right) !\leq {1\over \sqrt{h\left( n-
K^2 \right) } }$.
Then, for every
$L\{ P_t |t\in T\} \in {\cal L} _{\cal C} $
with $n\geq N$ elements, we obtain the following.

Recall that automorphisms map
order-autonomous sets to order-autonomous sets.
Therefore,
for every automorphism $\Phi \in {\rm Aut} (P)$,
all endomorphisms of all sets $P_s $ translate into
pairwise distinct
order-preserving maps from $P_s $ to $\Phi [P_s ]$.
In case
$A(P_t:t\in T) >K^2 $,
this provides the following estimates.

In case
$A(P_t:t\in T) >K^2 $ and
there is an $s\in T$ with
$|P_s |>K$, we have
\begin{eqnarray*}
{
|{\rm Aut} (L\{ P_t |t\in T\} )|
\over
|{\rm End} (L\{ P_t |t\in T\} )|
}
& \leq &
{
|{\rm Aut} (P_s )|
\over
|{\rm End} (P_s )|
}
\leq h(|P_s | )
<\varepsilon
\end{eqnarray*}

In case
$A(P_t:t\in T) >K^2 $ and, for all
$t\in T$ we have
$|P_t |\leq K$,
let $R$ be the set of
indices $r\in T$ such that
$|P_r |>1$.
Then $|R|\geq K$
and we have
\begin{eqnarray*}
{
|{\rm Aut} (L\{ P_t |t\in T\} )|
\over
|{\rm End} (L\{ P_t |t\in T\} )|
}
& \leq &
\prod _{r\in R}
{
|{\rm Aut} (P_r )|
\over
|{\rm End} (P_r )|
}
\\
& \leq &
(h(2))^K
<\varepsilon
\end{eqnarray*}

This leaves the
case
$A(P_t:t\in T) \leq K^2 $.
Because
$T\in {\cal C} $ is indecomposable or a chain or
$T\in {\cal C} $ does not contain any nontrivial order-autonomous antichains and
all $P_t $ are antichains,
automorphisms map sets $P_t $ to sets $P_t $ and every automorphism
induces a corresponding automorphism on $T$.
Hence we obtain the following.
\begin{eqnarray*}
{
|{\rm Aut} (L\{ P_t |t\in T\} )|
\over
|{\rm End} (L\{ P_t |t\in T\} )|
}
& \leq &
{
|{\rm Aut} (T)|\cdot A(P_t:t\in T)!
\over
|{\rm End} (T )|
}
\\
& \leq &
h\left(
|T|
\right)
A(P_t:t\in T)!
\\
& \leq &
h\left(
|T|
\right)
\left( K^2 \right) !
\\
& \leq &
\sqrt{h\left(
n
-
K^2
\right) }
<\varepsilon
\end{eqnarray*}
\qed

\begin{theorem}
\label{ACw<=12}

The
Automorphism Conjecture is true for ordered sets of
width $w\leq 11$.

\end{theorem}

{\bf Proof.}
By
Proposition \ref{manyinmaxlock}
and Lemma \ref{wid12AClem}, the
Automorphism Conjecture is true for
coconnected
ordered sets of
width $w\leq 11$ without nontrivial order-autonomous antichains.
Now apply Proposition \ref{onelexiter} twice,
using the coconnected
ordered sets of width $\leq 11$ without order-autonomous antichains
and the antichains as the base class ${\cal C}$ in the first application.
The first application leads to a class that contains all
coconnected
ordered sets of
width $w\leq 11$. The second leads to a class
that contains all
ordered sets of
width $w\leq 11$.
\qed

\begin{remark}
\label{betterlb}

{\rm
A natural way to improve the lower bound given by
Theorem 1 in \cite{DRSW} (on page 20 of \cite{DRSW})
is to find a significant number of chains of length $h$ in $P$.
To do so, we could start with an element of rank $h$ and
build chains by
going through every lower cover of
every element in the construction.
However, not every lower cover of an element of rank $k$
has rank $k-1$ and, in the most extreme case,
every element of rank $k$ has exactly one lower cover of rank $k-1$.
Thus, there is no easy way to improve the
lower bound for the number of endomorphisms.

Moreover, even the most optimistic estimate, for the above approach,
which is $w^{n\over w} =2^{{\lg (w)\over w} n} $,
indicates that
the additional factor obtained would
not lead to the Automorphism Conjecture for all ordered sets
of any given fixed width.
With the arguments given here, this improvement could
possibly
establish the
Automorphism Conjecture for ordered sets of width $\leq 15$.
\qex

}

\end{remark}

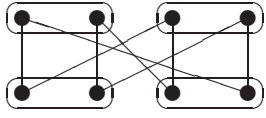
\begin{figure}

\centerline{
\unitlength 1mm 
\linethickness{0.4pt}
\ifx\plotpoint\undefined\newsavebox{\plotpoint}\fi 
\begin{picture}(36,16)(0,0)
\put(5,5){\circle*{2}}
\put(15,5){\circle*{2}}
\put(25,5){\circle*{2}}
\put(35,5){\circle*{2}}
\put(5,15){\circle*{2}}
\put(15,15){\circle*{2}}
\put(25,15){\circle*{2}}
\put(35,15){\circle*{2}}
\put(5,15){\line(0,-1){10}}
\put(15,15){\line(0,-1){10}}
\put(25,15){\line(0,-1){10}}
\put(35,15){\line(0,-1){10}}
\put(5,5){\line(2,1){20}}
\put(15,5){\line(2,1){20}}
\put(35,5){\line(-3,1){30}}
\put(25,5){\line(-1,1){10}}
\put(10,15){\oval(14,4)[]}
\put(30,15){\oval(14,4)[]}
\put(10,5){\oval(14,4)[]}
\put(30,5){\oval(14,4)[]}
\end{picture}
}

\caption{An $8$-crown $C_8 $ with a dictated orbit structure ${\cal D}$
indicated by ovals.
For $(C_8 , {\cal D})$, we have
$
{\rm Aut}_{{\cal D}} (C_8 )
=
{\rm End}_{{\cal D}} (C_8 )
$
.
}
\label{8crownlocked}

\end{figure}

\begin{remark}
\label{EndDlimit}

{\rm
Figure \ref{8crownlocked}
shows an
interdependent orbit union
$(U, {\cal D})$
such that
$
{\rm Aut}_{{\cal D}} (U )
=
{\rm End}_{{\cal D}} (U )
$.
Indeed,
any ${\cal D}$-endomorphism that collapses an orbit
$D\in {\cal D}$ to a point must, following the orbits in a cyclic order,
collapse each next orbit to the unique adjacent point, which, upon
completing the 4-cycle, leads to a contradiction,
because an order relation is not preserved.

Therefore, although Proposition \ref{usingEndD}
is very useful here and
should have significant potential for future use,
it is limited by the fact that there are
situations in which the
factor
${
|{\rm Aut}_{{\cal N}|U} (U)|
\over
|{\rm End}_{{\cal N}|U} (U)|
}
$
devolves into the trivial factor $1$.
\qex
}

\end{remark}

\section{Conclusion}
\label{conclusec}

Despite the considerable technical details,
the methods and results presented here
are a possible roadmap to proving the Automorphism Conjecture.
Proposition \ref{getallfromiou} and Theorem \ref{pruneorbit6}
show how the automorphism group of an ordered set naturally decomposes
into certain subgroups, which
are directly manifested in the
combinatorial structure of the ordered set.

Call a subset $S$ of an ordered set $P$
{\bf pseudo-order-autonomous} iff
there is a dictated orbit structure ${\cal D}$ for $S$
such that
$(S, {\cal D} )$
is an interdependent orbit union and
every
${\cal D} $-orbit
$D$ of $S$ is order-autonomous in $(P\setminus S)\cup D$.
Clearly, by part \ref{unionplacement0} of
Proposition \ref{unionplacement}, this idea is a generalization of
interdependent orbit unions.
Similar to Lemma \ref{maxlockh1presrank},
we can prove that any max-locked
interdependent orbit union $(M, {\cal D} )$
satisfies
$|{\rm End}_{\cal D} (M)|\geq (w-1)^w $.
Hence,
an abundance of
pseudo-order-autonomous subsets that are
max-locked interdependent orbit unions
up to a fixed height
can be dispensed with
an argument like in
Proposition \ref{manyinmaxlock}.
Moreover, Lemma \ref{maxlockeddescr} guarantees that
``tall" max-locked interdependent orbit unions
$(M, {\cal D} )$ use a disproportionate number of
points to
merely generate a factor $w(M)!$ for $|{\rm Aut} (P)|$.

This leaves us with nontrivial flexible interdependent
orbit unions, and their
further analysis looks promising.
Order-autonomous antichains contained
in dictated orbits can be addressed with
a single application of Proposition \ref{onelexiter}.
By
Theorem \ref{pruneorbit6},
for $D_m $ not being a cutvertex of
${\cal O} (U, {\cal D} )$,
if
$|{\rm Aut } _{{\cal D}_m} (U_m )|\leq 2^{0.6 |U_m |} $
and
$|{\rm Aut } _{{\cal D}_Q } (Q)| \leq
2^{0.6 \left( |Q |-|{\cal D}_Q |+1\right) } $,
we obtain
$|{\rm Aut } _{{\cal D}} (U)|
\leq
2^{0.6 |U_m |}
2^{0.6 \left( |Q |-|{\cal D}_Q |+1\right) }
=
2^{0.6 \left( |U_m |+|Q |-|{\cal D}_Q |+1\right) }
=
2^{0.6 |U |}
$.
Thus, for ordered sets whose interdependent orbit unions are
constructed with stages that satisfy this upper bound,
the ratio of automorphisms to endomorphisms
converges to zero with
an exponential upper bound.
Because, for $n\leq 17$, we have
$\left( \left\lfloor {1\over 2} (n-2)\right\rfloor -1\right) !
\leq 2^{0.6 n} $, the exponential bound holds for all
flexible interdependent orbit unions with
17 or fewer points, where the bound is
easily verified directly for the
flexible forbidden configurations.

The first step to prove the above mentioned
exponential
bound would be
to prove a bound of
$\leq 2^{0.6 |U|}$ automorphisms
for interdependent orbit unions $(U,{\cal D} )$ with $2$ orbits.
The author's limited investigation,
for small $n$, of
ordered sets $U$ of height $1$ with $n$ minimal and $n$ maximal elements
such that there is an automorphism that is a single cycle on
the minimal elements and
that is a single cycle on
the maximal elements,
has shown that these sets
can only have
more than $2^{0.6 |U|}$
automorphisms when they contain
pseudo-order-autonomous
max-locked
ordered sets $M$
of height $1$.
Similarly, \cite{ConVer} shows that,
up to 47 vertices, all vertex-transitive and
edge-transitive graphs $G=(V,E)$
that are neither complete, nor complete bipartite, have
fewer than
$2^{{1\over 2} |V|}
<2^{0.6 (|V|+|E|)}
$
automorphisms.
Hence, the bound $2^{0.6 |U|}$
from
above may be more than simply ``conveniently chosen."

By
Lemma \ref{pruneorbit5},
after using duality as needed and deleting
comparabilities between sets $A_i ^j $, for $D_m $
not a cutvertex of ${\cal O} (U,{\cal D} )$,
the analysis of interdependent orbit unions
$(Q , {\cal D}_Q )$
can be reduced to a thorough analysis of
ordered sets
$Q$ of height $1$
with a
tight dictated orbit structure
${\cal D} =\{ B, T_1 , \ldots , T_\ell \} $
such that, for $j=1, \ldots , \ell $, we have
$B\upharpoonleft _{{\cal D}_Q} T_j $.
Because of the promising behavior for ordered sets of height $1$ so far,
the following conjecture appears natural.

\begin{conj}
\label{expboundconj}

For all
indecomposable ordered sets $U$
with $n$ elements,
we have that
$|{\rm Aut} _{{\cal N}  } (U)|\leq 2^{0.6 |U|}$, or, there are
at least
$O\left( {n\over \lg (n)} \right) $ elements in
non-singleton antichains $|A|$ on which the automorphism group induces all
$|A|!$ permutations.

\end{conj}

Along different lines,
the set
${\rm Aut } _{{\cal D}_Q } ^U (Q) $
is, by
Theorem \ref{pruneorbit6},
a normal subgroup of ${\rm Aut} _{\cal D} (U)$
and
the factor group
${\rm Aut} _{\cal D} (U)/{\rm Aut} _{{\cal D} _Q } ^U (Q)$
is isomorphic to
a subgroup of
${\rm Aut } _{{\cal D}_n} (U_n )$.
It would be interesting to explore if
one can efficiently
obtain interdependent orbit unions $(U,{\cal D})$
such that ${\rm Aut} _{\cal D} (U)$
is a prescribed group, and
such that
fewer points are used than in the so far most efficient
construction in
\cite{BMgroupAutP}.
(Note that, for a $D_n $ such that ${\cal O} (U,{\cal D})-D_n $ is
disconnected,
after proper relabeling,
the results about $(U_n , {\cal D} _n )$ apply to all non-singleton
components of ${\cal O} (U,{\cal D})-D_n $.)
Here a thorough analysis
of the situation in which all $|A_i ^j |=1$ would be needed.

\vspace{.1in}

{\bf Acknowledgement.}
The author very much thanks Frank a Campo for a thorough
reading of an early version of this paper, for the detection of
mistakes, many helpful suggestions, and revisions of some proofs,
in particular
Lemma \ref{allbutoneac} and
a significant shortening of the proof
of Lemma \ref{maxlock2lev}.

\end{document}